\documentclass{article}

\usepackage{hyperref}
\usepackage{cite}
\usepackage{enumitem}
\usepackage[ruled,longend]{algorithm2e}
\usepackage{todonotes}

\usepackage{fancyhdr}
\usepackage{caption}
\usepackage{subcaption}
\usepackage{float}
\usepackage{amsmath,amssymb,amsthm,graphicx,url}

\theoremstyle{plain}
\newtheorem{lemma}{Lemma}[section]
\newtheorem{proposition}[lemma]{\textbf{Proposition}}
\newtheorem{fact}[lemma]{\textbf{Fact}}
\newtheorem{theorem}[lemma]{\textbf{Theorem}}
\newtheorem{corollary}[lemma]{\textbf{Corollary}}

\theoremstyle{definition}
\newtheorem{definition}[lemma]{\textbf{Definition}}
\newtheorem{example}[lemma]{\textbf{Example}}

\newtheorem{remark}[lemma]{Remark}

\title{A tree-based algorithm for the integration of monomials 
in the Chow ring of the moduli space of stable marked curves of genus zero}

\author{Jiayue Qi \thanks{Doctoral Program ``Computational Mathematics'' W1214, Johannes Kepler 
University Linz.}\thanks{Research Institute for Symbolic Computation, Johannes Kepler University, Linz, Austria.}}
\date{}
\begin{document}
\maketitle

\abstract{The Chow ring of the moduli space of marked rational curves is generated by Keel's 
divisor classes. The top graded part of this Chow ring is isomorphic to the integers, generated by 
the class of a single point. In this paper, we give an algorithm for computing the intersection degree of 
tuples of Keel's divisor classes. This computation is a concrete but complicated 
algorithmic question in the field. Also, we give a simple complexity argument for the algorithm.
Additionally, we introduce three identities on multinomial coefficients, as well as proofs for them.}

\section{Introduction}\label{sec:introduction}

The moduli space of stable $n$-pointed curves of genus zero, denoted by $\overline{\mathcal{M}}_{0,n}$, 
is a renowned object in modern intersection theory; for example,
it is the base for the definition of Gromov-Witten invariants~\cite{GW}. It is a smooth irreducible projective 
compactification of stable $n$-pointed genus-zero curves. This moduli space was first constructed 
by Knudsen and Mumford in their series of papers~\cite{km1},~\cite{km2} and~\cite{km3}; they
introduced also other constructions for moduli spaces for curves with genus bigger than zero.

Chow rings are essential in intersection theory.
%, to indicate the intersection numbers of subvarieties
 Let~$X$ be a projective variety 
of dimension~$k$. The Chow ring $A^{\bullet}(X)$ of~$X$ is the direct 
sum of $k+1$ abelian groups~$A^r(X)$, each of which is composed of all the cycles 
(formal sums for subvarieties modulo algebraic equivalence)
of codimension~$r$; by convention, we set $A^r(X)=\{0\}$ when $r>k$. 
%Particularly, $A^k(X)$ is the Chow group of cycles of dimension zero and is homomorphic to the integer 
%additive group, via the degree isomorphism.

In this paper, we work in the Chow ring of the moduli space $\overline{\mathcal{M}}_{0,n}$ of stable $n$-pointed curves of genus zero
and we denote it by $A^{\bullet}(\overline{\mathcal{M}}_{0,n})$. 
Since $\overline{\mathcal{M}}_{0,n}$ is of dimension $n-3$, we know that 
$A^{\bullet}(\overline{\mathcal{M}}_{0,n})=\bigoplus^{n-3}_{r=0}{A^r(\overline{\mathcal{M}}_{0,n})}$,
where $A^r(\overline{\mathcal{M}}_{0,n})$ denotes the 
Chow group of codimension~$r$. There is an isomorphism between $A^{n-3}(\overline{\mathcal{M}}_{0,n})$ and the integer additive group~$\mathbb{Z}$.
We have $A^r(\overline{\mathcal{M}}_{0,n})=~\{0\}$ when $r>n-3$, and $A^{n-3}(\overline{\mathcal{M}}_{0,n})\cong \mathbb{Z}$ --- we use the 
symbol~$\int$
to denote this degree isomorphism, following standard convention. The integer under this isomorphism is 
called the {\em integral value} or {\em value} of the given element in $A^{n-3}(\overline{\mathcal{M}}_{0,n})$.  
A set of generators of the group $A^1(\overline{\mathcal{M}}_{0,n})$ was given in Keel's paper~\cite{keel}, where each generator 
is indexed by a bi-partition of $\{1,\ldots,n\}$ and this set also generates the whole ring. They are the classes of the boundary divisors 
of $\overline{\mathcal{M}}_{0,n}$. 

%Computations in $A^{\bullet}(n)$ soon become rather hard as $n$ increases: numerical experiments
%indicate that the logarithm of the number of linearly independent elements in the ambient ring 
%grows quadratically. And the rank of the Chow group $A^1(\overline{\mathcal{M}}_{0,n})$ of codimension one which 
%actually generates the whole ring via multiplication, grows very fast, since it is exponential 
%in $n$. 
Note that the Chow ring $A^{\bullet}(\overline{\mathcal{M}}_{0,n})$ can be described as a quotient of a polynomial ring, hence
we can talk about monomials in this ring.
We will introduce an algorithm --- called {\em the forest algorithm} ---
for computing the integral value of a monomial in the generators of $A^{\bullet}(\overline{\mathcal{M}}_{0,n})$. 
This problem showed up during the study of counting the realization
of Laman graphs (minimally-rigid graphs) on a sphere~\cite{laman}, when we wanted to improve the algorithm
given in~\cite{laman}. 
%With the help of this algorithm, we invented another algorithm for the same goal. 
%However, by efficiency 
%it does not seem faster than the one provided in~\cite{laman}. 
We see this problem fundamental, 
standing on its own; we find the algorithm elegant and concise,
may as well be helpful for other similar or even further-away problems. Therefore, we formulate
it on its own.
We consider the situation when this monomial is of degree $n-3$, otherwise
we define its value to be zero. The input is $n-3$ such generators, hence is linear in~$n$ and the output 
is an integer. Our algorithm is quadratic in~$n$.

 Despite the existing work~\cite{KM94, Kau98}, the combinatorics
involved is complicated enough to make computations in $A^{\bullet}(\overline{\mathcal{M}}_{0,n})$ difficult. 
Although the situation for genus zero is relatively well understood, in practice, large computations can become intractable. 
In~\cite{Tav17}, Tavakol provides another set of generators for $A^{\bullet}(\overline{\mathcal{M}}_{0,n})$.
In this paper~\cite{keel}, Keel describes a quadratic relation between the generators of the Chow ring of 
$\overline{\mathcal{M}}_{0,n}$ that determines when a product is zero; also, he describes a linear relation between 
the generators. 
The property of Keel's quadratic relation naturally indicates the first step of our algorithm. We check if any of the two factors of the input
monomial fulfills this relation: if yes, return zero; otherwise, we consider an equivalent characterization 
of the given monomial, in a specific tree --- {\em loaded tree}. 
This characterization dates back to~\cite{Buneman1971}, where they character a phylogenetic tree using 
the split representation. The loaded trees we consider in this paper differ from phylogenetic trees only
in the sense that they allow edges to have more-than-one multiplicity. However, the conversion algorithm we use to transfer a monomial to 
a loaded tree comes from~\cite[Section 2.2, Definition 23]{QS:representation}. 
This first part is quadratic in~$n$ in the worst case.
In the second step, we transfer this tree 
via three steps to a forest. Next, we compute the integral value of the given monomial directly from the obtained forest.
The second part is linear with respect to~$n$.

To prove the correctness of the forest algorithm, references indicates that there may be another way, at least for the base case proof,
namely one can use Lemma~3 of~\cite{Cav16}
to re-express the input monomial, and then use Lemma~1.5.1 of~\cite{Koc01} to evaluate the resulting $\psi$-monomial. 
Another potential method can be derived from~\cite{Fab99}.
We will not go into details upon these approaches in this paper. 
We give a rather direct algebraic proof for the general case of the forest algorithm, in Section~\ref{sec:from_algebra_to_geometry}.
For the base case proof of the correctness of the forest algorithm, we introduce an equivalent characterization of the {\em linear reduction} 
on the monomial, via loaded trees.
Linear reduction is a reduction step on the given monomial, using the linear relation to substitute a factor and 
then eliminate all items that are zero because of the quadratic relation; more details can be seen in Section~\ref{sec:algebraic_linear_reduction}.
The equivalent characterization mentioned above is manifested mainly via an operation on the loaded trees
called {\em vertex splitting}, for more details see Section~\ref{subsec:vertex_splitting}.
This characterization can also lead to another algorithm for our focused question, however with rather bad time complexity.
We praise the idea behind since it gives another way of understanding the linear reduction, with much insight. 
Another reason is that, this proof delivers to us three new identities on multinomial coefficient, which we find elegant and meaningful.
In the paper, we also give the proof for those identities, see Section~\ref{sec:sun_like_trees},~\ref{subsec:identity_equivalence}, 
and~\ref{sec:proof_identity}.
Our proof addresses the problem in a graphical and more combinatorial way. 
%gives a better picture in the combinatorial 
%aspect, and may confluent with the just-mentioned method and reveal in a concise and clever way of the combinatorial 
%structure of the algorithm. 

%A quadratic relation between the generators was introduced in Keel's paper~\cite{keel}. This relation is called 
%{\em Keel's quadratic relation} in our paper. With the help of this 
%relation, we know that if any two factors of the given monomial fulfills the relation, then the whole monomial has value zero. 

%The tree characterization of the given monomial is inspired by 
%\cite[Section~2.2.]{QS:representation} and was first introduced in~\cite{ACM_communication}. 
% Note that the same problem was considered in~\cite{vertex_splitting}, where an 
% equivalent characterization for the algebraic reductions (in the ambient ring) for the input monomial is given, as some operations
% on the tree representation. In fact, we were motivated by that characterization, then
%we started to try-out many examples using the algorithm provided in~\cite{vertex_splitting}. Eventually and excitingly, 
%we discovered our algorithm which as an algorithm, is much more efficient and neat, comparing to the one given
%in~\cite{vertex_splitting}. 
%However, as a characterization,~\cite{vertex_splitting} wins the match for sure.
Our paper can be considered
as a proper extension for~\cite{ACM_communication}. We consider the same problem, describe the same algorithm, 
as~\cite{ACM_communication}; but with all the
detailed proofs provided. 

\section{Preliminaries}
Since the main problem this paper focus at is exactly the same with that of~\cite{ACM_communication},
the preliminaries will be very much similar. However, for completeness, we introduce everything from scratch.

Let~$n\in \mathbb{N}$, $n\geq 3$, and let~$N$ be any set with cardinality~$n$. 
Note that in most situations (in this paper) we take~$N:=\{1,\ldots,n\}$ as default.
Denote by~$\overline{\mathcal{M}}_{0,n}$ the moduli 
space of stable~$n$-pointed curves of genus zero. 
\begin{definition}
A bi-partition~$\{I,J\}$ of~$N$ where the cardinalities of~$I$ and~$J$ are both at 
least~$2$ is called a {\bf cut} (w.r.t.~$N$); we call~$I$ and~$J$ {\bf parts} of this cut.
\end{definition}
There is a hypersurface (boundary divisor)~$D_{I,J}\subset \overline{\mathcal{M}}_{0,n}$ for each cut~$\{I,J\}$ and its class in the Chow ring of 
$\overline{\mathcal{M}}_{0,n}$
is denoted by~$\delta_{I,J}$ --- note that~$D_{I,J}=D_{J,I}$, and as well~$\delta_{I,J}=\delta_{J,I}$.    
The Chow ring of~$\overline{\mathcal{M}}_{0,n}$ is a graded ring and we denote it by~$A^{\bullet}(\overline{\mathcal{M}}_{0,n})$. We have 
$$A^{\bullet}(\overline{\mathcal{M}}_{0,n})=\bigoplus^{n-3}_{r=0}{A^r(\overline{\mathcal{M}}_{0,n})}.$$
\begin{definition}
The homogeneous components of the Chow ring~$A^{\bullet}(\overline{\mathcal{M}}_{0,n})$
are defined as the Chow groups (of~$\overline{\mathcal{M}}_{0,n}$);~$A^r(\overline{\mathcal{M}}_{0,n})$ 
is the {\bf Chow group of codimension~$r$}.
\end{definition}
It is known that $A^r(\overline{\mathcal{M}}_{0,n}):=\{0\}$ for~$r > n-3$, 
and $A^{n-3}(\overline{\mathcal{M}}_{0,n})\cong \mathbb{Z}$. There is a canonical isomorphism (the degree isomorphism)
$$\int: A^{n-3}(\overline{\mathcal{M}}_{0,n}) \to \mathbb{Z}$$ that takes the class of a point to~$1$.
We can extend this map so that it is defined on the whole ring by setting the value of 
all other elements to be zero:
$\int:A^{\bullet}(\overline{\mathcal{M}}_{0,n})\to \mathbb{Z}$, $\int(M)=0$ if~$M\notin A^{n-3}(\overline{\mathcal{M}}_{0,n})$.

\begin{definition}
The set
$G_N:=\{\delta_{I,J}\, \mid \, \{I,J\} \text{ is a cut w.r.t. }N\}$ generates the group $A^1(\overline{\mathcal{M}}_{0,n})$, and also the 
whole ring~$A^{\bullet}(\overline{\mathcal{M}}_{0,n})$~\cite{keel}. We call each such generator a {\bf Keel's generator} or a {\bf Keel's factor}.
\end{definition}
Then~$\prod_{i=1}^{n-3}{\delta_{I_i,J_i}}$ is an element 
in~$A^{n-3}(\overline{\mathcal{M}}_{0,n})$. We define the {\bf value} of~$M:= \prod_{i=1}^{n-3}{\delta_{I_i,J_i}}$ to be  
$\int(\prod_{i=1}^{n-3}{\delta_{I_i,J_i}})$. 
The problem we deal with in this paper
is \textbf{\em to compute the value of a given monomial
$$M= \prod_{i=1}^{n-3}{\delta_{I_i,J_i}}.$$}

Keel introduced a quadratic relation between Keel's factors in~\cite{keel}; 
we call it {\em Keel's quadratic relation}. 
% Denote by $\epsilon_{ij|kl}:=\sum_{i,j\in I,k,l\in J}{\delta_{I,J}}$. Then the following 
% equalities which we call {\bf Keel's linear relation} (\cite[Section 4, Theorem 1.(2)]{keel}) hold:  
% $\epsilon_{ij|kl}=\epsilon_{il|kj}=\epsilon_{ik|jl}$. 
\begin{definition}
We say that two generators
~$\delta_{I_1,J_1}, \delta_{I_2,J_2}$ fulfill {\bf Keel's quadratic relation} (\cite[Section 4, Theorem~1.(3)]{keel}) if the following 
 four conditions hold:
 $$I_1\cap I_2\neq \emptyset; I_1\cap J_2\neq \emptyset;
 J_1\cap I_2\neq \emptyset; J_1\cap J_2\neq \emptyset.$$
 In this case, we have~$\delta_{I_1,J_1}\cdot \delta_{I_2,J_2}=0$, the two corresponding boundary divisors have empty intersection.
 \end{definition}
For example, when~$n=5$,
$\delta_{12,345}\cdot \delta_{14,235}=0$ since these two factors fulfill the Keel's quadratic relation. 
Note that we use abbreviated notations for the index
of the generators, for instance~$\delta_{12,345}$ represents~$\delta_{\{1,2\},\{3,4,5\}}$. We
will use this abbreviation also in the later context.
%Motivated by this quadratic relation, we realize that whenever two factors of the given monomial
%fulfill this relation, the value of the monomial is zero. 
Hence the first step of our algorithm 
is to check whether there are two factors fulfilling this relation --- if yes, return zero. There are in total~$n-3$ input 
factors, so we need to check~${n-3 \choose 2}$ many pairs in the worst case. The checking for each
pair of generators is linear in~$n$. Hence the algorithm in this step is quadratic in~$n$ in the worst case.

\begin{definition}
We call those monomials of which no two factors fulfill the Keel's quadratic relation {\bf tree monomials},
since there exists a one-to-one correspondence between these monomials and a specific type of trees which 
we call {\em loaded trees} (see Definition~\ref{def:loaded_tree}).
\end{definition}
Note that in the case of a tree monomial, the set-theoretical intersection
$\bigcap_{\delta_{I,J}\in M}{D_{I,J}}$
is non-empty (\cite{Gia16}), and is a codimension~$k$ boundary stratum.
Then, how should we compute the value of a tree monomial? The following theorem 
indicates the first thing to check, when we have a tree monomial at hand. 
This follows from the boundary of~$\mathcal{M}_{0,n}$ 
being simple normal crossings.
%We postpone the proof of it to Section~\ref{sec:proof_sun_like}.
\begin{theorem}\label{thm:clever}
 If~$M\in A^{n-3}(\overline{\mathcal{M}}_{0,n})$ is a tree monomial that is the product of~$n-3$ distinct factors, 
 then~$\int(M)=1$. In this case, the corresponding boundary divisors intersect transversally
at a single point. We call this type of tree monomials {\bf clever monomials}.
\end{theorem}
So the first step of our algorithm is: going through all pairs, checking whether they fulfill Keel's 
quadratic relation and whether they are distinct. It is quadratic in~$n$.

%Actually, we can combine the first two steps as one step, serving as the first part of our algorithm. This is because
%going through all the pairs of generators once is sufficient: we can check whether the pair fulfills the 
%Keel's quadratic relation or not and whether the generators in the pair are distinct at the same time. Hence the first part of our algorithm 
%is polynomial with respect to $n$.

\subsection{Loaded trees}\label{sec:loaded_trees}

In this section, we explain the one-to-one correspondence between tree monomials and loaded trees.
%and we show with examples how to transfer between these two different 
%representations. 

\begin{definition}[\cite{ACM_communication}, Definition 0.1.]\label{def:loaded_tree}
 A {\bf loaded tree with~$n$ labels and~$k$ fringes} is a tree~$T=(V,E)$ together with a 
 labeling function~$h: V\to 2^N$ where~$N$ is any set of cardinality~$n$ and an edge multiplicity function~$m: E \to \mathbb{N}^+$ such that the
 following three conditions hold:
 \begin{enumerate}
  \item $\{h(v)\}_{v\in V, h(v)\neq \emptyset}$ form a partition of~$N$; elements in~$N$ are called the {\bf labels}
  of~$T$.
 \item $\sum_{e\in E}{m(e)}=k$. Note that some edges can have a multiplicity more than one ---  those are called {\bf multi-edges}.
 \item For every~$v\in V$, $\deg(v)+|h(v)|\geq 3$, where~$\deg(v)$ is the degree of vertex~$v$ --- note that here 
 a multi-edge only contributes one to the degree 
 of its incident vertices, as a single edge does.
\end{enumerate}
We say that the loaded trees~$T_1$ and~$T_2$ are {\bf of the same type} if and only if they have the same set of labelings and 
the same number of fringes.
\end{definition}

We define the {\em monomial of a loaded tree} as follows. 
\begin{definition}
Denote by~$G_N:=\{\delta_{I,J}\mid \{I,J\}\text{ is a cut w.r.t. }N\}$ the set of Keel's factors (w.r.t.~$N$). 
Let~$\Psi: E\to G_N$ be the function that assigns to each edge its corresponding Keel's factor, where~$E$ is the 
edge set of the loaded tree~$LT=(V,E,h,m)$ with labeling set~$N$. Then
$$M_T:=\prod_{e\in E}{\Psi(e)^{m(e)}}$$
is called the {\bf the monomial~$LT$}. 
\end{definition}

\begin{definition}
 A loaded tree is {\bf proper} if its number of labels is its number of fringes plus three.
 Monomials of proper loaded trees are {\bf proper monomials}.
\end{definition}

\begin{figure}
\centering
\begin{subfigure}{.4\textwidth}
  \centering
  \includegraphics[width=0.9\linewidth]{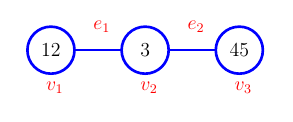}
  \caption{}
\end{subfigure}%
\begin{subfigure}{.4\textwidth}
  \centering
  \includegraphics[width=0.9\linewidth]{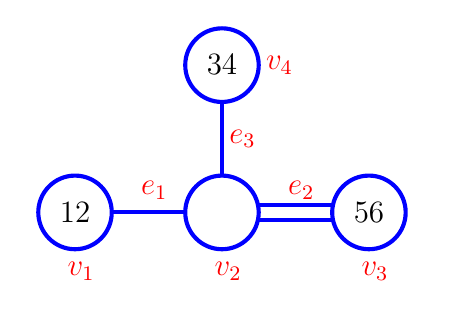}
  \caption{}
\end{subfigure}
\caption{On the left is a loaded tree with~$5$ labels and~$2$ fringes; its monomial is~$\delta_{12, 345}\cdot \delta_{123,45}$.
On the right is a loaded tree with~$6$ labels and~$3$ fringes; its monomial is 
$\delta_{12,3456}\cdot \delta_{34,1256}\cdot \delta^2_{56,1234}$. Note that we use an abbreviated notation
for the labeling set of vertices, for instance~$12$ refers to~$\{1,2\}$, for a neater picture.}
\label{fig:loaded_trees}
\end{figure}

 We see that a loaded tree uniquely determines its monomial. 
We see two examples of loaded trees and their monomials in Figure~\ref{fig:loaded_trees}. 
Note that in the example, we use an abbreviated notation for 
the labeling set of vertices shown on the picture. We will keep using it in the later context, for neater pictures.

The following result tells us the existence of a one-to-one correspondence between tree monomials and loaded trees. 
 \begin{lemma}\label{lem:correspondence}
  There is a one-to-one correspondence between tree monomials~$M=\prod^k_{i=1}{\delta_i}$, where~$\delta_i\in G_N$ 
  for all~$1\leq i\leq k$, and loaded trees with labeling set~$N$ and~$k$ fringes. 
 \end{lemma}
\begin{remark}
 Lemma~\ref{lem:correspondence} follows immediately 
 from~\cite[Section 2.2]{QS:representation}. 
 Usually we denote by~$T_M$ the loaded tree of monomial~$M$. 
 The corresponding loaded tree of a clever monomial is actually the dual tree of the intersection point of the 
 divisors in moduli space~$\overline{\mathcal{M}}_{0,n}$.
 \begin{definition}
 We call the corresponding loaded tree for clever monomials {\bf clever trees}.
 \end{definition}
 The corresponding loaded tree of a non-clever monomial --- when each edge is a single edge --- is the dual tree of the boundary stratum 
  of the intersection of corresponding boundary divisors. When the given monomial is non-clever, 
  the corresponding intersection has excess dimension; this is the case our main algorithm deals with. 
 \end{remark}

It is trivial to obtain the monomial of a loaded tree, while the other direction not. The algorithm for this harder direction
is described in~\cite{ACM_communication}. Note that it specifies the ambient group of the input monomial, but
the same algorithm also works for a monomial of any other degree. The idea of this algorithm comes from Section 2.2 of
the paper~\cite{QS:representation}. 
%We will not go into details of this algorithm.
For completeness, we illustrate this algorithm here, see 
Algorithm~\ref{alg:monomial_to_tree}.

 \begin{algorithm}
\thispagestyle{empty}
\caption{monomial to tree (Algorithm 1 in~\cite{ACM_communication})}
\label{alg:monomial_to_tree}
\SetKwInOut{Input}{input}
\SetKwInOut{Output}{output}

\Input{a tree monomial~$M$ in~$A^{n-3}(\overline{\mathcal{M}}_{0,n})$}
\Output{a loaded tree with~$n$ labels and~$n-3$ fringes}
\texttt{\\}
$C \gets$ collection of any cut that corresponds to some factor of~$M$\;
$P \gets$ collection of all the parts of cuts in~$C$\;
$c \gets$ any element~$c=\{I,J\}\in C$\;

\For{each element~$p\in P\setminus \{I,J\}$}
   {\If{$p\subset I$ or~$p\subset J$}
       {$c:=c\cup \{p\}$}}
$H \gets$ the Hasse diagram of elements in~$c$ with respect to set containment order\;
Consider~$H$ as a graph~$(V,E)$\;
\For{each vertex~$V$ of~$H$}
    {Define labeling set~$h(V)$ as its corresponding element in~$c$\;
    Update the labeling set:~$h(V):=h(V)\setminus \bigcup_{V_1<V \text{ in } H}h(V_1)$}

$E:=E\cup \{\{I,J\}\}$\;
Attach this labeling function~$h$ to~$H$\;
Set the multiplicity function value~$m(e)$ for each edge~$e$ 
as the power of its corresponding factor in~$M$\;

\Return{$H=(V,E,h,m)$}         
  
\end{algorithm}

Let us see an example on constructing the corresponding loaded tree of a given monomial, using Algorithm~\ref{alg:monomial_to_tree},
so as to have an intuitive comprehension.
\begin{example}\label{eg:monomial_to_tree}
 Consider the tree monomial 
 $$\delta^3_{123,456789}\cdot \delta_{12345,6789}\cdot \delta_{1234589,67}\cdot \delta_{1234567,89}.$$ 
 Obviously we have the labeling set $N:=\{1,2,3,4,5,6,7,8,9\}$.
 We collect the parts in set 
 $$P:= \{\{1,2,3\},\{4,5,6,7,8,9\},\{1,2,3,4,5\},\{6,7,8,9\},$$ 
 $$\{1,2,3,4,5,8,9\},\{6,7\},\{1,2,3,4,5,6,7\},\{8,9\}\},$$
 and we pick any cut $c=\{\{1,2,3,4,5\},\{6,7,8,9\}\}$
 from the set of cuts. After collecting all parts of the cuts in the subscripts of factors in the given monomial which are either contained 
 in $\{1,2,3,4,5\}$ or $\{6,7,8,9\}$, 
 we obtain $c=\{\{1,2,3,4,5\},\{6,7,8,9\},\{1,2,3\},\{6,7\},\{8,9\}\}$, then we construct the corresponding Hasse diagram for~$c$, 
 see Figure~\ref{fig:hasse}. The output loaded tree~$T_M$ is shown in Figure~\ref{fig:from_hasse}. 
 It is easy to see that if we go back from the tree constructing monomial, we again obtain~$M$.
 
  \begin{figure}[H]
\centering
\begin{subfigure}{.45\textwidth}
  \centering
  \includegraphics[width=0.8\linewidth]{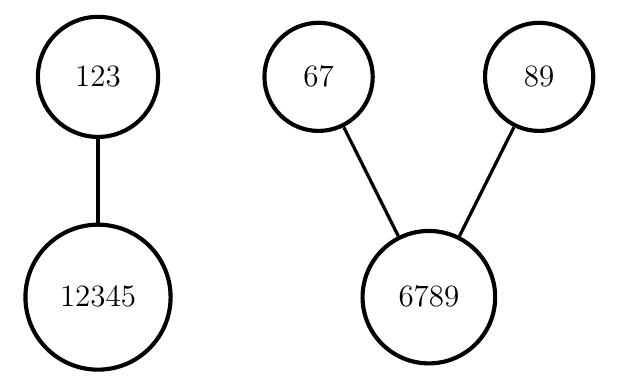}
  \caption{}
  \label{fig:hasse}
\end{subfigure}%
\begin{subfigure}{.45\textwidth}
  \centering
  \includegraphics[width=0.9\linewidth]{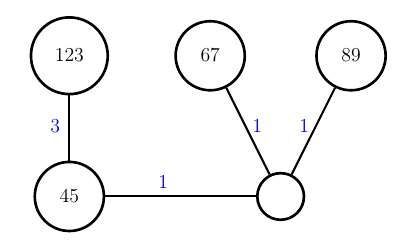}
  \caption{}
  \label{fig:from_hasse}
\end{subfigure}
\caption{On the left is the Hasse diagram of set $\{\{1,2,3,4,5\},\{6,7,8,9\},$ 
  $\{1,2,3\},\{6,7\},\{8,9\}\}$ with respect to set containment order. On the right is 
  the corresponding loaded tree of monomial 
$\delta^3_{123,456789}\cdot \delta_{12345,6789}\cdot \delta_{1234589,67}\cdot \delta_{1234567,89}$. 
Multiplicity function values are written in blue.}
\end{figure}
 \end{example}
 
 \begin{definition}
 If a loaded tree has no fringes, then its monomial has no factors; we
call such a monomial an {\bf empty monomial}. 
\end{definition}
We extend this one-to-one correspondence a little by including 
a single tree and a single monomial: the loaded tree with~$3$ labels and no fringes corresponds to
the empty monomial; the loaded tree with~$k$ labels and no fringes has no corresponding monomial if~$k\neq 3$.
With this extension proposed, we can now define the {\em value of a loaded tree}.
\begin{definition}
The {\bf value of a loaded tree} is the value of its corresponding monomial.
\end{definition}
 Hence our goal can be expressed in other words now: {\bf compute the value of a loaded tree with~$n$ labels and 
$n-3$ fringes, where~$n\geq 3$.} Now we see that by the correspondence, improper
loaded trees have value zero. 
Later we will see, the extension above is done so as to guarantee that the loaded tree has the same value 
with its monomial, which also ought to hold by definition.

This tree representation is the foundation for our algorithm, and serves as the second part of our algorithm ---
transferring the tree monomial to its corresponding loaded tree; this step is at most quadratic in~$n$.
%However, all the main contents so far are already mentioned in~\cite{vertex_splitting}. 
In the next section, we introduce our graphical algorithm for computing the value of a loaded tree,
i.e., the third part of our algorithm.
%, which chiefly reflect the efficiency and conciseness of our algorithm.

\section{The forest algorithm}
In this section, we illustrate a graphical algorithm called {\em the forest algorithm} (Algorithm~\ref{alg:forest}), 
computing the value of a proper loaded tree. We already know that 
the value of a clever tree is one --- we will see later that the case of a clever tree can be viewed as a special case for the forest algorithm. 
%We postpone the correctness proof of the algorithm to later sections. 
%Note that the algorithm is the same as 
%described in~\cite{ACM_communication}, only some terms (names of the trees) are modified. 

First, we introduce the construction of a {\em weighted tree} from a loaded tree.
\begin{definition}
 Let~$LT=(V,E,h,m)$ be a loaded tree. 
 Define the weight function~$w: V\cup E \to \mathbb{N}$ as $w(e)=m(e)-1$ for $e\in E$, $w(v)=|h(v)|+\deg(v)-3$ 
for $v\in V$.
Then, the tuple $T:=(V,E,w)$ is called the {\bf weighted tree of~$LT$}.
\end{definition}
\begin{remark}
From the third item of Definition~\ref{def:loaded_tree},
we see that $w(v)\geq 0$ for all $v\in V$ and the multiplicity of any edge is naturally in~$\mathbb{N}^+$.
Hence the weight assigned to each vertex and each edge in the weighted tree is non-negative.
\end{remark}

\begin{proposition}
 Let $(V,E,w)$ be a weighted tree of some proper loaded tree, then we have $\sum_{v\in V}{w(v)}=\sum_{e\in E}{w(e)}$.
 This identity is called {\bf the weight identity}.
\end{proposition}
\begin{proof}
  \begin{flalign*}
\sum_{v\in V}{w(v)} &=  \sum_{v\in V}{(\deg(v)+|h(v)|-3)}&&\\
                    &= \sum_{v\in V}{\deg(v)}+\sum_{v\in V}{|h(v)|}-3\cdot |V|&&\\
                    &= 2\cdot |E| + n -3\cdot |V|&&\\
                    &= 2\cdot |E| +n -3\cdot |E|-3&&\\
                    &= n-3-|E|\\
\sum_{e\in E}{w(e)} &= \sum_{e\in E}{(m(e)-1)}&&\\
                    &= \sum_{e\in E}{m(e)}-|E|&&\\
                    &= n-3-|E|&&
\end{flalign*}
Note that~$|E|$ refers to the number of edges, not fringes.
\end{proof}

Now we illustrate the forest algorithm.

 \begin{algorithm}
\thispagestyle{empty}
\caption{The forest algorithm}
\label{alg:forest}
\SetKwInOut{Input}{input}
\SetKwInOut{Output}{output}

\Input{ a loaded tree $LT=(V,E,h,m)$ with~$n$ labels and~$n-3$ fringes }
\Output{ the value of the input loaded tree, i.e.~$\int(LT)$ }
\texttt{\\}

$T\gets$ $(V,E,w)$, the weighted tree of~$LT$, where $w:V\cup E\to \mathbb{N}$ is the weight function\;
$S\gets \sum_{e\in E}{w(e)}$\;
$sign \gets (-1)^S$\;
Apply the following operation on the tree~$T$: replace each edge by a length-two edge with a vertex in the middle inheriting 
the weight of the replaced edge (see Figure~\ref{fig:weighted_tree} to Figure~\ref{fig:redundancy_tree} for such a construction), 
obtaining a new tree $T_1=(V_1,E_1,w_1)$ where $w_1:V_1\to \mathbb{N}$. The tree~$T_1$ is called {\bf the redundancy tree of~$LT$}\;
Delete from~$T_1$ vertices with zero-weight and their adjacent fringes, and obtain the {\bf redundancy forest of~$LT$} 
denoted by~$RF$\;
Apply a recursive formula on~$RF$ --- the details of which will be displayed later --- obtaining the absolute 
value~$A$ of~$LT$\;

\Return{$sign\cdot A$}         
\end{algorithm}

\begin{remark}
The sign of the tree value
is~$-1$ to the power of the edge weight sum (or equivalently, the vertex weight sum).
For the argument of this claim on sign, see Remark~\ref{rem:sign}.
The word ``redundancy'' in Algorithm~\ref{alg:forest} reflects some 
information on how far the given loaded tree is away from a clever tree (which is the 
most concise one).
\end{remark}

%\cite[Section 2.]{vertex_splitting}.
%From Remark 3.8. and the reduction chain algorithm (Section 5.2) of~\cite{vertex_splitting},
%we know that if the value of the given loaded tree is non-zero, then 

%This is because each recursion step contributes a negative sign to the value, and from the linear reduction
%(\cite[Section 3]{vertex_splitting}) we know that each recursive step reduces the edge weight sum by one. 
%More details will be explained in Section~\ref{sec:proof_sun_like}.

%\item {\bf redundancy forest:}
%We say that the tuple $(V,E,w)$ is a {\bf redundancy tree} if $(V,E)$ is a tree and $w:V\to \mathbb{N}$
%is a function defined on the vertex set. In Step 3 of the algorithm, we obtain a 
%weighted tree. Then, we replace each edge by two edges with a vertex in the middle ---
%inheriting the weight of the replaced edge --- connecting them. We see that in this way, we actually
%get a redundancy tree. We call the so-gained tree the {\em redundancy tree of the given loaded tree (or, of the given weighted tree)}.
%A {\bf redundancy forest} is defined to be a forest in which each tree is a redundancy tree.
%From the redundancy tree we obtain in Step 5, we will obtain a redundancy forest by
%deleting all vertices of zero weight-value and their incident edges; this so-obtained forest is then 
%called {\em the redundancy forest of the given loaded tree / weighted tree / redundancy tree}.

After we obtain the redundancy forest of the given loaded tree, we apply a recursive formula to each redundancy tree in the forest,
so as to obtain the absolute value of the loaded tree. 
\begin{definition}
Let~$RF$ be the redundancy forest of the loaded 
tree~$LT$, define the {\bf value of~$RF$} (denoted by~$\int(RF)$) as the product of the values of all the trees in the forest.
\end{definition}
\begin{definition}
Define the {\bf value of a redundancy tree $RT=(V,E,w)$} recursively as follows.
Pick any leaf $l\in V$, compare the 
weight of~$l$ with that of its unique parent~$l_1$: if $w(l)>w(l_1)$, return~$0$; otherwise, 
$$\int(RT):={w(l_1)\choose w(l)}\cdot \int(RT_1),$$ where $RT_1=(V_1,E_1,w_1)$ is the redundancy tree defined as 
follows. Delete leaf vertex~$l$ and its incident edge from~$RT$, and then replace the weight of~$l_1$ by $w(l_1)-w(l)$.
We have $V_1=V\setminus \{l\}$, $E_1=E\setminus \{l,l_1\}$, $w_1(l_1)=w(l_1)-w(l)$ and $w_1(v)=w(v)$
for all $v\in V_1\setminus \{l_1\}$. When~$RT$ is a degree-zero vertex, $\int(RT):=0$ if it has non-zero weight and 
$\int(RT):=1$ otherwise. If~$RT$ is a null graph ---  the graph that contains no vertices or edges --- then 
$\int(RT):=1$. 
\end{definition}
\begin{theorem}
 The absolute value of loaded tree~$LT$ equals to the value of the redundancy forest~$RF$. 
\end{theorem}
\begin{proof}
 This states that Algorithm~\ref{alg:forest} is correct. For the proof, see Section~\ref{sec:correctness}.
\end{proof}

Let us see an example, on how to obtain the value of a given loaded tree.

\begin{example}\cite[Example 0.4.]{ACM_communication}\label{eg:forest}
Figure~\ref{fig:loaded_tree} depicts~$LT$ --- a loaded tree with~$14$ labels and~$11$ fringes,
while Figure~\ref{fig:weighted_tree} shows the weighted tree~$WT$ of~$LT$.
We obtain that the edge weight sum of~$WT$ is $2+4+0+1=7$ (which is the same as its vertex weight sum
$1+1+4+0+1=7$). Then we obtain that the sign of~$\int(LT)$ is $(-1)^7=-1$.
Figure~\ref{fig:redundancy_tree} shows the redundancy tree~$RT$ of~$LT$, 
and Figure~\ref{fig:redundancy_forest} describes the corresponding redundancy forest~$RF$.
Then apply the recursive formula (Figure~\ref{fig:formula_new}) on~$RF$, we obtain 
$$\int(RF)=
 [{1 \choose 1} \times 1] \times [{4 \choose 1}\times {4 \choose 3}\times {2 \choose 1} \times {1 \choose 1} \times 1]=32.$$ 
 Hence the value of the loaded tree~$LT$ shown in Figure~\ref{fig:loaded_tree} is~$-32$.

  \begin{figure}
\centering
\begin{subfigure}{0.5\textwidth}
  \centering
  \includegraphics[width=0.9\linewidth]{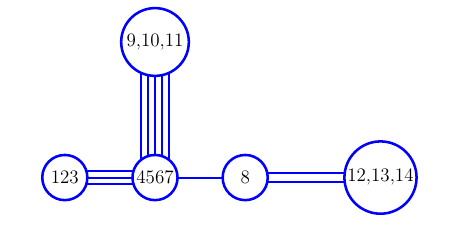}
  \caption{}
  \label{fig:loaded_tree}
\end{subfigure}%
\begin{subfigure}{0.5\textwidth}
  \centering
  \includegraphics[width=0.9\linewidth]{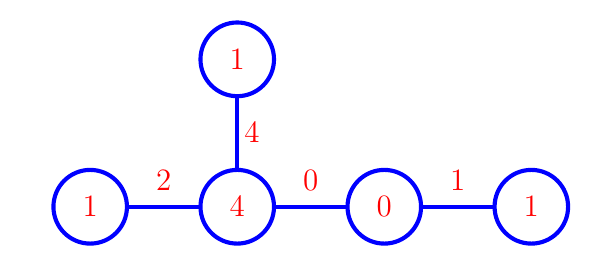}
  \caption{}
  \label{fig:weighted_tree}
\end{subfigure}
\caption{On the left is a loaded tree~$LT$ with~$14$ labels and~$11$ fringes. Labels are tagged in black.
On the right is the weighted tree~$WT$ of the loaded tree~$LT$ described in Figure~\ref{fig:loaded_tree}. 
Weight function values are marked in red. }
\end{figure}

  \begin{figure}
\centering
\begin{subfigure}{0.5\textwidth}
  \centering
  \includegraphics[width=0.9\linewidth]{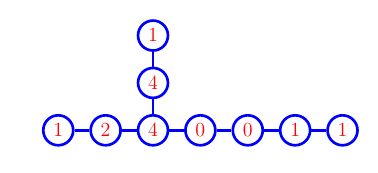}
  \caption{}
  \label{fig:redundancy_tree}
\end{subfigure}%
\begin{subfigure}{0.5\textwidth}
  \centering
  \includegraphics[width=0.9\linewidth]{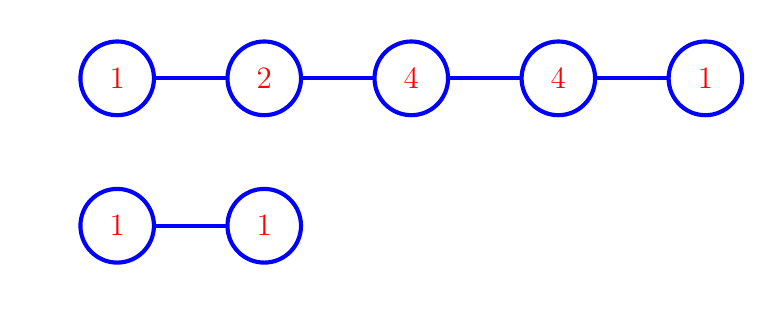}
  \caption{}
  \label{fig:redundancy_forest}
\end{subfigure}
\caption{On the left is the redundancy tree~$RT$ of~$LT$ described in Figure~\ref{fig:loaded_tree}.
Weight values are marked in red.
On the right is the redundancy forest obtained from~$RT$ shown on the left.}
\end{figure}

 \begin{figure}
  \centering
  \includegraphics[width=0.6\linewidth]{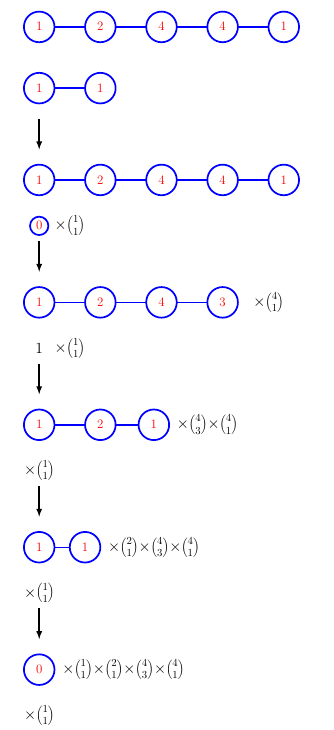}
\caption{This picture shows the process of applying the recursive formula on the two trees in the forest~$RF$ illustrated in
Figure~\ref{fig:redundancy_forest}.}
\label{fig:formula_new}
\end{figure}
\end{example}

%Let us continue with our running example.
%\begin{example}\cite[Example 0.4.]{ACM_communication}
%Figure~\ref{fig:formula_new} shows the process of applying the recursive formula
%on each tree in the redundancy forest $RF$ depicted in Figure~\ref{fig:redundancy_forest}. 
% \begin{figure}
%  \centering
%  \includegraphics[width=0.6\linewidth]{formula_new.png}
%\caption{This picture shows the process of applying the recursive formula on the two trees in the forest $RF$ drawn in
%Figure~\ref{fig:redundancy_forest}.}
%\label{fig:formula_new}
%\end{figure}
%In the last step, we encounter a single vertex with weight zero, hence the tree has value one.
%Combine it with the illustration in the picture we obtain that 
%$$\int(RF)=
% [{1 \choose 1} \times 1] \times [{4 \choose 1}\times {4 \choose 3}\times {2 \choose 1} \times {1 \choose 1} \times 1]=32.$$ 
% Product of $32$ and $-1$ tells us that the value of the loaded tree $LT$ shown in Figure~\ref{fig:loaded_tree} is $-32$.
%\end{example}
\begin{proposition}
 Algorithm~\ref{alg:forest} is an algorithm; we call it {\bf the forest algorithm}. 
\end{proposition}
\begin{proof}
Consider the whole procedure, from inputting a proper loaded tree, to finally obtaining the tree value. Termination is trivial since 
the input tree is finite, so does its redundancy forest.
In the recursion formula calculus, each step strictly reduces the size of the redundancy forest. 
Because of the following identity of binomial coefficients, we know that the recursive formula gives 
us the same value, no matter in which sequence we consider and delete the leaf vertices of a redundancy tree,
when the two leaves share a common parent:
$${c\choose a}\cdot {c-a \choose b} \equiv {c\choose b}\cdot {c-b \choose a},$$
where~$c$ is the weight of the common parent of two leaves whose weights are~$a\,,b$ respectively. 

In the case when they do not, 
we can easily argue it by no matter in which order we delete these two leaves, the obtained value 
equals $${a_1 \choose a}\cdot {b_1 \choose b}\cdot rf',$$ where~$rf'$ is the value of the new redundancy forest 
after deleting the two leaves,~$a$,~$b$ are the weights of the two leaves respectively, and~$a_1$,~$b_1$ 
are the weights of their parents respectively.

Therefore, the above process is indeed an algorithm. 
\end{proof}
\begin{remark}
It is not hard to see that
the complexity of the forest algorithm is linear with respect to the number of vertices of the 
input loaded tree.
If the input loaded tree is a clever tree, then each vertex of its redundancy tree has weight zero.
Therefore, its redundancy forest is a null graph and hence the input loaded tree has value one. 
This indicates that the value of clever trees can also be handled by the forest algorithm, as a special case.
\end{remark}
Based on the forest algorithm, let us consider again the extension of the one-to-one correspondence between loaded trees and tree monomials:
a loaded tree with a single vertex and~$3$ labels has the null graph as its redundancy forest, hence has value one;
a loaded tree with a single vertex and~$k$ ($k\neq 3$) labels has the single vertex with non-zero weight as its 
redundancy forest, hence has zero value. 
And the correctness on these base cases directly come from the definition of the degree map~$\int$.
Given the correctness of the forest algorithm (see Section~\ref{sec:correctness}), our previous extension stands defensible.
In order to verify the base case of the forest algorithm, we need the concept of ``linear reduction'': algebraic and graphical.

\section{Linear reduction}
In this section, first we introduce an algebraic reduction of a monomial in~$A^{n-3}(\overline{\mathcal{M}}_{0,n})$, using Keel's linear 
relation and Keel's quadratic relation. Then we give the equivalent graphical characterization for it. 
\subsection{Algebraic linear reduction}\label{sec:algebraic_linear_reduction}
 Keel's linear relation was originally 
proved in~\cite[Theorem 1.(2)]{keel}; we state exactly the same content as follows, but in different notations.

\begin{fact}[Keel's linear relation,~\cite{keel} Theorem 1.(2)]\label{def:keel_linear_relation}
 Denote by $$\epsilon_{ij\mid kl}:=\sum_{i,j\in I,k,l\in J}{\delta_{I,J}}.$$ Then we have the equality relations 
 $$\epsilon_{ij\mid kl}=\epsilon_{il\mid kj}=\epsilon_{ik\mid jl}.$$ We call it {\bf Keel's linear relation}.
\end{fact}

Let us see a concrete example on it.
\begin{example}
 When~$n=6$, we have $\epsilon_{12\mid 35}=\epsilon_{13\mid 25}=\epsilon_{15\mid 23}$, i.e., 
 $$\delta_{12,3456}+\delta_{124,356}+\delta_{126,345}+\delta_{1246,35}$$
 $$=\delta_{13,2456}+\delta_{134,256}+\delta_{136,245}+\delta_{1346,25}$$
 $$=\delta_{15,2346}+\delta_{145,236}+\delta_{156,234}+\delta_{1456,23}$$
\end{example}

\begin{remark}
 From the example above we see that we can substitute some~$\delta_{I,J}$, say~$\delta_{12,3456}$,
 by $\epsilon_{13\mid 25}-(\epsilon_{12\mid 35}-\delta_{12,3456})$. Basically we can replace~$\delta_{I,J}$
 by a sum of~$(2^{n-3}-1)$ many~$(\pm)\delta_{I',J'}$s.
\end{remark}
\begin{definition}
A {\bf linear reduction} on a non-clever tree monomial~$M$ in~$A^{\bullet}(\overline{\mathcal{M}}_{0,n})$ is defined as follows.
Replace a factor which has power more than one by Keel's linear reduction, and then eliminate all value-zero terms using 
Keel's quadratic relation. 
\end{definition}
Let us see some examples. 

 \begin{example}\label{eg:one_tree_generated}
  Apply Keel's linear relation to the tree monomial $\delta^2_{12,3456}\cdot \delta_{1234,56}$, 
  replacing one occurrence of~$\delta_{12,3456}$
 by $\epsilon_{13\mid 25}-(\epsilon_{12\mid 35}-\delta_{12,3456})$:
  $$\delta^2_{12,3456}\cdot \delta_{1234,56}= \delta_{12,3456}\cdot \delta_{1234,56}\cdot (\epsilon_{13\mid 25}-\delta_{124,356}-\delta_{126,345}
 -\delta_{1246,35}).$$ 
 Since any summand~$\delta_{I,J}$ of~$\epsilon_{13\mid 25}$ has $1,3\in I$ and $2,5\in J$, any such~$\delta_{I,J}$ together with~$\delta_{12,3456}$ fulfills 
 Keel's quadratic relation and we obtain that $\delta_{12,3456}\cdot \delta_{I,J}=0$. Consequently, we have 
 $\delta_{12,3456}\cdot \epsilon_{13\mid 25}=0$. Hence we have:
  $$\delta^2_{12,3456}\cdot \delta_{1234,56}= \delta_{12,3456}\cdot \delta_{1234,56}\cdot (-\delta_{124,356}-\delta_{126,345}
 -\delta_{1246,35}).$$ 
 
 One can check that the two pairs
 $(\delta_{1234,56},\delta_{126,345})$ and $(\delta_{1234,56}, \delta_{1246,35})$ both fulfill Keel's quadratic relation. 
 Hence both products are zero.
 Therefore, we have 
 $$\delta^2_{12,3456}\cdot \delta_{1234,56}= -\delta_{12,3456}\cdot \delta_{1234,56}\cdot \delta_{124,356}.$$ 
 
 A linear reduction on the given monomial is now accomplished. 
\end{example}

\begin{example}\label{eg:two_reductions_needed}
Let $M:=\delta^3_{123,4567}\cdot \delta_{12345,67}\in A^{4}(\overline{\mathcal{M}}_{0,7})$ be the given monomial.
We use Keel's linear relation, replacing~$\delta_{123,4567}$ via $\epsilon_{12\mid 46}=\epsilon_{14\mid 26}$. 
Then, we obtain
$$\delta^3_{123,4567}\cdot \delta_{12345,67}=\delta^2_{123,4567}\cdot \delta_{12345,67}\cdot (\epsilon_{14\mid 26}-(\epsilon_{12\mid 46}-\delta_{123,4567})).$$
Then we see that each summand of~$\epsilon_{14\mid 26}$ fulfills Keel's quadratic
relation together with~$\delta_{123,4567}$, hence 
we have $\epsilon_{14\mid 26}\cdot \delta_{123,4567}=0$. Hence we have:
\begin{flalign}
\delta^3_{123,4567}\cdot \delta_{12345,67} &= \delta^2_{123,4567}\cdot \delta_{12345,67}\cdot (-(\epsilon_{12\mid 46}-\delta_{123,4567})) &&\\\nonumber
            &= \delta^2_{123,4567}\cdot \delta_{12345,67}\cdot (-\delta_{12,34567}-\delta_{125,3467}  &&\\\nonumber
            & \ \  \ \ -\delta_{127,3456} -\delta_{1235,467}-\delta_{1237,456}-\delta_{1257,346} &&\\\nonumber
            & \ \  \ \ -\delta_{12357,46}) &&
\end{flalign}
Then we need to exclude those summands of~$\epsilon_{12\mid 46}$ which fulfill 
Keel's quadratic relation with any factor(s) of~$M$ --- here it refers to~$\delta_{12345,67}$ and~$\delta_{123,4567}$.
After the exclusion, we obtain that 
\begin{flalign}\label{eq:two_monomials}
\delta^3_{123,4567}\cdot \delta_{12345,67} &= \delta^2_{123,4567}\cdot \delta_{12345,67}\cdot (-\delta_{12,34567}-\delta_{1235,467}) &&\\\nonumber
                                          &= -\delta^2_{123,4567}\cdot \delta_{12345,67}\cdot \delta_{12,34567} &&\\\nonumber
                                          & -\delta^2_{123,4567}\cdot \delta_{12345,67}
                                          \cdot\delta_{1235,467}&&
\end{flalign}
\qed
\end{example}

\begin{remark}[sign]\label{rem:sign}
We see that in the linear reduction process, whenever we replace one occurrence of~$\delta_{I,J}$ by some 
$\epsilon_{ik\mid jl}-(\epsilon_{ij\mid kl}-\delta_{I,J})$, we can directly omit~$\epsilon_{ik\mid jl}$, since any 
summand of it fulfills the Keel's quadratic relation with~$\delta_{I,J}$ and there is at least one occurrence of 
$\delta_{I,J}$ still left in the remaining part of the monomial. Hence from 
now on we will only say that we replace~$\delta_{I,J}$ by $-(\epsilon_{ij\mid kl}-\delta_{I,J})$. 
From this analysis we also observe that
whenever we do one-time linear reduction, we obtain a negative sign on the right hand side.
Therefore, how many times of linear reduction we use decides the sign for the value of the given
monomial --- odd times gives a negative sign while even times leads to a positive sign.
\end{remark}

When we consider the linear reduction of a tree monomial, we see a parallel process of the linear reduction, on the tree representation of the given monomial:
first we should decide on which multi-edge of the loaded tree to reduce, then we should pick up 
a quadruple to do the reduction. In the sequel, we introduce the next step: how to directly tell which loaded trees
or tree monomials will be generated on the right hand side of the equation in the algebraic linear reduction 
but using a graphical method which is called {\em vertex splitting}.

\subsection{Vertex splitting}\label{subsec:vertex_splitting}
In this section, we consider the parallel process of the linear reduction on a tree monomial, on the tree representation. 
Before we can explain this concept, we need to define the {\em branches of a vertex} first.
\begin{definition}
 Let~$v$ be a vertex in the loaded tree~$T$. Removing vertex~$v$ but not its incident edges
 gives us degree of~$v$ many parts. Each part 
 is a structure of a tree but lacking a vertex.
 We make~$\deg(v)$ many copies of vertex~$v$, and 
 concatenate it to each of these parts at the place where a vertex is missing. 
 \begin{enumerate}
\item Then we obtain~$\deg(v)$ many trees, we call them {\bf branches of~$v$}.
\item We call the copy of~$v$ the {\bf special vertex} in each of these trees.
\item We say that a branch with the special vertex~$v$ is {\bf attached to} some vertex~$v'$ (of some tree~$T'$) if we add an extra edge~$e'$ between 
~$v$ and~$v'$, and then contract it (two vertices merged into a new vertex~$v_1$), and then set the labeling set of~$v_1$ to be the labeling 
set of~$v'$.
\end{enumerate}
\end{definition}
\begin{remark}
 We see that the operation of attaching a branch to another tree~$T'$, has in principle nothing to do with the starting tree~$T$.

Note that each branch corresponds to a proper cluster of~$v$. Each cluster
 of~$v$ corresponds either to a branch of~$v$ or to a label of~$v$.
 In our concrete operation, we can think of the branch simply as 
 the structure of a tree but lacking that special vertex. When we attach a branch to some vertex~$v'$ of another tree~$T'$, we can just view it as 
to concatenate~$T'$ to the branch at vertex~$v'$, to the endmost of that branch where a vertex is missing. 
From now on, for convenience, 
we will use this simplified concept of branch and branch-attaching.
 \end{remark}

\begin{figure}
\centering
\begin{subfigure}{.32\textwidth}
  \centering
  \includegraphics[width=0.95\linewidth]{m6_with_notation.png}
  \caption{}
\end{subfigure}
\begin{subfigure}{.32\textwidth}
  \centering
  \includegraphics[width=0.95\linewidth]{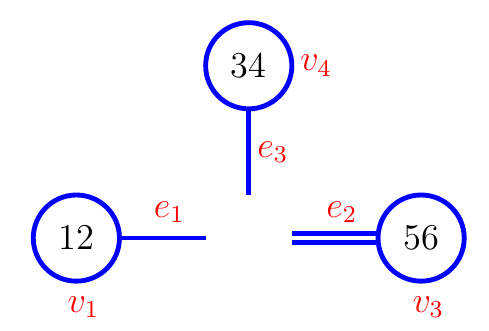}
  \caption{}
\end{subfigure}
\begin{subfigure}{.34\textwidth}
  \centering
  \includegraphics[width=0.75\linewidth]{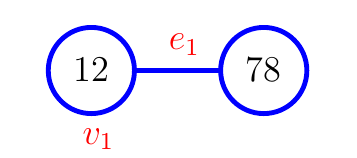}
  \caption{}
\end{subfigure}
\caption{In the middle are the (simplified) branches of vertex~$v_2$ of the loaded tree on the left. 
Vertex~$v_2$ has three branches. 
On the right is the loaded tree that is obtained from attaching the branch of~$v_2$ containing labels~$1$ and~$2$
to a single vertex with labeling set~$\{7,8\}$. }
\label{fig:branches_etc}
\end{figure}

In order to better understand the concept, we see an example. Let us turn our focus back to vertex~$v_2$ in the loaded tree in 
Figure~\ref{fig:loaded_trees}(b); we illustrate it again here in  Figure~\ref{fig:branches_etc}(a).
See Figure~\ref{fig:branches_etc}(b) for its branches. See Figure~\ref{fig:branches_etc}(c) for an example of branch-attaching.

Now we explain {\bf the vertex splitting} process as a series of operations on a loaded tree, see Algorithm~\ref{alg:vertex_splitting}.
Although we express it in an algorithmic environment --- for the sake of consistency in the format ---
note that this process is not deterministic, there are some degrees of freedom.
Recall that the weight function is defined as $\deg(v)+|h(v)|-3$ for 
each vertex~$v$. 

 \begin{algorithm}
\thispagestyle{empty}
\caption{Vertex splitting}
\label{alg:vertex_splitting}
\SetKwInOut{Input}{input}
\SetKwInOut{Output}{output}

\Input{ a non-clever loaded tree $T=(V,E,h,m)$; 
a multi-edge $e=\{v_1,v_2\}$ (with the corresponding cut $\{I,J\}$); a quadruple set $Q:=\{i,j,k,l\}$ 
such that $i,j\in I$, $k,l\in J$ --- w.l.o.g. assume that~$v_1$ is in Component-$I$.  }
\Output{ NULL or a loaded tree $\tilde{T}=(\tilde{V},\tilde{E},\tilde{h},\tilde{m})$. }
\texttt{\\}

\If{both~$v_1$ and~$v_2$ have zero-weight}
  {\Return NULL}
  \Else
   {
   \begin{enumerate}
    \item   Pick one of the incident vertices of~$e$ with non-zero weight --- assume w.l.o.g. that vertex~$v_1$ is chosen\;
   \item Construct the structure~$S$ of two vertices~$v'_1$ and~$v''_1$ connected by a single edge $e'=\{v'_1,v''_1\}$.
   (The idea is that this structure will replace~$v_1$ later on.)\;
   \item Set the labeling sets of~$v'_1$ and that of~$v''_1$ so that $\tilde{h}(v'_1)\cup \tilde{h}(v''_1)=h(v_1)$ 
   and~($\star_1$) $h(v_1)\cap (Q\cap I) \subset \tilde{h}(v'_1)$;
   \subitem ($\star_2$) However, note that if~$v_1$ has no other branches except for the ones that
   contain any label in $Q\cap I$ or edge~$e$, then $|\tilde{h}(v'_1)|\geq 1$ must hold\;
   Attach the branches of~$v_1$ to the structure~$S$ at~$v'_1$ or~$v''_1$ as follows:
   \begin{enumerate}
    \item Among the branches of~$v_1$, those containing any label in~$Q\cap I$ are
    attached to~$v'_1$ and the one that contains labels in $Q\cap J$ is attached to~$v''_1$ --- it is not hard to check that the two labels
	in $Q\cap J$ are in the same branch of~$v_1$. The branch 
	containing~$e$ should be modified slightly: the multiplicity of edge~$e$ 
	in this branch gets reduced by one, then gets attached to~$v''_1$.
      \item  The other branches of~$v_1$ in~$T$ can be either attached to~$v'_1$ or~$v''_1$.
      \item ($\star_3$) However, note that if $\tilde{h}(v''_1)=\emptyset$, then we must attach at least one branch to~$v'_1$, except for the one containing~$e$.
   \end{enumerate}
   \end{enumerate}
   }
\end{algorithm}

\begin{remark}\label{rem:vertex_weight_non_zero}
By~($\star_1$) we want to note that~$v'_1$ and~$v''_1$ are symmetric. It is also no problem if here we distribute to~$v'_1$ 
instead of~$v''_1$ the labels in $h(v_1)\cap (Q\cap I)$; then we just need to interchange the roles of~$v'_1$ and~$v''_1$, then 
conduct the remaining steps of the process. The process in Algorithm~\ref{alg:vertex_splitting} can be roughly viewed as~$v_1$ 
 being replaced by the structure~$S$. 
\end{remark}

 \begin{lemma}
 The output~$\tilde{T}$ in Algorithm~\ref{alg:vertex_splitting} is a loaded tree, and is of the same type with the input tree~$T$. 
 Thence when~$T$ is proper,~$\tilde{T}$ is also proper. 
 \end{lemma}
 \begin{proof}
 First, let us look more into details of conditions~($\star_2$) and~($\star_3$). 
 In~($\star_2$), there are two branches containing labels in $Q\cap I$ for~$v_1$, plus the branch of~$e$; these are 
 three branches of~$v_1$. If~$v_1$ has no other branches, then it must have at least one label that can be freely distributed 
 (to~$v'_1$ or~$v''_1$) since $w(v_1)\geq 1$. Hence, we can always require $|\tilde{h}(v''_1)|\geq 1$, under this situation. 
 In~($\star_3$), since $|\tilde{h}(v''_1)|=0$,~($\star_2$) tells us that~$v_1$ has at least one branch that is free
 to be distributed. This is why we can ensure the described arrangement. 
 Analogously, one can verify that the operation in~($\star_3$) can be guaranteed. We leave it to the readers as an exercise.

  Hence the requirements~($\star_2$) and~($\star_3$) guarantee that vertex~$v''_1$ and~$v'_1$ fulfill the third condition 
 in Definition~\ref{def:loaded_tree}. Plus, obviously $$\{\tilde{h}(\tilde{v}\in\tilde{V})\}_{\tilde{v}\in\tilde{V},\; \tilde{h}(\tilde{v})\neq \emptyset}$$
 form a partition of the labeling set of~$T$. Hence~$\tilde{T}$ is a loaded tree.
 It is not hard to see that~$T$ and~$\tilde{T}$ have the same set of labels and the same number 
 of fringes. Hence,~$\tilde{T}$ is of the same type with~$T$. 
 \end{proof}
 
 \begin{proposition}
  The process stated in Algorithm~\ref{alg:vertex_splitting} is an algorithm.
 \end{proposition}
\begin{proof}
 From the above analysis we see that the process terminates, and returns a loaded tree or NULL.
\end{proof}

 The following result states some relation between the weight functions of~$T$ and~$\tilde{T}$. 
 \begin{lemma}
 Let~$w$ be the weight function for~$T$ and let~$\tilde{w}$ be that of~$\tilde{T}$. Then we have
 $\tilde{w}(v'_1)+\tilde{w}(v''_1) = w(v_1)-1$.
 \end{lemma}
 \begin{proof}
 \begin{flalign*}
  \tilde{w}(v'_1)+\tilde{w}(v''_1) &= (\deg(v'_1)+|\tilde{h}(v'_1)|-3) + (\deg(v''_1)+|\tilde{h}(v''_1)|-3) &&\\
                    &=  (\deg(v'_1)+\deg(v''_1)) + (|\tilde{h}(v'_1)| + |\tilde{h}(v''_1)|)-6  &&\\
                    &=  (\deg(v_1)+2) + |h(v_1)| -6  &&\\
                    &=  \deg(v_1) + |h(v_1)| -4  &&\\
                    &=  ( \deg(v_1) + |h(v_1)| -3) -1 &&\\
                    &= w(v_1) -1 &&
\end{flalign*}
\end{proof}
  We see from the above reasoning that the vertex-splitting process indeed requires the split vertex to have non-zero
  weight. 
  \begin{proposition}\label{prop:two_vertices_weight_zero}
  If there
 is a multi-edge $e=\{v_1,v_2\}$ in some proper loaded tree, both~$v_1$ and~$v_2$ being
 weight-zero leads to the value of this tree being zero. 
 \end{proposition}
 \begin{proof}
  We postpone it to Section~\ref{sec:the_missing_proofs}.
 \end{proof}
 If~$v_1$ has zero weight, while~$v_2$ has non-zero weight,
 we can simply exchange the names of the two vertices and then continue the vertex-splitting process.
\begin{definition}
  We say that in Algorithm~\ref{alg:vertex_splitting}, vertex~$v_1$ {\bf is split into} 
 vertices~$v'_1$ and~$v''_1$. 
\end{definition}

 \begin{remark}\label{rem:finite}
It is not hard to check, that the weight sum of edges 
 of~$\tilde{T}$ is always one less than that of~$T$. Therefore, if we recursively apply this process,
 then after finitely many steps, we will obtain only clever trees. This observation provides
 us an idea on calculating the value of a given tree monomial.
\end{remark}

In the above algorithm, output is just one loaded tree. We observe that we actually have
some freedom at several steps:
\begin{enumerate}
 \item If both~$v_1$ and~$v_2$ have non-zero weights, then we can split either of them.
 \item We could also have some freedom on how to set up the labeling function 
 for~$v'_1$ and~$v''_1$, respectively --- as long as condition~($\star_2$) is fulfilled.
 \item Also, we have some freedom on the arrangements of branches of~$v'_1$ and those of~$v''_1$
  --- as long as condition~($\star_3$) is fulfilled.
\end{enumerate}
When we consider all these freedom, and collect all the possibly generated loaded trees, 
we obtain the tree-version linear reduction algorithm ---
this algorithm 
does the same thing and should give us the same result 
as the algebraic linear reduction which is introduced in Section~\ref{sec:algebraic_linear_reduction}.
We will see some examples on the vertex splitting process (Algorithm~\ref{alg:vertex_splitting}) in Section~\ref{sec:tree_version_reduction_algorithm},
after the tree-version linear reduction is introduced.

\subsection{Tree-version linear reduction algorithm}\label{sec:tree_version_reduction_algorithm}
In this subsection, we explain the tree-version linear reduction algorithm, see 
Algorithm~\ref{alg:tree_version_reduction}. 

\begin{algorithm}
\thispagestyle{empty}
\caption{tree-version linear reduction}\label{alg:tree_version_reduction}
\SetKwInOut{Input}{input}
\SetKwInOut{Output}{output}

\Input{a proper non-clever loaded tree~$T_M$ (the corresponding loaded tree of a proper monomial~$M$); 
a multi-edge $e=\{v_1,v_2\}$ (with corresponding cut $\{I,J\}$) of~$T_M$; 
a quadruple set $\{i,j,k,l\}$ such that 
$i,j\in I$ and $k,l\in J$.}
\Output{loaded trees whose corresponding monomials are the ones on the right hand side of the equation, after a step of a linear reduction
on~$M$ to reduce one occurrence of~$\delta_{I,J}$ which uses Keel's
linear reduction on the relation $\epsilon_{ij\mid kl}=\epsilon_{ik\mid jl}$. }

$w_1 \gets$ weight of~$v_1$\;
$w_2 \gets$ weight of~$v_2$\;
\If{$w_1=0$ and $w_2=0$}
    {\Return~$\emptyset$}
$ST\gets$ the set of all loaded trees that can be obtained from~$T_M$ with edge
$e$ and the quadruple set $\{i,j,k,l\}$, after the 
 vertex splitting process (Algorithm~\ref{alg:vertex_splitting}) --- either by splitting vertex~$v_1$, or by splitting vertex~$v_2$ --- in the set~$ST$\;

\Return{$ST$}
  
\end{algorithm}

\begin{remark}
 Let~$M$ be a proper tree monomial. First, transfer it to its corresponding loaded tree~$T_M$
using Algorithm~\ref{alg:monomial_to_tree}. Then, apply the tree-version linear reduction 
to it via the quadruple $\{i,j,k,l\}$, obtaining~$ST$. 
Then, write down the negative sum of all monomials~$M'$ such that $T_{M'}\in ST$. 
What we get now is already the result of one step of a good linear reduction using the relation $\epsilon_{ij\mid kl}=\epsilon_{ik\mid jl}$.
The above algorithm characterize the parallel process of a linear reduction on the tree representation for a tree monomial.
\end{remark}

\begin{remark}
Note that this algorithm is also applicable to non-proper monomials (trees). However, we focus on proper 
monomials (trees), and 
we always obtain one or more proper monomials (loaded trees) after each linear reduction.
\end{remark}

\begin{theorem}\label{thm:correctness_tree_version}
 Algorithm~\ref{alg:tree_version_reduction} is correct.
\end{theorem}
\begin{proof}
 We postpone this correctness proof to Section~\ref{sec:the_missing_proofs}.
\end{proof}

\begin{figure}
\centering
\includegraphics[width=0.35\linewidth]{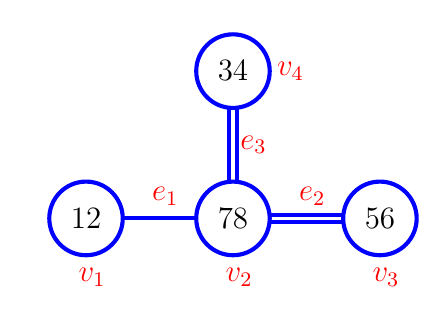}

\caption{This a proper loaded tree~$T$ with~$8$ labels and~$5$ fringes. We want to 
reduce edge~$e_2$ with the quadruple set $\{1,7,5,6\}$, via splitting the
vertex~$v_2$.}
\label{fig:vertex_splitting}

\end{figure}

\begin{figure}[!htb]
\minipage{0.32\textwidth}
  \includegraphics[width=\linewidth]{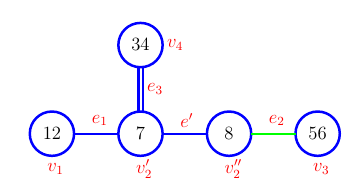}

\endminipage\hfill
\minipage{0.32\textwidth}
  \includegraphics[width=\linewidth]{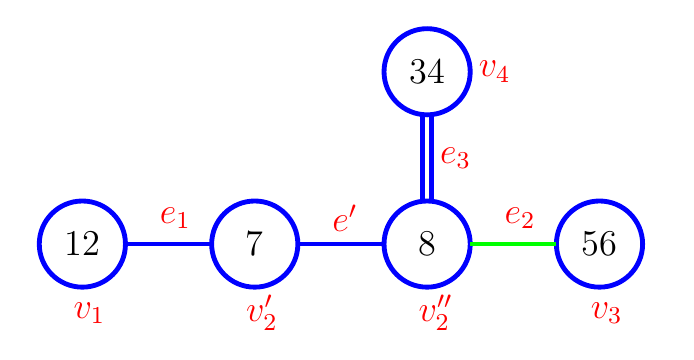}
 
\endminipage\hfill
\minipage{0.32\textwidth}
  \includegraphics[width=\linewidth]{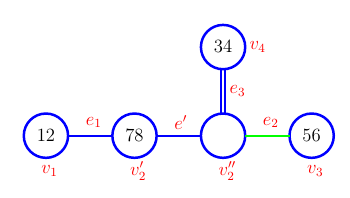}

\endminipage

\caption{The set of these three loaded trees is the output of Algorithm~\ref{alg:tree_version_reduction} --- applied to
the loaded tree~$T$ with multi-edge~$e_2$ (in Figure~\ref{fig:vertex_splitting}) and the quadruple set $\{1,7,5,6\}$. 
The edge which gets reduced is marked in green. Vertex~$v_2$ is split into~$v'_2$ and~$v''_2$. The new edge is denoted
as~$e'$.
Each of them is a loaded tree with~$8$ labels and~$5$ fringes.}
\label{fig:vertex_split}
\end{figure}

Now we see an example for a better understanding of the above algorithm, and as well of the
vertex splitting algorithm (Algorithm~\ref{alg:vertex_splitting}) stated earlier.
\begin{example}\label{eg:three_trees_generated}
Consider the proper loaded tree~$T$ in Figure~\ref{fig:vertex_splitting}. Let us follow 
 the above vertex-splitting process, see what we will obtain.
 We pick a multi-edge~$e_2$ and the quadruple set $\{1,7,5,6\}$.
 We calculate the weight of its incident vertices, find
 out that the
 weight of~$v_3$ is zero while that of~$v_2$ is~$2$. 
 Hence we can only choose~$v_2$ to split (into~$v'_2$ and~$v''_2$), and note that~$v_2$ is in Component-$\{1,2,3,4,7,8\}$ of~$T$. 
 Denote by~$\tilde{T}=(\tilde{V},\tilde{E},\tilde{h},\tilde{m})$ any tree in the output of Algorithm~\ref{alg:tree_version_reduction}.
 
 We should bipartition
 the labeling set of~$v_2$ such that 
 $$\{7,8\}\cap (\{1,2,3,4,7,8\}\cap \{1,7,5,6\})=\{7\}\subset h_2(v'_2)$$ 
 holds.
 So if we follow the instructions in Algorithm~\ref{alg:tree_version_reduction}, we have two options in this step: 
\begin{enumerate} 
\item $\tilde{h}(v'_2)=\{7\}$ and $\tilde{h}(v''_2)=\{8\}$;
\item $\tilde{h}(v'_2)=\{7,8\}$ and $\tilde{h}(v''_2)=\emptyset$.
\end{enumerate}
 
 Then, based on the above two choices, we distribute the branches of~$v_2$ to be the branches of~$v'_2$ or those of~$v''_2$. 
 From the requirements in the vertex-splitting process, we know that 
 the branch containing~$1$ should be attached to~$v'_2$ and the branch containing~$5$ 
 or~$6$ should be attached to~$v''_2$; edge multiplicity of~$e_2$ should be reduced by one.
 The remaining branch --- the branch containing 
 labels~$3$ and~$4$ can be attached to either~$v'_2$ or~$v''_2$ in the first label-distribution option;
 but it must be attached to~$v''_2$ in the second label-distribution option because of the~($\star_3$) requirement. 
 Hence, we obtain in total three loaded trees, after applying Algorithm~\ref{alg:tree_version_reduction}.
 
 Figure~\ref{fig:vertex_split} shows these loaded trees. 
It is not hard to check that each of them is still a loaded tree, and is of the same 
 type with~$T$. 
\end{example}

In the sequel, we look back on our examples in Section~\ref{sec:algebraic_linear_reduction}.
We will apply the tree-version linear reduction algorithm and see if we will
obtain the same result as if we apply the algebraic method, in the ``corresponding tree space''.

\begin{example}
See Figure~\ref{fig:one_tree_generated}(a) for the corresponding loaded tree~$T$ of monomial 
$M:=\delta^2_{12,3456}\cdot \delta_{1234,56}$ in Example~\ref{eg:one_tree_generated}. 
Now we apply the tree-version
linear reduction algorithm to it. 
We pick the multi-edge~$e_1$ to reduce.
Its corresponding cut is $\{\{3,4,5,6\},\{1,2\}\}$.
The quadruple set we pick here is $\{3,5,1,2\}$. 
By an easy calculation we know that the weight of~$v_1$
is zero and that of~$v_2$ is~$1$. Therefore, we can only split vertex~$v_2$ 
--- no freedom of choice here.
And~$v_2$ is in Component-$\{3,4,5,6\}$ of~$T$.
 Denote by~$\tilde{T}=(\tilde{V},\tilde{E},\tilde{h},\tilde{m})$ any loaded tree in the output set.

 First we split vertex~$v_2$ 
 into~$v'_2$ and~$v''_2$ such that 
 $$\{3,4\}\cap (\{3,4,5,6\}\cap \{3,5,1,2\})=\{3\}\subset \tilde{h}(v'_2).$$
Observe that~$v_2$ has only two branches, one contains~$5$, the other
contains edge~$e_1$. Therefore, by condition~$(\star_2)$, ~$\tilde{h}(v''_2)$ 
should contain at least one label. Hence we have that $\tilde{h}(v''_2)=\{4\}$.
We see that there is also no freedom of choice in this step.

The next step is to arrange the branches of~$v_2$ to be the branches of~$v'_2$ 
or those of~$v''_2$. The branch containing~$e_1$ should be modified --- 
multiplicity of edge~$e_1$ should get reduced by one --- and 
then gets attached to~$v''_2$. The branch containing label~$5$ should be attached 
to~$v'_2$. We see that there is also no freedom of different options in this step. We
obtain only one loaded tree, which is also a clever tree.

Hence we only need one time linear reduction for calculating the value of~$T$,
the sign for the result is then $(-1)^1=-1$.
 In Figure~\ref{fig:one_tree_generated}(b) we see the output tree. 
 We can easily obtain that the corresponding monomial 
 of this new tree is $\delta_{12,3456}\cdot \delta_{124,356}\cdot \delta_{1234,56}$, 
 which coincides with the result 
 we obtain in Example~\ref{eg:one_tree_generated}, with the algebraic linear reduction.
 
 \begin{figure}
\centering
\begin{subfigure}{.5\textwidth}
  \centering
  \includegraphics[width=0.7\linewidth]{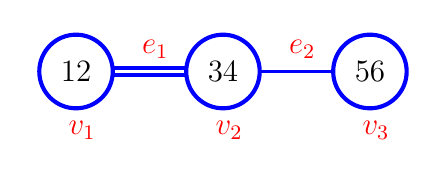} 
  \caption{}
\end{subfigure}%
\begin{subfigure}{.5\textwidth}
  \centering
  \includegraphics[width=\linewidth]{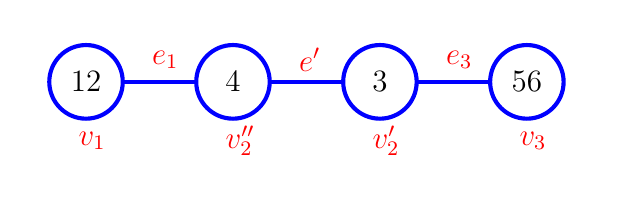}
  \caption{}
\end{subfigure}
\caption{On the left is the corresponding loaded tree~$T_M$ of monomial $M = \delta^2_{12,3456}\cdot
\delta_{1234,56}$. On the right is the only tree that we obtain after reducing edge~$e_1$ with the quadruple set
$\{3,5,1,2\}$, applying the tree-version linear reduction algorithm.
 The new edge is denoted by~$e'$. }
\label{fig:one_tree_generated}
\end{figure}

\end{example}

In the sequel, we apply the tree-version linear reduction algorithm to the monomial in Example~\ref{eg:two_reductions_needed}.

\begin{example}\label{eg:half}
 See Figure~\ref{fig:m7} for the corresponding loaded tree~$T_M$ of monomial 
 $M:=\delta^3_{123,4567}\cdot \delta_{12345,67}$ in Example~\ref{eg:two_reductions_needed}. This is a proper tree.
 The edge $e_1=~\{v_1,v_2\}$ to be reduced is the corresponding edge of cut $\{\{1,2,3\},\{4,5,6,7\}\}$. 
 We pick~$1,2$ from two distinct clusters of~$v_1$ and~$4,6$ from two distinct
 clusters of~$v_2$.
 
 We see that we can either split~$v_1$ or~$v_2$ since they both have non-zero weight.
 When we split 
 vertex~$v_1$, there is only one loaded tree~$\tilde{T}_1$
 that can be obtained, see the loaded tree in Figure~\ref{fig:m7'}. When we split 
 vertex~$v_2$, there is only one loaded tree~$\tilde{T}_2$
 that can be obtained, see the loaded tree in Figure~\ref{fig:m7''}. Hence after 
 we apply one time tree-version linear reduction algorithm to~$T$,
 we obtain in total two loaded trees:~$\tilde{T}_1$ and~$\tilde{T}_2$.
 Their corresponding monomials are 
 $\delta_{12,34567}\cdot \delta^2_{123,4567}\cdot \delta_{12345,67}$
 and $\delta^2_{123,4567}\cdot\delta_{1235,467}\cdot \delta_{12345,67}$, respectively. 
 This 
 coincides with the resulting monomials in Equation~\ref{eq:two_monomials}.

 \begin{figure}
\centering
\includegraphics[width=0.35\linewidth]{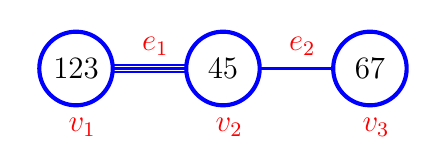}
\caption{This is the corresponding loaded tree~$T_M$ of monomial 
$M=\delta^3_{123,4567}\cdot \delta_{12345,67}$. We want to reduce edge~$e_1$
with the quadruple set $\{1,2,4,6\}$.}
\label{fig:m7}
\end{figure}

\begin{figure}
\centering
\begin{subfigure}{.5\textwidth}
  \centering
  \includegraphics[width=0.92\linewidth]{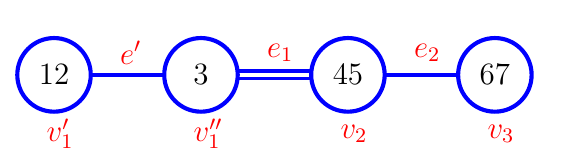} 
  \caption{loaded tree~$\tilde{T}_1$}
  \label{fig:m7'}
\end{subfigure}%
\begin{subfigure}{.5\textwidth}
  \centering
  \includegraphics[width=0.92\linewidth]{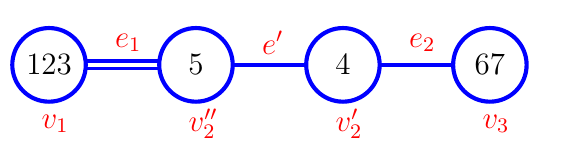}
  \caption{loaded tree~$\tilde{T}_2$}
  \label{fig:m7''}
\end{subfigure}
\caption{The set of these two loaded trees is the output set of Algorithm~\ref{alg:tree_version_reduction}
applied to tree~$T_M$ in Figure~\ref{fig:m7} with edge~$e_1$ and quadruple set $\{1,2,4,6\}$. When we choose to 
split vertex~$v_1$, we obtain the loaded tree in Figure~\ref{fig:m7'}.
When we choose to 
split vertex~$v_2$, we obtain the loaded tree in Figure~\ref{fig:m7''}.
 }
\label{fig:m7'_and_m7''}
\end{figure}

\end{example}

We see that in the above two examples, we obtain exactly the same result as when we did it
via algebraic approach. 
Note that we should choose
the Keel's quadruple set coinciding with the one used in the algebraic reduction; 
then we will always get the same result as in the algebraic reduction.
This is an equivalent characterization of the algebraic linear reduction, in a graphical way; 
and actually much more efficient, compared to the algebraic reduction. 
This characterization naturally also leads to a complete algorithm that is equivalent to 
the whole process of the algebraic linear reduction until only clever monomials are in the 
expansion on right hand side of the equation. This naturally lead to another algorithm 
on computing the value of a given proper monomial. However, that algorithm is exponential in the 
complexity, way slower than the forest algorithm. We do not present the details of it in this paper. 
But the idea of tree version linear reduction is inspiring for us, when we want to compute 
the value of a ``sunlike tree'' which is one of the base cases of the forest algorithm correctness proof.

Now we know what is a linear reduction --- both algebraically and graphically. 
Actually, a potential reason why we do the linear reduction, is that 
we want to find a way to reduce the given monomial to a sum of clever monomials. Then, we can get the value 
of the given monomial directly. With this intention in mind, we want to know what is a ``good linear reduction'', namely
the one that reduce the power of some factor without giving another factor a higher power. If we have such a 
method, we can then guarantee that after finite number of steps, we reach a point where all factors have power one.
Recall that in the linear reduction, the only freedom we have is on the quadruple. We will introduce a 
criterion with the help of the tree representation of the given monomial, in the upcoming section.

\subsection{Good linear reduction}\label{subsec:good_linear_reduction}
Actually we can pick any quadruple $(i,j,k,l)$ to do the linear reduction, as long as $i,j\in I$ and $k,l\in ~J$.
However, we want to make sure that the quadruple we pick fulfills the {\em summand distinction property}.
\begin{definition}
We say that a quadruple $(i,j,k,l)$ fulfills the {\bf summand distinction property} w.r.t.~$M$ and~$\delta_{I,J}$.,
if no summand in $-(\epsilon_{ij\mid kl}-\delta_{I,J})$ divides~$M$, where~$M$ is a product of Keel's factors.
\end{definition}
We want this property for the chosen quadruple,
simply because it means any summand in $-(\epsilon_{ij\mid kl}-\delta_{I,J})$
is distinct from all factors of~$M$. In this way, we can do a ``good'' linear reduction --- what do we mean by ``good''
will be explained later. So why do we need this property on the quadruple?

Suppose that we want to replace one occurrence of~$\delta_{I_1,J_1}$ in 
$$M=\delta^p_{I_1,J_1}\cdots \delta_{I_q,J_q}$$ by 
$-(\epsilon_{ij\mid kl}-\delta_{I_1,J_1})$ and $(i,j,k,l)$ fulfills the summand distinction property.
 W.l.o.g., assume that 
 $$-(\epsilon_{ij\mid kl}-\delta_{I_1,J_1})=\sum_{s=1}^t{\delta_s}.$$
Then we obtain the following equation:
\begin{flalign}
 \delta^p_{I_1,J_1}\cdots \delta_{I_q,J_q} 
 &= \delta^{p-1}_{I_1,J_1}\cdots \delta_{I_q,J_q} \cdot (-(\epsilon_{ij\mid kl}-\delta_{I_1,J_1})) &&\\\nonumber
 &= \delta^{p-1}_{I_1,J_1}\cdots \delta_{I_q,J_q} \cdot \sum_{s=1}^t{\delta_s} &&\\\nonumber
 &= \delta^{p-1}_{I_1,J_1}\cdots \delta_{I_q,J_q} \cdot \delta_1 + \cdots + \delta^{p-1}_{I_1,J_1}\cdots \delta_{I_q,J_q} \cdot \delta_t &&
\end{flalign}
Now we focus on any monomial on right hand side of the equation, say $M_r:= \delta^{p-1}_{I_1,J_1}\cdots \delta_{I_q,J_q} \cdot \delta_r$.
Since~$\delta_r$ does not divide~$M$,~$\delta_r$ is distinct from all factors of~$M$. Therefore, no factor can have a higher power after 
the replacement. If we follow this rule for each linear reduction, then
 after finitely many times of linear reduction, we can obtain a (maybe negative) sum of clever monomials, 
which, in value, is equal to the given monomial. The analysis above also tells us that we should also obey this
rule when doing the algebraic reduction, otherwise we may not finish the reduction (until a sum of only clever monomials)
after finitely many steps. 

How can we guarantee that the chosen quadruple satisfies summand distinction property? 
By a {\em proper choice} of the Keel's linear quadruple.
In order to explain what is this ``proper choice'', we need the concept of {\em cluster} first.

\begin{definition}
Let $LT=(V,E,h,m)$ be a loaded tree with labeling set $N=\{1,\ldots,n\}$. 
We say that~$cl\subset N$ is a {\bf cluster of vertex~$v\in V$} if and only if one of the following two conditions holds:
\begin{enumerate}
\item $cl$ is a one-element-subset of~$h(v)$.
\item $cl$ is the collection of 
labels in one component of the graph obtained by removing vertex~$v$ and all its incident edges.
\end{enumerate}
Note that here~$h$ denotes the labeling function of the loaded tree to which~$v$ belongs.
If a cluster has cardinality one, we say that it is a {\bf singleton}; 
otherwise, we say that it is a {\bf proper cluster}.

\end{definition}
\begin{remark}
 We observe that collecting all parts of cuts of a loaded tree gives us exactly the collection
 of proper clusters of all vertices.
\end{remark}

\begin{remark}
 
It is not hard to check that the above two cases are disjoint for any loaded tree.
When a cluster of vertex~$v$ fulfills the first condition, it contributes one to the 
cardinality of~$h(v)$.
When a cluster of vertex~$v$ fulfills the second condition, it contributes one to the degree of~$v$.
Recall the expression $\deg(v)+|h(v)|$ in the third item of Definition~\ref{def:loaded_tree}:
For a vertex, each incident edge corresponds to a cluster of it, and each of its labels corresponds to a cluster of it as well.
\end{remark}

\begin{example}
For a better idea of this definition, let us see what are the clusters for vertices of the loaded tree in Figure~\ref{fig:loaded_trees}(b). 
To make the reading easier, we paste the figure also here.
 \begin{figure}
\centering
\includegraphics[width=0.5\linewidth]{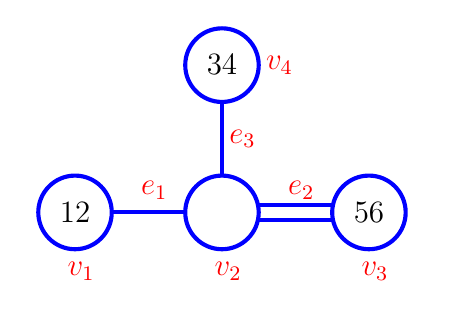}
\label{fig:m6_recall}
\end{figure}

\begin{itemize}
 \item Clusters for~$v_1$: $\{1\}$, $\{2\}$, $\{3,4,5,6\}$. 
 \subitem Singletons for~$v_1$: $\{1\}$, $\{2\}$. 
  \subitem      Proper clusters for~$v_1$: $\{3,4,5,6\}$.
 \item Clusters for~$v_2$: $\{1,2\}$, $\{3,4\}$, $\{5,6\}$. 
 \subitem Singletons for~$v_2$: none.
   \subitem     Proper clusters for~$v_2$: $\{1,2\}$, $\{3,4\}$, $\{5,6\}$.
 \item Clusters for~$v_3$: $\{1,2,3,4\}$, $\{5\}$, $\{6\}$.
 \subitem Singletons for~$v_3$: $\{5\}$, $\{6\}$.
    \subitem   Proper clusters for~$v_3$: $\{1,2,3,4\}$.
 \item Clusters for~$v_4$: $\{1,2,5,6\}$, $\{3\}$, $\{4\}$. 
       \subitem Singletons for~$v_4$: $\{3\}$, $\{4\}$.
      \subitem  Proper clusters for~$v_4$: $\{1,2,5,6\}$.
\end{itemize}
 \end{example}
 Now we are prepared for the concept of a ``proper choice'' of Keel's linear quadruple.

\begin{definition}
\begin{enumerate}
\item Assume w.l.o.g. that when we remove edge $e=\{v_1,v_2\}$, vertex 
$v_1$ is in the component where all labels collected to be~$I$ and~$v_2$ 
is in the component where all labels collected to be~$J$.  We call the corresponding components {\bf Component-$I$} and 
{\bf Component-$J$},
respectively.
\item We choose the quadruple $(i,j,k,l)$ such that $i,j\in I$ are from two distinct clusters of~$v_1$ and $k,l\in J$ are from
two distinct clusters
of~$v_2$. We call this way of choosing $i,j,k,l$ a {\bf proper choice}. 
\item And we call the corresponding quadruple
$(i,j,k,l)$ a {\bf proper quadruple} of the edge~$e$, or of the cut~$\{I,J\}$;
we call $\{i,j,k,l\}$ a {\bf proper quadruple set} of~$e$ or of~$\{I,J\}$. 
\end{enumerate}
Note that we are always able to choose a proper quadruple for any edge of some loaded tree, since apart from the cluster 
connected 
by edge~$e$ to~$v_2$ (or~$v_1$),~$v_1$ (or~$v_2$) has at least two more clusters, 
by the 
third condition of Definition~\ref{def:loaded_tree}. 
\end{definition}
It is not hard to see from the definition that we only talk about 
``proper quadruple'' when the given monomial is a tree monomial. To avoid confusion and unnecessary trouble, we only 
focus on non-clever tree monomials when we talk about proper quadruples.

\begin{example}
Let us continue with focusing on the loaded tree in Figure~\ref{fig:loaded_trees}(b).
For this loaded tree, suppose that we want to replace one occurrence of edge~$e_2$ (i.e., the cut $\{\{1,2,3,4\},\{5,6\}\}$). 
We should choose~$i,j$ from
$\{1,2,3,4\}$ and~$k,l$ from~$\{5,6\}$ for Keel's linear reduction. We see that $(1,3,5,6)$ is a proper choice, but 
neither $(1,2,5,6)$ nor $(3,4,5,6)$ is. 
\end{example}

Now let us look back on Example~\ref{eg:one_tree_generated} and Example~\ref{eg:two_reductions_needed}. Please verify that 
that we did choose the proper quadruples for both examples. We claim that any proper quadruple fulfills the summand distinction
property. For a better coordination of the structure, we postpone the proof to Section~\ref{sec:the_missing_proofs}.
\begin{proposition}\label{prop:summand_distinction_property}
Let~$M$ be a non-clever proper tree monomial and let~$\delta_{I,J}$ be the power-higher-than-one factor to be reduced. 
Let $Q=\{a,b,c,d\}$ be a proper quadruple set w.r.t.~$M$ and~$\delta_{I,J}$ such that $a,b\in I$ and $c,d\in J$. 
 Then, no summand in $-(\epsilon_{ab\mid cd}-\delta_{I,J})$ divides~$M$. 
 \end{proposition}
In the next section, we fill in the gaps of the proofs that are missing for this section.

\subsection{The missing proofs}\label{sec:the_missing_proofs}
In this section, we will settle down all the missing proofs for this section. After going through the previous context, we collect in total three
proof-missing statements:
\begin{enumerate}
 \item Proposition~\ref{prop:summand_distinction_property}: Any proper quadruple fulfills the summand distinction property (w.r.t. the ambient monomial and Keel's factor).
 \item Proposition~\ref{prop:two_vertices_weight_zero}: If a proper loaded tree has a multi-edge $e=\{v_1,v_2\}$ where the weights of~$v_1$,~$v_2$
 are both zero, then the tree has value zero.
 \item Theorem~\ref{thm:correctness_tree_version}: Correctness of the tree-version linear reduction algorithm (Algorithm~\ref{alg:tree_version_reduction}).
\end{enumerate}

Let us start with the first item listed above.

\begin{proof}[Proof of Proposition~\ref{prop:summand_distinction_property}:]
Let~$M:=\delta_{I_1,J_1}^{r_1}\cdots\delta^{r_t}_{I_t,J_t}$, where $I_1=I$, $J_1=J$ --- for the consistency of notations.
By Keel's linear relation (Fact~\ref{def:keel_linear_relation}),
$$\delta_{I_1,J_1}=\epsilon_{ac\mid bd}-(\epsilon_{ab\mid cd}-\delta_{I_1,J_1})$$ holds. 
By Remark~\ref{rem:sign}, we have 
$$M=\delta_{I_1,J_1}^{r_1-1}\cdots\delta^{r_t}_{I_t,J_t}\cdot(-(\epsilon_{ab\mid cd}-\delta_{I_1,J_1})).$$
W.l.o.g., assume that $\epsilon_{ab\mid cd}-\delta_{I_1,J_1}=\sum_{i=1}^k{\delta_i}$, then we have
$$M=-\delta_{I_1,J_1}^{r_1-1}\cdots \delta_{I_t,J_t}^{r_t}\cdot \delta_1-\cdots 
-\delta_{I_1,J_1}^{r_1-1}\cdots \delta_{I_t,J_t}^{r_t}\cdot \delta_k.$$
One observes that any summand on the right hand side of the above equation does not divide~$M$
if and only if~$\delta_i$ is distinct from any Keel's factor of~$M$, for any $1\leq i\leq k$. 
For simplicity, from now on we use~$I,J$ instead of~$I_1,J_1$. 
Denote by $T_M=(V,E,h,m)$ the tree of~$M$. Let $e=\{v_1,v_2\}$ be the corresponding multi-edge
of cut $\{I,J\}$ --- w.l.o.g. assume that~$v_1$ is in Component-$I$. Pick any~$\delta_i$ ($1\leq i \leq k$)
and denote by $\delta_{I',J'}:=\delta_i$, w.l.o.g. assume that $a,b\in I'$ and $c,d\in J'$.

First of all, since~$\delta_{I',J'}$ is a summand in $\epsilon_{ab\mid cd}-\delta_{I,J}$, it is clear by the definition
of~$\epsilon_{ab\mid cd}$ that $\delta_{I',J'}\neq \delta_{I,J}$. Then, because of the symmetry of~$a,b$ and~$c,d$ 
in the definition of a proper quadruple, it suffices to argue that $\delta_{I',J'}\neq \delta_{\tilde{I},\tilde{J}}$
for any proper choice of~$a,b$, where~$\delta_{\tilde{I},\tilde{J}}$ refers to the corresponding Keel's factor of any edge
of~$T_M$ in Component-$I$. There are three cases, since 
we choose~$a,b$ from two distinct proper branches of~$v_1$. It is not hard to visualize any edge~$\delta_{\tilde{I},\tilde{J}}$ in Component-$I$.

\begin{enumerate}
 \item[] Case 1). $a,b\in h(v_1)$. Now, we see that either $a,b,c,d\in \tilde{I}$, or $a,b,c,d\in\tilde{J}$, 
 where~$\delta_{\tilde{I},\tilde{J}}$ is the corresponding Keel's factor for any edge
 in Component-$I$; however, we know that $a,b\in I'$, $c,d\in J'$. Hence obviously $\delta_{I',J'}\neq \delta_{\tilde{I},\tilde{J}}$.
 \item[] Case 2). $a\in h(v_1)$, $b\notin h(v_1)$. In this case, we see that either $a,c,d\in \tilde{I}$, or $a,c,d\in\tilde{J}$, 
 where~$\delta_{\tilde{I},\tilde{J}}$ is the corresponding Keel's factor for any edge in Component-$I$. 
 For the same reason as in Case 1), we get $\delta_{I',J'}\neq \delta_{\tilde{I},\tilde{J}}$.
 \item[] Case 3). $a\notin h(v_1)$, $b\in h(v_1)$. Analogous to Case 2).
 \item[] Case 4). $a\notin h(v_1)$, $b\notin h(v_1)$. By an easy visualization of~$T_M$, we see that 
 $|\{a,b,c,d\}\cap \tilde{I}|\in\{1,3\}$, where~$\delta_{\tilde{I},\tilde{J}}$ is the corresponding Keel's factor for any edge in Component-$I$.
 For the same reason as in Case 1), we obtain that $\delta_{I',J'}\neq \delta_{\tilde{I},\tilde{J}}$.
\end{enumerate}
Hence,~$\delta_{I',J'}$ is distinct from any factor of~$M$. That is to say, $\{a,b,c,d\}$ fulfills the summand distinction property 
w.r.t.~$M$ and~$\delta_{I,J}$.
\end{proof}

From the above proof, we actually gained more information than needed.
Using the same notations as in the above proof, denote by
$$M_i:=\delta_{I_1,J_1}^{r_1-1}\cdots \delta_{I_t,J_t}^{r_t}\cdot\delta_i,\; 1\leq i\leq k.$$
Then we see that the set of labels of~$M_i$, i.e. $I_i\sqcup J_i$ for any $1\leq i\leq k$, 
is the same as that of~$M$. Also, the number of factors in~$M$ is the same as that in~$M_i$, 
$1\leq i\leq k$ --- both are~$\sum_{i=1}^t{r_i}$. So if we assume that~$M_i$ is a tree monomial, we get 
that~$T_M$ and~$T_{M_i}$ are of the same type. Observe that $\frac{M}{\delta_{I_1,J_1}}\cdot \delta_i=M_i$;
hence, comparing the trees~$T_M$ with~$T_{M_i}$ ($1\leq i\leq k$), we
can list the following differences:

\begin{enumerate}
 \item The multiplicity of edge~$e$ is reduced by one, since the power of~$\delta_{I_1,J_1}$ 
 is reduced by one.
 \item A new edge~$e'$ (corresponding to the factor~$\delta_i$ in~$M_i$) of multiplicity one is generated.
 \item Any other edge (except for~$e$) stays unchanged, in the sense that both its corresponding cut and its multiplicity remain unchanged.
\end{enumerate}
It seems that all the properties of~$T_{M_i}$ above would be true if~$T_{M_i}$ was a tree obtained from~$T_M$
after applying the tree-version linear reduction on edge~$e$ and quadruple $\{a,b,c,d\}$. Naturally, one raises 
the question: is it true that~$T_{M_i}$ is among the trees obtained after applying the tree-version linear 
reduction on~$e$ with $\{a,b,c,d\}$ on~$T_M$? Also, is it true conversely? Namely,
if a tree~$T$ can be obtained from~$T_M$ after the tree-version linear reduction with the corresponding parameters
--- meaning that it can be obtained from~$T_M$ when apply to it with edge~$e$ and quadruple $\{a,b,c,d\}$ after a vertex
splitting process --- then $T=T_{M_i}$ for some $1\leq i\leq k$. To sum up, this states the correctness of the 
tree-version linear reduction.

In the sequel, we prove the correctness of tree-version linear reduction algorithm. We
will keep the notations from the above analysis.
\begin{proof}[Proof of Theorem~\ref{thm:correctness_tree_version}:]
We prove the correctness of the tree-version linear reduction algorithm. We need to prove that any monomial whose corresponding tree 
is in the output set of the algorithm is a survival monomial in the ambient linear reduction process; while as any survival monomial
in the ambient linear reduction process is a monomial of a tree that is in the output of the tree-version linear reduction.

Let~$\tilde{T}$ be a tree among the trees in the output of the algorithm, for the input tree~$T_M$ with edge~$e$ and quadruple $\{a,b,c,d\}$. 
Then we see that~$M_{\tilde{T}}$, compared to~$M$: the factor~$\delta_{I_1,J_1}$ of~$\{I_1,J_1\}$ should get reduced by one, while all other factors stay unchanged, 
and one new factor~$\delta'$ that is a summand of~$\epsilon_{ab\mid cd}$ is added to the product, with power one. 
This shows that~$M_{\tilde{T}}$ is one summand on the right hand side of the equation, when we do a linear reduction 
using Keel's relation $\epsilon_{ab\mid cd}=\epsilon_{ac\mid bd}=\epsilon_{ad\mid bc}$ to reduce an occurrence of~$\delta_{I_1,J_1}$
in~$M$. 

And actually, depending on how we distribute the labels and arrange the branches of the split vertex to the two adjacent vertices 
of the new edge, we actually can manifest any survival monomial in the ambient linear reduction process, by considering the corresponding
monomial of~$\tilde{T}$. Because of the one-to-one correspondence (Lemma~\ref{lem:correspondence}), any monomial~$M_i$ ($1\leq i\leq k$)
is the corresponding monomial of a tree in the output of the algorithm, for the input tree~$T_M$ with edge~$e$ and quadruple $\{a,b,c,d\}$. 
\end{proof}

From the correctness of tree-version linear reduction algorithm, we naturally obtain the proof of Proposition~\ref{prop:two_vertices_weight_zero}.
\begin{proof}[Proof of Proposition~\ref{prop:two_vertices_weight_zero}:]
 In this case, when we want to replace one occurrence of the corresponding factor of~$e$ in the given monomial~$M_T$ to do the 
 linear reduction, we see that it is not possible to split any adjacent vertices of~$e$ via 
 the vertex splitting process (Algorithm~\ref{alg:vertex_splitting}). Hence we obtain NULL via Algorithm~\ref{alg:vertex_splitting}.
 By Theorem~\ref{thm:correctness_tree_version},~$M_T$ equals to zero
 in~$A^{\bullet}(\overline{\mathcal{M}}_{0,n})$. Therefore, the value~$M_T$ (or equivalently, of~$T$) is zero, since the zero element 
 of the ring is mapped to the zero element in~$\mathbb{Z}$ under the isomorphism.
\end{proof}

In the next section, we introduce ``sun-like trees''.

\section{Sun-like trees}\label{sec:sun_like_trees}
\begin{figure}[H]
\centering
\includegraphics[width=0.6\linewidth]{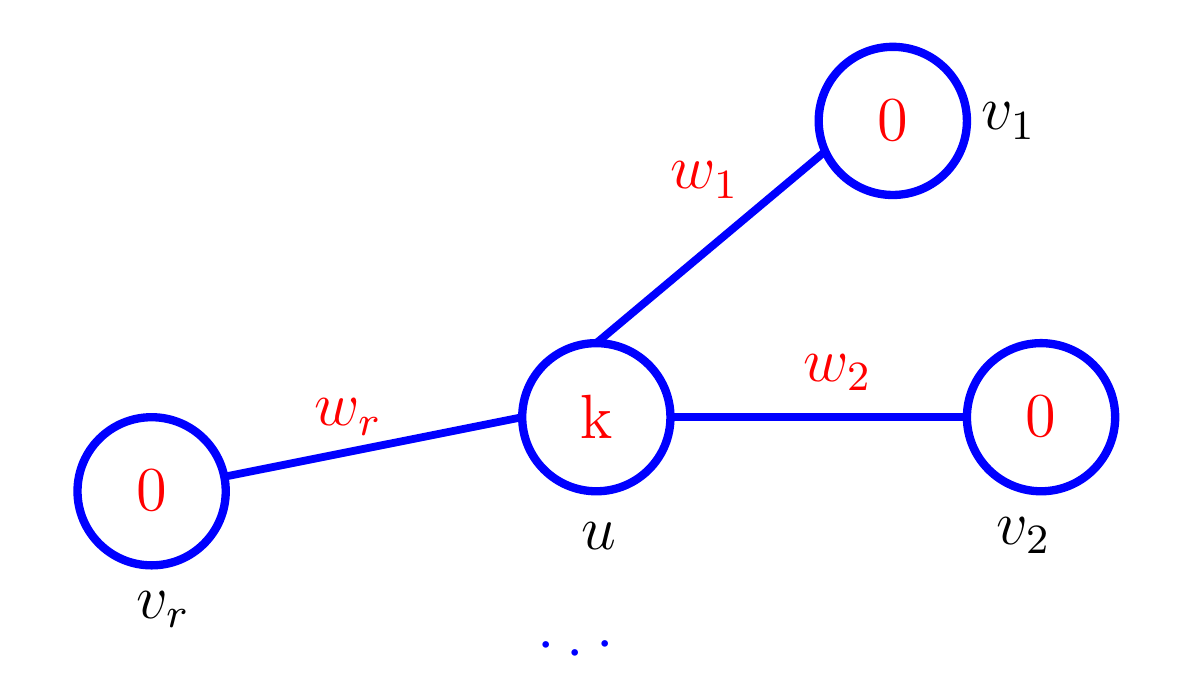}
\caption{The weighted tree of a sun-like tree, where the weights for vertices and edges
are marked in red. Note that~$w_i>0$ for all $1\leq i\leq r$ and that $k=\sum_{i=1}^r{w_i}$.}\label{fig:sun}
\end{figure}

In this section, we introduce a specific type of proper trees, and compute their values.
This serves as the base case for the correctness of forest algorithm.
\begin{definition}
We call a proper loaded tree {\bf sun-like} if it has domination number equal to one --- there exists a vertex~$u$
such that all other vertices are neighbors of it --- and all adjacent vertices of~$u$ have zero weight, all
edges have positive weights.

We call this vertex~$u$ the {\bf central vertex}. 
\end{definition}
Let~$k$ be the weight for the central vertex and $w_1, \ldots,w_r$ the weights for its incident edges, respectively.
By the weight identity for proper loaded trees, we know that $k=\sum_{i=1}^r{w_i}$. See Figure~\ref{fig:sun}
for a visualization. 

The main result for this chapter is the following theorem.
\begin{theorem}\label{thm:sun_like}
  Let $T=(V,E,h,m)$ be a sun-like proper loaded tree with~$u$ its central vertex. Let $e_1, \ldots, e_r$
  be the edges of~$T$, where $e_i=\{u,v_i\}$ for $1\leq i\leq r$. 
 Denote by~$w$ the weight function for~$T$, assume that $w(u)=k$ and $w(e_i)=w_i\geq 1$ for $1\leq i\leq r$.
 Then we have 
 $$\int(T)=(-1)^{\sum_{i=1}^r{w_i}}\cdot{k\choose w_1,\ldots,w_r}.$$ 
\end{theorem}

As a preparation before the proof of the above theorem, we need the ``edge-cutting lemma''.
We now focus on the {\em single-edge cutting} operation. 
\begin{definition}[single-edge cutting]
Let~$T$ be a loaded tree and $e=\{u,v\}$ be a single-edge of~$T$ with multiplicity~$r$ and corresponding 
factor~$\delta_{I_1,I_2}$. We construct two other loaded trees~$T_1$ 
and~$T_2$ by cutting off edge~$e$ and adding one more label~$x$ and~$y$ to~$u$ and~$v$, respectively.
This process is called {\bf single-edge cutting}.
\end{definition}
\begin{remark}
We see that~$T_1$ and~$T_2$ are still loaded trees and 
the weights of vertices~$u$ and~$v$ stay unchanged before and after the single-edge-cutting. 
\end{remark}
The following result tells us that $\int(T) = \int(T_1) \cdot \int(T_2)$.

\begin{proposition}[edge-cutting lemma]\label{prop:edge_cutting}
Let~$T$ be a loaded tree and~$T_1$ and~$T_2$ are the two trees obtained from a single-edge-cutting operation 
 executed on edge $e=\{v_1,v_2\}$ of~$T$.
 Then, we have that $\int(T) = \int(T_1) \cdot \int(T_2)$.
\end{proposition}

We postpone the proof of the edge-cutting lemma to Section~\ref{sec:from_algebra_to_geometry}.

\begin{corollary}\label{cor:unbalanced}
 Let~$T$ be a proper tree with a single edge~$e$ and~$T_1$,~$T_2$ be two loaded trees 
 obtained from~$T$ via a single-edge-cutting operation on edge~$e$. If~$T_1$ or~$T_2$ is not proper,
 then~$\int(T)=0$.
 \end{corollary}
\begin{proof}
Assume w.l.o.g. that~$T_1$ is improper. 
 Assume that~$T$ has~$n$ labels and~$n-3$ fringes,~$T_1$ has~$n_1$ labels and~$k_1$ 
 fringes, ~$T_2$ has~$n_2$ labels and~$k_2$ fringes. By construction, we obtain: $n_1+n_2=n+2$, $k_1+k_2=n-4$.
 Therefore, we obtain that $n_1+n_2=k_1+k_2-6$. Hence~$T_1$ is proper if and only
 if~$T_2$ is proper. Suppose w.l.o.g. that $k_1>n_1-3$. Then we know that~$T_1$ corresponds
 to a monomial in the Chow group~$A^{k_1}(n_1)$. From the fact stated in Section~\ref{sec:introduction} 
 we have that $A^{k_1}(n_1)=\{0\}$ and hence $\int(T_1)=0$. 
 By Proposition~\ref{prop:edge_cutting}, we obtain that 
 $$\int(T)=\int(T_1)\cdot \int(T_2)=0\cdot \int(T_2)=0.$$
\end{proof}
\begin{definition}
If both~$T_1$ and~$T_2$ in above stated process are proper, we say that~$T$
{\bf is balanced with respect to edge~$e$}.

A loaded tree is {\bf balanced} if and only if it is balanced with respect to 
any single edge. 
\end{definition}
If we do this operation on each 
single edge of loaded tree~$T$, then no matter in which sequence we choose to cut those
edges, we always obtain a same set of loaded trees.
 With the similar reasoning, the edge-cutting lemma has a generalized 
 version; so does Corollary~\ref{cor:unbalanced}. We conclude them in the following lemma,
 and we omit the proof.

\begin{lemma}\label{lem:n_new_trees}
 Let~$T$ be a loaded tree and $T_1,\ldots, T_n$ are the trees obtained from edge-cutting operations
 on all single edges of~$T$.
 Then, we have that 
 $$\int(T) = \int(T_1) \cdot \int(T_2)\cdots \int(T_n).$$ 
 And if~$T_i$ is not proper for any $1\leq i\leq n$, $\int(T)=0$.
\end{lemma}

The next result holds consequently.
\begin{theorem}
 Unbalanced loaded trees have value zero.
\end{theorem}
\begin{proof}
 By Lemma~\ref{lem:n_new_trees}, straightforward.
\end{proof}

Now let us look back on the vertex-splitting process. Whenever we split a vertex,
one single edge~$e'$ is generated in the structure~$S$. Therefore, in the tree-version linear reduction algorithm,
we can further filter the loaded trees in the output set --- those that are not balanced
w.r.t. edge~$e'$ can be directly removed from the set. 
\begin{definition}
We call those trees that remain in the output set after the above-described removal
the {\bf survival trees} (w.r.t. the input proper loaded tree and the chosen proper quadruple) 
and their corresponding monomials the {\bf survival monomials}. 
\end{definition}
Note that the balancing condition can be much more complicated if we consider it on the monomials alone
without the loaded-tree representation --- this is another argument for praising the tree representation,
an elegant language!

As another preparation before the proof of Theorem~\ref{thm:sun_like}, we introduce an identity on multinomial 
coefficients.
For any~$r$-many positive-integer parameters $m_1,m_2,\ldots, m_r$, define $s:=\sum_{i=1}^r{m_i}$. 
Denote a set of~$r$-many indeterminates as  
$$X:=\{x_1,x_2,\ldots,x_r\}.$$ 
Define $T:=\{B\mid B\subset X, x_1\in B\}$ and 
$$\mathcal{B}:=\{(B_1,B_2)\mid B_1\in T, B_2 = X\setminus B_1\}.$$
Define 
$$g: X\to \{m_1-1,m_2,\ldots, m_r\}$$
by $g(x_1):= m_1-1$, $g(x_i):= m_i$ for $i>1$.
For convenience in the later writing, we introduce the following notation.
Define for $B\subset X$, 
$$S(B):=\sum_{x\in B}{g(x)},\; \; 
{S(B) \choose B}:= \frac{S(B) !}{\prod_{x\in B}{(g(x)!)}}.$$ 
Based on the above preparation, the identity we want to 
introduce can be formulated as follows.
\begin{theorem}\label{thm:identity}
 $${s \choose {m_1,m_2,\ldots, m_r}}=\sum_{(B_1,B_2)\in
 \mathcal{B}}{{{s-r+1} \choose {S(B_2)-|B_2|}}\cdot{S(B_1) \choose B_1}\cdot{S(B_2) \choose B_2}},$$ 
 where~$|B_2|$ is the cardinality of~$B_2$.
\end{theorem}
The proof for this identity is postponed to Section~\ref{sec:proof_identity}, in order not to distract 
our main rhythm. In the next section, we provide the proof of Theorem~\ref{thm:sun_like} using the above identity.

\subsection{The value of a sun-like tree}\label{sec:proof_sun_like}
In this section, we give the proof of Theorem~\ref{thm:sun_like}.
First, the sign is taken care of by Remark~\ref{rem:sign}.
Hence we only need to verify that Theorem~\ref{thm:sun_like} computes the correct absolute value of the given sun-like tree.

\begin{itemize}
 \item{\bf Base case:} We prove by induction on~$k$. When $k=0$, since the tree is proper and all weights of the edges
 are positive, we know that in this case the tree has no edge. By Theorem~\ref{thm:clever}, we know that
 the tree has absolute value~$1$. Theorem~\ref{thm:sun_like} holds in this case.
 
 \item{\bf Proof idea for the general case --- apply Algorithm~\ref{alg:tree_version_reduction}:} 
 When $k\geq 1$, by Remark~\ref{rem:sign} and Algorithm~\ref{alg:tree_version_reduction},
 we know that the absolute value of~$T$ equals to 
 the sum of the values of the survival trees w.r.t.~$T$ and any proper quadruple.
 We hence choose a multi-edge and a proper quadruple, then consider what survival trees remain after the vertex splitting process and the balancing condition
 checking. Since in this process, one new weight-zero edge is generated
 and an old edge's weight is reduced by one, the edge weight sum is reduced by 
 one. For any tree in the output, we can cut-off the zero-weighted edge generated in the vertex-splitting process. 
 Then for the obtained two loaded trees, we can apply the induction hypothesis.
 
 \item{\bf Get prepared for the input data:} 
 Consider the sun-like tree $T=(V,E,h,m)$ with the weight function~$w$ in Figure~\ref{fig:sun}, denote by~$e_i$ the edge $\{u,v_i\}$. 
 We choose edge~$e_1$ to reduce. 
 In order to choose a proper quadruple, we should first
 figure out how many labels does the central vertex~$v$ have. 
 By the definition of the weight function, we have $w(u)=k=\deg(u)+|h(u)|-3$. Clearly $\deg(u)=r$ and we also have 
 $k=\sum_{i=1}^r{m_i}\geq r$ because of the weight identity for proper loaded trees. Hence $$|h(u)|=k-\deg(u)+3=k-r+3\geq 3.$$
 Hence~$u$ has at least three labels; denote by $c,d$ two of them.
 Since vertex~$v_i$ has degree one and weight zero for any $1\leq i\leq r$,
 it has two labels; denote them by $a_i,b_i$. We see that $\{c,d,a_1,b_1\}$
 is a proper quadruple set for the multi-edge~$e_1$. And obviously we can only split vertex~$u$, since~$v_1$ has zero weight. 
 We will use as input these data for Algorithm~\ref{alg:tree_version_reduction}.

 \item{\bf Arranging the labels --- part one:} 
 Suppose that the statement is true for all sun-like trees with the weight of the central vertex less or equal 
 to~$k-1$, where $k\geq 0$. Now,
 let us consider the case when the central vertex has weight~$k$. Let $T_1=(V_1,E_1,h_1,m_1)$ be any loaded 
 tree in the output (of Algorithm~\ref{alg:tree_version_reduction} applied to~$(T,e_1,\{c,d,a_1,b_1\})$). We split vertex~$u$ into~$u'$ and~$u''$, 
 denote by~$e'_1$ the edge $\{u',u''\}$. By Step 5. of vertex-splitting, 
 we know that $c,d\in h_1(u'')$. Hence, we have the freedom on the arrangements of all other labels of~$u$ to be 
 either the labels of~$u'$ or those of~$u''$. How many choices do we have? The answer will be revealed a bit later. 
 
 \item{\bf Arranging the branches:} We also have the freedom on arranging the branches. Denote by~$B_1$ the set of branches
 that will be attached to~$u'$ and by~$B_2$ for those that will be attached to~$u''$. 
 Denote by $x_1,\ldots,x_r$ the branches corresponding to $e_1,\ldots,e_r$, respectively.
 Let $X:=\{x_1,\ldots,x_r\}$, we see that $X=B_1\cup B_2$.
 By Step 6. of vertex-splitting, we know that $x_1\in B_1$, since this branch contains labels~$a_1$ and~$b_1$.
 For each bipartition of~$X$ into~$B_1$ and~$B_2$ such that $x_1\in B_1$, we 
 need to consider the distribution of labels
 so that the obtained tree is balanced w.r.t.~$e'_1$. 
 For this, we need to introduce several notations, so as to express the 
 arrangements of branches. 
 
% \item{\bf A formal model for the arrangements of branches:} 
% Define a function $g:X\to \{m_1-1,m_2,\ldots,m_r\}$ by $x_1\mapsto m_1-1$, $x_i\mapsto m_i$ for $2\leq i\leq r$. 
% Define $S(A):=\sum_{x\in A}{g(x)}$ and
% $${S(A) \choose A}:= \frac{S(A) !}{\prod_{x\in A}{(g(x)!)}}.$$
% Let $B_1$ be any arrangement of branches for vertex $u'$. Then $B_2=B\setminus B_1$. 
 
 \item{\bf Balancing condition considered:} 
 Since we require loaded tree~$T_1$ to be balanced with respect to edge~$e'_1$.
 After we cut-off the edge~$e'_1$, the tree containing~$u'$ should be proper; denote it by~$T'_1$. So does the tree 
 containing~$u''$, denote it by~$T''_1$ --- the weight identity should hold for both trees. From this, we obtain that 
the weight of~$u'$ in~$T'_1$ is~$S(B_1)$ and that the weight of~$u''$ in~$T''_1$ is~$S(B_2)$. 
Since the single-edge-cutting 
operation does not change the weight of vertex, we know that $w_1(u'')=S(B_2)$ and $w_1(u')=S(B_1)$ 
already hold in~$T_1$. 
And obviously, the degree of~$u''$ in~$T_1$ is $|B_2|+1$, while that of~$u'$ in~$T_1$ is $|B_1|+1$.

\item{\bf Arranging the labels --- part two:}
Now we can figure out how many labels should we distribute to~$u''$ (so as to guarantee the survival of~$T_1$). In~$T_1$:
$$w_1(u'')=\deg(u'')+|h_1(u'')|-3.$$
Therefore, we have: 
\begin{flalign*}
 |h_1(u'')|&=  w_1(u'')-\deg(u'')+3&&\\
                    &= S(B_2)-(|B_2|+1)+3&&\\
                    &= S(B_2)-|B_2|+2&&
\end{flalign*}
However, we already know that $c,d\in h_1(u'')$. Therefore, we should distribute $S(B_2)-|B_2|$ 
many labels of~$u$ to 
vertex~$u''$, so as to guarantee the survival of~$T_1$. Then naturally the remaining labels went to the labeling set of 
vertex~$u'$.

\item{\bf The induction step:} Then we cut off the single edge~$e'_1$ in~$T_1$, obtaining two 
proper trees~$T'_1$ and~$T''_1$, since~$T_1$ is balanced w.r.t.~$e'_1$. By Remark~\ref{rem:vertex_weight_non_zero}, we have
$w_1(u')+w_1(u'')=w(u)-1$. Hence $0\leq w_1(u'')<w(u)=k$ and $0\leq w_1(u')<w(u)=k$. Therefore, we can use 
the induction hypothesis on~$T'_1$ and~$T''_1$ if $m_1-1\neq 0$. 

\item{\bf Special case for the induction step:} 
When $m_1-1=0$, we can cut off edge~$e_1$ in~$T''_1$, then apply the induction hypothesis on the cut remainder that
contains~$u'$.
The other cut remainder is a single vertex with three labels (two of which are~$a_1$ and~$b_1$) --- this loaded tree 
has value one. Then by the edge-cutting lemma (Proposition~\ref{prop:edge_cutting}), we see that the value 
of~$T'_1$ equals to the value of the cut remainder that contains~$u'$. Then we can use the induction
hypothesis on this tree. Consider the following property of multinomial coefficients:
$${S\choose s_1,\ldots,s_p}={S\choose 0,s_1,\ldots, s_p}$$ 
for $S=\sum_{i=1}^p{s_i}$
and $s_i\in\mathbb{N}^+$. We see that our induction hypothesis can also apply in this case ---  the value of~$T'_1$
is not influenced by whether $m_1-1$ is zero or not.

\item{\bf Express the value of each loaded tree in the output:} For each given~$B_1$, among all the $k-r+3$ labels of~$u$,
$c,d$ are pre-fixed to belong to~$u''$. We should choose $S(B_2)-|B_2|$ many labels for~$u''$, 
from $k-r+1$ many labels
of~$u$. Then by the edge-cutting lemma, we know that for this arrangement,
the value of the obtained survival tree is the product of the values of two smaller (in the sense that they each 
has less weights on the 
central vertex) trees, after
cutting off the edge~$e'_1$ --- which is exactly ${S(B_2)\choose B_2}\cdot {S(B_1)\choose B_1}$.

\item{\bf Summing over these values:} 
That is to say, whenever the arrangement of the branches is fixed, because of the balancing condition requirement,
the label distribution is also fixed. There are ${k-r+1 \choose S(B_2)-|B_2|}$ many ways of 
label distributions. Since permutation or renaming the labels does not influence the value of the 
loaded tree, we know that for any two label distributions, the two trees have the same value: 
${S(B_2)\choose B_2}\cdot {S(B_1)\choose B_1}$ ---  the product of the values of 
the two cut remainders w.r.t.~$e'_1$.
Summing over the different arrangements of the branches, we obtain that the sum of the values of all the survival trees in the output 
is precisely:
$$\sum_{(B_1,B_2)\in \mathcal{B}}{{{k-r+1} \choose {S(B_2)-|B_2|}}
\cdot{S(B_1) \choose B_1}\cdot{S(B_2) \choose B_2}},$$
where~$\mathcal{B}$ represents the set of all pairs of bipartitions of~$X$ into~$B_1$,~$B_2$ such that $x_1\in B_1$.

\item{\bf Show-time for the identity:} 
Hence the only thing that is needed for our proof is the following identity:
$$\sum_{(B_1,B_2)\in \mathcal{B}}{{{k-r+1} \choose {S(B_2)-|B_2|}}
\cdot{S(B_1) \choose B_1}\cdot{S(B_2) \choose B_2}}={k \choose {m_1,\ldots, m_r}}.$$
With Theorem~\ref{thm:identity}, we conclude the proof.
\end{itemize}

\subsection{Equivalent characterization}\label{subsec:identity_equivalence}
In the last section, we see that we can prove Theorem~\ref{thm:sun_like}
using Theorem~\ref{thm:identity}.
In this section, we show that we can also prove the identity
given that Theorem~\ref{thm:sun_like} holds. 

Let us reconsider the proof steps in Section~\ref{sec:proof_sun_like}.
All the steps until the last do not depend on the identity, therefore
we can still use those analysis for the proof in this section. 
We keep using the notations from the last section, suppose that 
Theorem~\ref{thm:sun_like} holds. Then we know that 
if we apply Algorithm~\ref{alg:tree_version_reduction}
to~$T$ with multi-edge~$e_1$ and corresponding quadruple set $\{i,j,k,l\}$
same as in the last section, the update the output such that only balanced trees w.r.t. the newly generated edge~$e'$
remain in the output set. The sum of the values of output loaded trees
can be expressed as:
$$\sum_{(B_1,B_2)\in \mathcal{B}}{{{k-r+1} \choose {S(B_2)-|B_2|}}
\cdot{S(B_1) \choose B_1}\cdot{S(B_2) \choose B_2}},$$
where~$\mathcal{B}$ represents the set of all pairs of bipartitions of~$X$ into~$B_1$,~$B_2$ such that $x_1\in B_1$.

Correctness of Algorithm~\ref{alg:tree_version_reduction}
tells us that the above sum equals to~${k \choose m_1,\ldots,m_r}$.
That is to say, we have
$${k \choose {m_1,m_2,\ldots, m_r}}=\sum_{(B_1,B_2)\in \mathcal{B}}{{{k-r+1} \choose {S(B_2)-|B_2|}}
\cdot{S(B_1) \choose B_1}\cdot{S(B_2) \choose B_2}}.$$

Hence, Theorem~\ref{thm:sun_like} is an equivalent characterization of the identity.
Algorithm~\ref{alg:tree_version_reduction},
or the vertex-splitting process plays the essential role, during the process of proving their equivalence. We believe that the identity
indicates some complicated structural information of the vertex-splitting process, in an algebraic way.

Actually, when we choose the Keel's linear quadruple for the reduction, we have some freedom; since any proper quadruple will do the work. 
However, with this variant, we can gain two more identities on multinomial coefficients. 
The stand of these two identities follow from the equivalence proved above.

 \subsection{Two more identities}
 Let us keep the notations in Section~\ref{sec:proof_sun_like}.
Looking at the structure of a sun-like tree, we see that it is a generic choice to choose edge~$e_1$ to reduce.
Then, to choose the quadruple,~$a_1,b_1$ is the necessary and the only choice as the only two distinct clusters of~$v_1$.
However, among the clusters of~$u$, we actually have three distinct choices which then lead to three different identities on
multinomial coefficients:
\begin{enumerate}
 \item We can choose two labels from the labeling set of~$u$. This then leads to the same identity as described in Theorem~\ref{thm:identity},
 since this is the same case as we analyzed in Section~\ref{sec:proof_sun_like}.
 \item We can choose only one label from the labels of~$u$, and another label from some proper cluster of~$u$. 
 \item We can choose both two labels from two distinct proper clusters of~$u$. 
\end{enumerate}
In the second case, w.l.o.g. assume that we choose one label~$a_2$ from the branch~$x_2$ of~$u$, let~$c$ be the label we choose 
from the labels of~$u$. Then the quadruple we use is $\{c,a_2,a_1,b_1\}$. Let~$u''$ be the new label that contains~$c$, then 
we know that the branch~$x_2$ is also attached to~$u''$. Recall that $|h(u)|=k-r+3$ --- we can freely distribute all but~$c$ among 
the labels of~$u$, which are $k-r+2$ many. Recall that $|h(u'')|=S(B_2)-|B_2|+2$ --- but~$c$ is already fixed as a label of~$u''$;
we only need to distribute $S(B_2)-|B_2|+1$ many labels to~$u''$. The other part of the analysis stays unchanged. 
Hence by the correctness of tree-version linear reduction, we get that the absolute value of the given sun-like tree equals to:
$$\sum_{(B_1,B_2)\in \mathcal{B'}}{{{k-r+2} \choose {S(B_2)-|B_2|+1}}
\cdot{S(B_1) \choose B_1}\cdot{S(B_2) \choose B_2}},$$
where~$\mathcal{B'}$ represents the set of all pairs of bipartitions of~$X$ into~$B_1$,~$B_2$ such that $x_1\in B_1$ and $x_2\in B_2$.
This indicates that the central vertex of the given sun-like tree must have at least two degree. 
By Theorem~\ref{thm:sun_like}, we obtain that the absolute value of the given sun-like tree is~${k \choose {m_1,m_2,\ldots, m_r}}$.
Hence we have 
$${k \choose {m_1,m_2,\ldots, m_r}}=\sum_{(B_1,B_2)\in \mathcal{B'}}{{{k-r+2} \choose {S(B_2)-|B_2|+1}}
\cdot{S(B_1) \choose B_1}\cdot{S(B_2) \choose B_2}},$$
where $r\geq 2$.

In the third case, with analogous analysis as above, we obtain that 
$${k \choose {m_1,m_2,\ldots, m_r}}=\sum_{(B_1,B_2)\in \mathcal{B''}}{{{k-r+3} \choose {S(B_2)-|B_2|+2}}
\cdot{S(B_1) \choose B_1}\cdot{S(B_2) \choose B_2}},$$
where~$\mathcal{B''}$ represents the set of all pairs of bipartitions of~$X$ into~$B_1$,~$B_2$ such that $x_1\in B_1$ and $x_2,x_3\in B_2$, 
hence $r\geq 3$ is required.

With these two identities on multinomial coefficients, we conclude the chapter. 
In the next chapter, we present the correctness proof of the forest algorithm.

 \section{Correctness}\label{sec:correctness}
 In this section, we prove the correctness of the forest algorithm. As a preparation, we need to introduce 
 three types of edge-cutting on loaded trees.

\subsection{Three types of edge-cutting}
In this section, we introduce/describe three different types of edge-cutting on the loaded trees, and using these graphical operations,
express the main theorem in algebraic language.
In order not to interrupt the story later on, we need to introduce the concept of {\em star-cut} first.
 \begin{definition}[star]
 Let $T=(V,E)$ be a tree with $|V|\geq 3$. If there exists $v\in V$ such that all other vertices are neighbors of~$v$, 
 then we call this tree a {\bf star}.
\end{definition}
Consider a tree $T=(V,E)$. Pick any edge $e=\{u,v\}\in E$. We remove edge~$e$, attaching respectively a 
new vertex~$u_1$ to~$u$,~$v_1$ to~$v$ via a new edge~$e_1$ and~$e_2$.
 Then we obtain two new trees~$T_1$,~$T_2$. 
 \begin{definition}
 If we obtain at least one star after applying edge cut on some edge of~$T$, we 
call this edge cut a {\bf star-cut}.
\end{definition}

\begin{proposition}\label{prop:star-cut}
 Star-cut exists for any tree with no less than three vertices.
\end{proposition}
\begin{proof}
 Let~$T$ be an arbitrary tree with no less than three vertices. Define~$L$ to be the set of leaves of~$T$. 
 Then define an equivalence relation on~$L$ by
 $v_1\sim v_2$ if and only if $N(v_1)=N(v_2)$ for any two leaves $v_1,v_2$ of~$T$, where~$N(v)$ is the set of neighbors of vertex~$v$. 
 It is not hard to see that there is a 1-1 correspondence 
 between the set of support vertices, i.e., vertices that are adjacent to at least one leaf of~$T$ and the set of equivalence classes we define above. 
 
If~$T$ is a star, then the proposition holds
 since we can apply edge-cut to any edge of~$T$ and we will get a star after it. 
 Otherwise we delete all leaves of~$T$, obviously we obtain a nontrivial
 tree~$T_1$ --- here non-trivial means that it is not a single vertex. So~$T_1$
 must have an vertex~$u$ with degree~$1$, w.l.o.g., assume $\{u,u'\}\in E(T_1)$, then we also have $\{u,u'\}\in E(T)$. 
 Obviously,~$u$ is a support vertex of~$T$. If we apply edge-cut to edge $\{u,u'\}$ in~$T$, we get a star centered at 
 vertex~$u$ --- it is a star-cut.
 \end{proof}

\begin{figure}[H]
\centering
\includegraphics[width=0.5\linewidth]{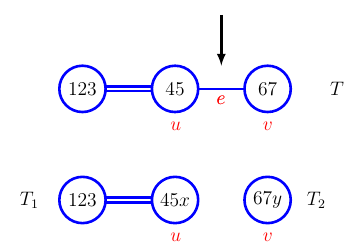}
\caption{The single-edge-cut operation applied to a tree~$T$ (above) w.r.t. 
edge~$e=\{u,v\}$, which gives two trees~$T_1$,~$T_2$ (below). Note that a new label is
added to~$u$ and~$v$ respectively.}\label{fig:single_cut_not_draw}
\end{figure}

Let~$T$ be a loaded tree and $e=\{u,v\}$ be any edge of~$T$ with multiplicity~$r$ and corresponding 
factor~$\delta_{I_1,I_2}$. When $r=1$, it is the single-edge cutting introduced in the paragraph before Proposition~\ref{prop:edge_cutting}.
Proposition~\ref{prop:edge_cutting} says us that $\int(T)=\int(T_1)\cdot \int(T_2)$.
See Figure~\ref{fig:single_cut_not_draw} for an illustration of such a cut. 

\begin{figure}[H]
\centering
\includegraphics[width=0.8\linewidth]{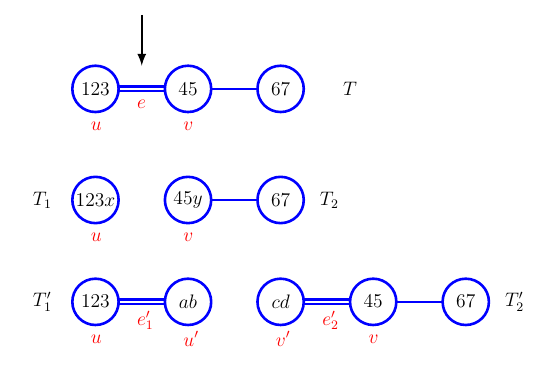}
\caption{The multi-edge-cut operation applied to a tree~$T$ (first row)
w.r.t. 
edge $e=\{u,v\}$, which gives two trees~$T'_1$,~$T'_2$ (third row).}\label{fig:multi_cut_not_draw}
\end{figure}

When $r>1$,
let~$s_1$ be the number of fringes of~$T_1$ and~$s_2$ be that of~$T_2$.
Let $N:=I_1\cup I_2$ and $|N|=n$. A quick calculation reveals that it can never happen that 
both~$T_1$ and~$T_2$ are proper, which indicates that~$\int(T)$ is always zero, if we want an analogous relation
as that of Proposition~\ref{prop:edge_cutting}.
Therefore, we need to modify our construction. We first construct~$T_1$ and~$T_2$ from~$T$ as in 
the single-edge cutting case. Then,
we construct~$T'_1$ (from~$T_1$) by removing the 
label~$x$ and attaching to~$u$ a new vertex~$u'$ via an edge~$e_1'$ connecting~$u$ and~$u'$, where the labeling set 
of~$u'$ is $\{a,b\}$ and the multiplicity of~$e_1'$ is set to be $|I_1|-s_1-1$; note that~$a$,~$b$ are 
two new labels not in~$N$. The construction
of~$T'_2$ via~$T_2$ is done analogously.
%In this way, we are able to obtain proper trees by adding new edges at the 
%end of the obtained trees, given that $T$ is a proper loaded tree.
\begin{definition}
Cutting off an edge of~$T$, obtaining~$T'_1$ and~$T'_2$ as stated above, is called a {\bf multi-edge cutting} operation.
\end{definition}
See Figure~\ref{fig:multi_cut_not_draw} for an illustration of a multi-edge cut. 
\begin{definition}
We call the two obtained trees after any of the 
three types of edge-cutting introduced in this sections the corresponding {\bf cut-remainders}.
\end{definition}

In this section, a main task for us is to investigate the relation between the value of~$T$ and the values
of~$T'_1$ and~$T'_2$. In the sequel, we will introduce more algebraic notations, so as to express better the corresponding theorem.

First, notice that~$T'_1$ and~$T'_2$ are both in many cases strictly smaller than~$T$; in those cases their monomials live in 
different ambient Chow rings than~$T$. For this, we introduce a foot index for the integral symbol,
indicating the ambient space. For instance, let~$M$ be a monomial in~$A^{n-3}(\overline{\mathcal{M}}_{0,n})$, then~$M$ is in the 
Chow ring of~$\overline{\mathcal{M}}_{0,n}$. Let~$N:=\{1,\ldots,n\}$ be the labeling set of~$\overline{\mathcal{M}}_{0,n}$, we 
denote by~$\int_{\mathcal{M}_N}(M)$ for the value of~$M$. 
From now on, we will sometimes use this notation, so as to clarify the ambient Chow ring and the labels in the ambient variety
as well. 
%Besides, in order to avoid clashes of notations,
%we denote by $A^{\bullet}(X)$ instead of $A^*(X)$ for the Chow ring of the variety $X$, from now on.

Let~$M_T$ be the monomial of~$T$, let~$\delta_{I_1,I_2}$ be the corresponding factor of the edge~$e$ and let $r\geq 1$ 
be the multiplicity of 
$e$; note that $I_1\cup I_2=N$. 
Let~$\mu_1$ be the product of factors of 
fringes in Component-$I_1$ and let~$\mu_2$ be that in Component-$I_2$ ---  
note that all factors inherit
the multiplicities of their corresponding edges via their powers. Since~$M_T$ is a tree 
monomial, we obtain the following conclusion.
Let~$S$ be the multi-set of the factors of~$M_T$, note that factors with power higher than one appear more than once in the set. 
Then~$\mu_1$ is the product of the generators $\delta_{U,V}\in S$ such that $U\subsetneq I_1$, 
and~$\mu_2$ is the product of the generators $\delta_{U,V}\in S$ such that $V\subsetneq I_2$.
Then $M_T=\mu_1\cdot \mu_2\cdot \delta^r_{I_1,I_2}$ and thence 
$$\int(T)=\int_{\mathcal{M}_N}{\mu_1\cdot \mu_2\cdot \delta^r_{I_1,I_2}}.$$
 Let~$s_1$,~$s_2$ be the degrees of~$\mu_1$ and~$\mu_2$, respectively; then $\mu_1\in A^{s_1}(\mathcal{M}_N)$, 
 $\mu_2\in A^{s_2}(\mathcal{M}_N)$. Because~$T$ is a proper loaded tree, we have $r=|N|-s_1-s_2-3$.

 Let~$T_1$,~$T_2$ be the two trees obtained from~$T$ by a single-edge cutting operation on the edge corresponding to 
 the cut~$\{I_1, I_2\}$; w.l.o.g., let~$T_1$,~$T_2$ have labels $I_1\cup \{x\}$, $I_2\cup\{y\}$ respectively. 
 This applies to the case when $r=1$.
 Let~$T'_1$,~$T'_2$ be the two trees obtained from~$T$ by a multi-edge cutting operation on the edge corresponding
 to the cut $\{I_1, I_2\}$; w.l.o.g., denote by~$a,b$ the two new labels added to~$T'_1$,~$T'_2$. This applies to 
 the case when $r>1$. In order to obtain the monomial for~$T_1$,~$T_2$,~$T'_1$ and~$T'_2$, we just need to consider the following
 operations.  Replace each factor~$\delta_{U,V}$ in~$\mu_1$ by~$\delta_{U, V\setminus I_2\cup \{x\}}$, denote 
 by~$\gamma_1$ the obtained monomial. Analogously, replace each factor~$\delta_{U,V}$ in~$\mu_2$ by~$\delta_{U, V\setminus I_1\cup \{y\}}$, 
 denote by~$\gamma_2$ the obtained monomial. 
 One can see that $M_{T_1}=\gamma_1\in A^{s_1}(\mathcal{M}_{I_1\cup\{x\}})$, $M_{T_2}=\gamma_2\in A^{s_2}(\mathcal{M}_{I_2\cup\{y\}})$. 
 Because of the new edges~$e'_1$,~$e'_2$ added during the multi-edge cutting, an extra factor is needed other 
 than~$\nu_1$ ($\nu_2$) for~$M_{T'_1}$ ($M_{T'_2}$). We have 
 $$M_{T'_1}=\nu_1\cdot (\delta_{I_1,\{a,b\}})^{|I_1|-s_1-1},\; M_{T'_2}=\nu_2\cdot (\delta_{I_2,\{a,b\}})^{|I_2|-s_2-2}.$$
 Then the following theorem reveals to us the relation between the value of~$T$ and the 
 values of~$T'_1$ and~$T'_2$. 
 \begin{theorem}\label{thm:reveal}
 With the notations above, the following equation holds.
  \begin{flalign*}
\int_{\mathcal{M}_N}{\mu_1\cdot \mu_2\cdot (\delta_{I_1,I_2})^r} =  
&{r-1 \choose |I_1|-s_1-2}\cdot \int_{\mathcal{M}_{I_1\cup\{a,b\}}}{\nu_1\cdot (\delta_{I_1,\{a,b\}})^{|I_1|-s_1-1}}&&\\
                      &\cdot \int_{\mathcal{M}_{I_2\cup\{a,b\}}}{\nu_2\cdot (\delta_{I_2,\{a,b\}})^{|I_2|-s_2-1}}.&&
\end{flalign*}
To say it in the expression of tree values, we have 
$$\int{T}={r-1 \choose |I_1|-s_1-2}\cdot \int{T'_1}\cdot \int{T'_2}.$$
 \end{theorem}
 As for the single-edge cutting case, we get the following theorem as a reformulation of Proposition~\ref{prop:edge_cutting}
 which says $\int{T}=\int{T_1}\cdot \int{T_2}$.
 \begin{theorem}\label{thm:edge-cutting}
$$\int_{\mathcal{M}_N}{\mu_1\cdot \mu_2\cdot \delta_{I_1,I_2}} =  
 \int_{\mathcal{M}_{I_1\cup\{x\}}}{\gamma_1}
                      \cdot \int_{\mathcal{M}_{I_2\cup\{y\}}}{\gamma_2}$$
 \end{theorem}
 \begin{remark}
  The above two theorems can be unified in the expression in the tree values: 
  $$\int{T}={r-1 \choose |I_1|-s_1-2}\cdot \int{T_1}\cdot \int{T_2},$$
 where~$T_1$,~$T_2$ are obtained after an edge-cut operation on~$T$. 
 In this sense, Theorem~\ref{thm:edge-cutting} can be viewed as a special
 case of Theorem~\ref{thm:reveal}.
  \end{remark}
 
 We postpone the proofs of Theorem~\ref{thm:edge-cutting} and Theorem~\ref{thm:reveal} to Section~\ref{sec:from_algebra_to_geometry}. 
 In the next part, we prove the correctness of forest algorithm with the help of 
 these two theorems.

 \subsection{Correctness of the forest algorithm}\label{subsec:correctness_forest_algorithm}
In this section, we prove the correctness of the forest algorithm. Basically, the correctness proof is to proof that 
$\int{(RF)}=\int{(LT)}$, where~$RF$ is the redundancy forest of the loaded tree~$LT$. 
First, we address the part that computes the absolute value. 
At the end of the section, we address the part of the algorithm that gives the sign. 

Before anything else, we want to address that one can check: removing any weight-zero vertex and its adjacent
edges of a redundancy tree~$RT$ does not influence the output value of the recursive formula. 
Inspired by this idea, we introduce the following definition.
\begin{definition}
We define the {\bf value/integral value of a redundancy tree}~$RT$ as: $\int{(RT)}=\int{(RF)}$, where~$RF$ is the 
redundancy forest obtained from~$RT$ by removing all weight-zero vertices and their adjacent edges.
\end{definition}
So, we can consider a modified 
forest algorithm, where we skip the ``weight-zero-vertices deleting'' step, directly apply the recursive formula on the redundancy tree 
of the given loaded tree. The correctness proof will be easier in some cases for this modified version, compare to the original version.
However, the original version of forest algorithm 
can be much more efficient in the practical computations.
Therefore, for the convenience of the correctness proof, we do not differ these two versions of forest algorithms. It is indicated 
in the context which one we are considering. Bearing this in mind, we come to the other parts of the proof.

We base our consideration on the loaded trees. We view sun-like alike trees (loaded trees that have the same shape/structure as sun-like trees)
as our base case. Given a loaded tree, we can cut off the single and multi-edges. By the ``star-cut'' proposition (Proposition~\ref{prop:star-cut}),
the definition of single and multi-edge cutting operations, we can repeat the edge cutting process until only sun-like alike 
trees or single vertices are left. Hence, these are the base cases we need to deal with.

When \textbf{\em the loaded tree has a single vertex with weight zero}: this is a clever tree, hence has value one. 
When the loaded tree has a single vertex of nonzero weight: it is improper, hence has value zero. 
It is not hard to imagine the corresponding redundancy trees of the given loaded trees in the above two cases,
from where we see that the forest algorithm gives the correct results, in these two cases. 

\begin{figure}
\centering
\includegraphics[width=0.6\linewidth]{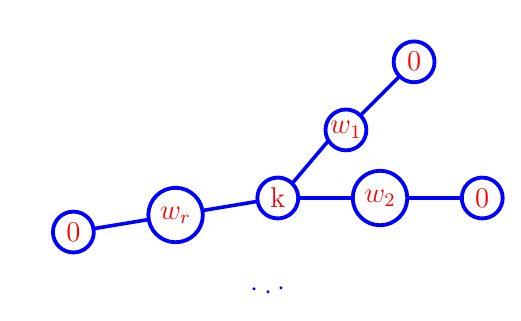}
\caption{This is the redundancy tree of a sunlike tree with weight~$k$ on the central vertex and weights $w_1,\ldots, w_r$ on the 
$r$ edges respectively, with weights marked in red.}\label{fig:rt_of_sunlike}
\end{figure}

When \textbf{\em the loaded tree~$LT$ is sun-like} with weight~$k$ on the central vertex and weights $w_1,\ldots, w_r$ on 
the~$r$ edges respectively: recall that sun-like trees are proper by definition, hence $k=\sum_{i=1}^r{w_i}$. 
Now consider the corresponding redundancy tree (depicted in Figure~\ref{fig:rt_of_sunlike}). By Theorem~\ref{thm:sun_like},
$|\int{(LT)}|={k \choose w_1,\ldots,w_r}$. From a basic property of multinomial coefficients, 
we have $${k\choose w_1,\ldots,w_r} = {k\choose w_1 }\cdot {k-w_1 \choose w_2}\cdots{w_r\choose w_r}.$$
Right hand side of the above equation is exactly what forest algorithm gives us. 
Hence the forest algorithm gives the correct result, in this case.

Consider \textbf{\em the case when an improper loaded tree~$LT$ is sun-like alike}, i.e. it has the same structure as a sun-like tree. 
Let~$k$ be the weight for the central vertex and let $w_1,\ldots, w_r$ be the weights for the~$r$ edges
of~$LT$. There are two sub-cases. First is when $k > \sum_{i=1}^r{w_i}$. According to the operations on the 
redundancy tree of~$LT$ described in the forest algorithm, in the end we will get a single degree-zero 
vertex with nonzero weight. Therefore the output from the forest algorithm in this case would be zero. 
Second is when $k < \sum_{i=1}^r{w_i}$. Then let us conduct the recursive formula on the adjacent 
vertices of the corresponding vertex~$v'$ of the central vertex~$v$ of~$LT$ in the redundancy tree~$RT$, sequentially. 
At some point of this process, vertex~$v'$ will have zero or negative weight. If it has zero weight, we can delete it,
which then leads to a situation of separate zero-degree vertices of nonzero weights. In either case, the forest algorithm
 outputs zero. Anyways, the loaded tree~$LT$ has zero value because of its improperness. Hence the forest algorithm
 is correct, in this case.
 
 So far, we dealt with all the base cases; in the sequel, we continue with the idea
 of proof by induction on the maximal number of edge-cuttings needed for the given loaded tree to reach 
 a level where all cut remainders belong to the base cases.

 Consider \textbf{\em the case when there exists a single-edge cut such that one part of the cut-remainders is an improper 
 sun-like tree or improper single vertex}. Then we see from the above analysis that the forest algorithm outputs 
 zero, while by Proposition~\ref{prop:edge_cutting} the loaded tree also has value zero.

 Now consider \textbf{\em the case when there exists a single-edge cut such that one of the cut-remainders~$LT_1$ is a weight-zero
 single vertex}. The corresponding redundancy tree~$RT$ (of the given loaded tree~$LT$) has a weight-zero leaf. The forest algorithm simply deletes this
 vertex and its adjacent edge. We then actually get the redundancy tree~$RT_2$ of the other cut-remainder~$LT_2$. 
 Output of the forest algorithm is~$\int{RT_2}$; while the value of~$LT$ (by Theorem~\ref{thm:edge-cutting}) is 
 $$\int{(LT_1)}\cdot\int{(LT_2)}=1\cdot \int{(LT_2)}=\int{(LT_2)}.$$ By induction, the forest algorithm is correct, in this case.

 Given a loaded tree~$LT$ with redundancy tree~$RT$, consider \textbf{\em the case when there exists a single-edge cut 
 such that one of the cut-remainders is a sun-like tree~$LT_1$}, and denote by~$LT_2$ the other cut-remainder. 
 Start the recursive formula on the leaves 
 of~$RT$ corresponding to the leaves of~$LT_1$. Let~$k$ be the weight of the central vertex and $w_1,\ldots,w_r$
be the weights of edges adjacent to the central vertex. The forest algorithm then will return 
$${k\choose w_1}\cdot{k-w_1 \choose w_2}\cdots{w_r\choose w_r}\cdot \int{(RT_2)},$$
where~$RT_2$ is the redundancy tree
of~$LT_2$. By induction, we have $\int{(RT_2)}=\int{(LT_2)}$. A property of multinomial coefficients says
$${k\choose w_1 }\cdot {k-w_1 \choose w_2}\cdots{w_r\choose w_r} = {k\choose w_1,\ldots,w_r}.$$ 
We then see that the forest algorithm returns ${k\choose w_1,\ldots,w_r}\cdot \int{(LT_2)}$, which by Theorem~\ref{thm:sun_like} is
$\int{(LT_1)}\cdot \int{(LT_2)}$, which by Theorem~\ref{thm:edge-cutting} is exactly the value of~$LT$.
Hence the forest algorithm is correct in this case.

So far, we have dealt with all cases where by a single-edge cut, we can reach one of the base cases as one of the 
cut-remainders. Now, we consider the remaining cases. 

First, we consider \textbf{\em the case when there is a leaf~$l$ with nonzero
weight whose unique incident edge $e=\{l,l_1\}$ is a multi-edge}. W.l.o.g., assume~$l_1$ is in Component-$I_1$
  and~$l$ is in Component-$I_2$ if we remove edge~$e$. Denote by~$s_1$ the number of fringes in Component-$I_1$ 
  and by~$s_2$ the number of fringes in Component-$I_2$. Now let us have a look at the formula in Theorem~\ref{thm:reveal}. 
  %Now we
 %cut off the multiple edge $e$, obtain two new trees $T'_1$ and $T'_2$; then we apply the formula in 
 %Theorem~\ref{thm:reveal}. 
 We see that, in this case, $r-1=m(e)-1=w(e)$. Since the given loaded tree~$LT$ is proper, we have $|I_1|+|h(l)|-3=m(e)+s_1$.
 By definition we have $|h(l)|+\deg(l)-3=w(l)$. Hence we get $|I_1|-s_1-2=m(e)-1-w(l)+\deg(l)-1=w(e)-w(l)+1-1=w(e)-w(l)$.
 Hence the binomial coefficient on the right hand side of the formula in Theorem~\ref{thm:reveal} is 
 $${w(e)\choose w(e)-w(l)}={w(e)\choose w(l)}.$$

 Note that $|I_2|=|h(l)|$, $s_2=0$ and $w(l)=|h(l)|+\deg(l)-3=|I_2|+1-3$, we obtain $|I_2|-s_2-1=w(l)+1$.
 Now we cut off the multi-edge~$e$, obtaining two new trees~$T'_1$,~$T'_2$. By the calculation, we 
 know that the new edge~$e'_2$
 incident to~$l$ in~$T'_2$ has multiplicity $w(l)+1$, hence its weight is~$w(l)$. Therefore~$T'_2$ is a proper loaded 
 tree with two vertices connected by an edge and the weights of two vertices are~$w(l)$ and~$0$, respectively. 
 Since $w(l)\neq 0$, the tree~$T'_2$ is a sun-like tree, by Theorem~\ref{thm:sun_like}, we know that its absolute
 value is ${w(l)\choose w(l)}=1$. Recall that the tree~$T'_1$ is obtained from~$LT$ by replacing vertex~$l$ by a weight-zero vertex, and 
 replacing edge~$e$ by an edge~$e'_1$ with multiplicity $w(e)-w(l)+1$ and hence the weight of~$e'_1$ is $w(e)-w(l)$.
 Theorem~\ref{thm:reveal} tells us that 
 $$|\int(LT)|={w(e)\choose w(l)}\cdot|\int(T'_1)|\cdot |\int(T'_2)|= {w(e)\choose w(l)}\cdot|\int(T'_1)|.$$
 Let~$RT_1$ be the redundancy tree of~$T'_1$. It is not hard to see that~$RT_1$ is exactly what we obtain 
 in the forest algorithm after removing the vertex corresponding to~$l$. Hence we see that the forest
 algorithm outputs ${w(e) \choose w(l)}\cdot |\int{(RT_1)}|$, which then by induction equals 
 ${w(e) \choose w(l)}\cdot |\int{(T'_1)}|$. This claims its correctness in this case.

% Notice that if some leaf has weight bigger than its parent vertex in the redundancy forest $RF$, then, when we cut off its unique incident edge
% in $LT$, the multiplicity of the edge $e'_1$ would be negative, which leads to the input loaded tree being value-zero.

The only case that is left is when all leaves have value zero and all edges are multi-edges. 
Now we do a multi-edge star cut on edge $e=\{v_1,v_2\}$ --- by Proposition~\ref{prop:star-cut} this operation is feasible.
Let~$LT$ be the given loaded tree. Let~$T'_1$ be the cut-remainder that is a sun-like 
tree with the central vertex~$v'_1$, and denote by~$T'_2$ the other cut-remainder with the vertex corresponding to~$v_2$
in~$LT$ denoted by~$v'_2$. Recall from how we define the multi-edge cut 
that~$T'_1$,~$T'_2$ are both proper. Denote by~$RT$ the redundancy tree of~$LT$, with the vertex corresponding to~$v_1$ in~$LT$ denoted 
by~$v''_1$. Now we delete all weight-zero leaves of~$RT$ that corresponds to the vertices in~$LT$ in the set 
$N(v_1)\setminus \{v_2\}$; denote by~$RT'$ the obtained redundancy tree with the vertex corresponding to~$v_1$ in~$LT$ 
denoted by~$v'$. Now we consider applying the recursive formula on the 
leaves adjacent to vertex~$v'$ sequentially.
 
By the fact that~$T'_1$ is proper, we know that when the forest algorithm already went through all 
leaves of~$v'$, the weight left in~$v'$ is $|I_1|-s_1-2$ (which is exactly the weight of the newly added edge incident to~$v'_1$ in the 
multi-edge cutting operation) --- recall that~$s_1$ is the number of fringes
of~$T_1$,~$I_1$ is the collection of labels of~$T_1$ and~$T_1$ is the loaded tree obtained in the intermediate
step of multi-edge cut applied to~$LT$. Hence, the binomial coefficient obtained in the next step in the 
forest algorithm is~${r-1 \choose |I_1|-s_1-2}$ (where~$r-1$ is the weight of the edge that is chosen to be cut), which coincides with the coefficient presented in 
Theorem~\ref{thm:reveal}. Note that here we consider the forest algorithm on~$LT$. 
Then, the forest algorithm would continue with popping up the binomial coefficient
$${w(v_2)\choose (r-1)-(|I_1|-s_1-2)},$$
where~$w(v_2)$ refers to the weight of~$v_2$ in~$LT$. Since~$LT$ is proper, we have that 
$|I_1|+|I_2|-3=s_1+s_2+r$; hence we have $(r-1)-(|I_1|-s_1-2)=|I_2|-s_2-2$. 
We see that $|I_2|-s_2-2$ is exactly the weight of the newly added edge incident to~$v'_2$ in~$T'_2$. 
Then, with the property of multinomial coefficients that 
$${k\choose w_1 }\cdot {k-w_1 \choose w_2}\cdots{w_r\choose w_r} = {k\choose w_1,\ldots,w_r},$$
and Theorem~\ref{thm:reveal}, it is not hard to see that the forest algorithm is correct in this case.

Before we conclude the correctness of the forest algorithm, there is still one matter that we need to address: the sign.
Assume~$T$ is a proper loaded tree with non-zero integral value. Let~$T_1$,~$T_2$ be the two trees obtained after 
applying a single-edge-cutting operation on~$T$; let~$T'_1$,~$T'_2$ be the two trees obtained after applying a 
multi-edge-cutting operation on~$T$. Then we know that all~$T_1$,~$T_2$,~$T'_1$,~$T'_2$ are proper loaded trees,
and an easy calculation tells us that the sum of edge/vertex weight sum of~$T_1$ and that of~$T_2$ equals that of~$T$; so 
does that of~$T'_1$ and~$T'_2$. Then, by Theorem~\ref{thm:sun_like}, we know that the sign of a given (proper) loaded 
tree is indeed~$-1$ to the power of its edge/vertex weight sum. 
 
 \subsection{From algebra to geometry}\label{sec:from_algebra_to_geometry}
 
 In this section, we prove Theorem~\ref{thm:edge-cutting} and Theorem~\ref{thm:reveal} which indicate
 the main geometric structure hidden beneath the forest algorithm,
 however using pure algebra. Bolded lines are serving as the indication of a road map, to help readers better understand the proof story line.
 
 We introduce an equivalence relation~$\sim$ in the Chow ring~$A^{\bullet}(\overline{\mathcal{M}}_{0,n})$ as follows. 
 For $a,b\in A^k(\overline{\mathcal{M}}_{0,n})$, we say that 
 $a\sim b$ if and only if $\int{c\cdot (a-b)}=0$ holds for any $c\in A^{n-3-k}(\overline{\mathcal{M}}_{0,n})$. 
 Then the elements that are in 
 the equivalence class of~$0$ form an ideal, denoted by~$I$. 
 \begin{definition}
The quotient ring~$A^{\bullet}(\overline{\mathcal{M}}_{0,n})/I$ 
 is called the {\bf numerical Chow ring}
 of~$A^{\bullet}(\overline{\mathcal{M}}_{0,n})$, denoted by~$A^{\bullet}_{num}(\overline{\mathcal{M}}_{0,n})$.
 \end{definition}
 Since in our ambient moduli space~$\overline{\mathcal{M}}_{0,n}$, the Chow equivalence and the numerical equivalence 
 are the same, we can conduct our proof in the numerical 
 Chow ring of~$\overline{\mathcal{M}}_{0,n}$. Throughout this section, 
 we consider the numerical Chow ring~$A^{\bullet}_{num}(\overline{\mathcal{M}}_{0,n})$, instead of the previous 
 ring~$A^{\bullet}(\overline{\mathcal{M}}_{0,n})$.
  We need this view angle on the ambient ring in order to use the following 
  result:
  $$A^1_{num}(X\times Y)\cong A^1_{num}(X)\bigoplus A^1_{num}(Y)$$ holds for any two smooth subvarieties~$X$,~$Y$ of~$\overline{\mathcal{M}}_{0,n}$.

 \textbf{\em We need to get familiar with the concepts of pushforward and pullback maps, so as to conduct the proof of Theorem~\ref{thm:edge-cutting}.}
 Let $f:X\to Y$ be a proper map between two smooth projective varieties. Then~$f$ induces the pushforward map
 $f_*:A^{\bullet}_{num}(X)\to A^{\bullet}_{num}(Y)$, which is a group homomorphism, and the pullback map 
 $f^*:A^{\bullet}_{num}(Y)\to A^{\bullet}_{num}(X)$, which is a ring homomorphism that preserves the degree of the ambient group where the 
 element lives. Let $\alpha\in A^{\bullet}_{num}(X)$, $\beta\in A^{\bullet}_{num}(Y)$, then the following adjoint formula on the integrals holds:
 $$\int_X(\alpha\cdot f^*(\beta)) = \int_Y(f_*(\alpha)\cdot \beta).$$ 
 
 Now we can come to the proof of Theorem~\ref{thm:edge-cutting}.
 
 \begin{proof}[Proof of Theorem~\ref{thm:edge-cutting}.]
 Let~$\delta_{I'_1,I'_2}$ be a factor of~$\mu_1$, then $I'_1\subsetneq I_1$ and $I_2\subsetneq I'_2$ hold.
 Then we can replace~$I_2$ by the symbol~$x$ in the monomial~$\delta_{I'_1,I'_2}$ and obtain a factor of~$\gamma_1$
 which is in the Chow ring~$A^{\bullet}_{num}(\mathcal{M}_{I_1\cup \{x\}})$. All factors of~$\mu_1$ correspond to elements
 in~$A^1_{num}(\mathcal{M}_{I_1\cup \{x\}})$.
 The product of these factors corresponds to $\gamma_1\in A^{\bullet}_{num}(\mathcal{M}_{I_1\cup\{x\}})$.
  By~\cite[Fact 2.]{keel}, we have 
  $$D_{I_1,I_2}\cong \mathcal{M}_{I_1\cup\{x\}}\times \mathcal{M}_{I_2\cup\{y\}}.$$
  Let $p_1:D_{I_1,I_2}\to \mathcal{M}_{I_1\cup \{x\}}$ be the natural projection map; 
  analogously we have $p_2:D_{I_1,I_2}\to \mathcal{M}_{I_2\cup \{y\}}$. 
   Then~$p_1^*(\gamma_1)$ is the product of all factors depending on the first group 
   of variables in $\mathcal{M}_{I_1\cup\{x\}}\times \mathcal{M}_{I_2\cup\{y\}}$, which is just~$\gamma_1$;
   the situation is analogous for~$p_2^*(\gamma_2)$.
  Denote by~$i$ the embedding 
  of~$D_{I_1,I_2}$ as a hypersurface into~$\mathcal{M}_N$. Then we have 
  $$\mu_1\cdot \mu_2\cdot \delta_{I_1,I_2}=i_*(p_1^*(\gamma_1)\cdot p_2^*(\gamma_2)).$$
  Hence we have
 \begin{flalign*}
\int_{\mathcal{M}_N}{\mu_1\cdot \mu_2\cdot \delta_{I_1,I_2}\cdot 1} &=  
 \int_{\mathcal{M}_N}{i_*(p_1^*(\gamma_1)\cdot p_2^*(\gamma_2))\cdot 1}&&\\
                      &=\int_{D_{I_1,I_2}}{p_1^*(\gamma_1)\cdot p_2^*(\gamma_2)\cdot i^*(1)}&&\\
                      &=\int_{D_{I_1,I_2}}{p_1^*(\gamma_1)\cdot p_2^*(\gamma_2)\cdot 1}&&\\
                      &= \int_{\mathcal{M}_{I_1\cup\{x\}}}{\gamma_1}\cdot \int_{\mathcal{M}_{I_2\cup\{y\}}}{\gamma_2}.&&
\end{flalign*}
 \end{proof}
 
 From now on, we try to prove Theorem~\ref{thm:reveal}. 
 \textbf{\em In order to prove Theorem~\ref{thm:reveal}, we need to explore deeper in~$\mathcal{M}_N$.}
 
 It is known that $D_{I_1,I_2}\cong \mathcal{M}_{I_1\cup\{x\}}\times \mathcal{M}_{I_2\cup\{y\}}$, where $x,y\notin N$ are new labels.
 Denote by~$p_{I_1,x,N}$ the projection from~$D_{I_1,I_2}$ to~$\mathcal{M}_{I_1\cup\{x\}}$
 and by~$p_{I_2,x,N}$ the projection from~$D_{I_1,I_2}$ to~$\mathcal{M}_{I_2\cup\{x\}}$.
Denote by~$i_{I_1,I_2}$ the embedding of~$D_{I_1,I_2}$ as a hypersurface into~$\mathcal{M}_N$.
 Let $i:=i_{I_1,I_2}$ and let $p_1:=p_{I_1,x,N}$, $p_2:=p_{I_2,y,N}$.
 We have $i^*(\mu_1\cdot \mu_2)=p^*_1(\gamma_1)\cdot p^*_2(\gamma_2)$. 
 Then apply the pushforward map on both sides. We obtain 
 $$\mu_1\cdot \mu_2\cdot\delta_{I_1,I_2}=i_*(p_1^*(\gamma_1)\cdot p_2^*(\gamma_2));$$
 this equation will be used in the proof of Theorem~\ref{thm:reveal}.
 We have:
 \begin{flalign*}
\int_{\mathcal{M}_N}{\mu_1\cdot \mu_2\cdot (\delta_{I_1,I_2})^r} &=  
\int_{\mathcal{M}_N}{(\delta_{I_1,I_2})^{r-1}\cdot i_*(p_1^*(\gamma_1)\cdot p_2^*(\gamma_2)) }&&\\
                      &=\int_{D_{I_1,I_2}}{i^*((\delta_{I_1,I_2})^{r-1})\cdot p_1^*(\gamma_1)\cdot p_2^*(\gamma_2)}&&\\
                      &=\int_{D_{I_1,I_2}}{(i^*(\delta_{I_1,I_2}))^{r-1}\cdot p_1^*(\gamma_1)\cdot p_2^*(\gamma_2)}.&&
\end{flalign*}
 \textbf{\em In the next step, we want to express the term~$i^*(\delta_{I_1,I_2})$ as a sum of two terms, so as to get the 
 coefficient~${r-1 \choose |I_1|-s_1-2}$. To do so, we introduce two new notations~$\beta_1$ and~$\beta_2$.}

Define 
$$\beta_1 = \beta_{x,I_1\cup\{x\}}:=[(p^{-1}_{I_1,x,I_1\cup \{a,b\}})^*\circ i^*_{I_1,\{a,b\}}](\delta_{I_1,\{a,b\}}).$$
Note that~$p_{I_1,x,I_1\cup\{a,b\}}$ is the projection from~$D_{I_1,\{a,b\}}$ to~$\mathcal{M}_{I_1\cup \{x\}}$. 
Recall that~$i_{I_1,\{a,b\}}$ is the embedding of~$D_{I_1,\{a,b\}}$ to~$\mathcal{M}_{I_1\cup\{a,b\}}$.
Consider the isomorphism $D_{I_1,\{a,b\}}\cong \mathcal{M}_{I_1\cup\{x\}}\times \mathcal{M}_{\{a,b,x\}}$; 
since~$\mathcal{M}_{\{a,b,x\}}$ is just a point, we obtain that $D_{I_1,\{a,b\}}\cong \mathcal{M}_{I_1\cup \{x\}}$.
Therefore the inverse of~$p_{I_1,x,I_1\cup\{a,b\}}$ exists. Since pullback is degree-preserving and 
$\delta_{I_1,\{a,b\}}\in A^1_{num}(\mathcal{M}_{I_1\cup\{a,b\}})$, we know that $\beta_{x,I_1\cup\{x\}}\in A^1_{num}(\mathcal{M}_{I_1\cup\{x\}})$. 
Analogously, we define 
$$\beta_2 = \beta_{x,I_2\cup\{x\}}:=[(p^{-1}_{I_2,x,I_2\cup \{a,b\}})^*\circ i^*_{I_2,\{a,b\}}](\delta_{I_2,\{a,b\}}) \in A^1_{num}(\mathcal{M}_{I_2\cup\{x\}}),$$
simply by 
replacing~$I_1$ by~$I_2$, in the definition of~$\beta_{x,I_1\cup\{x\}}$. 
Now we can express~$i^*(\delta_{I_1,I_2})$ as a sum of two summands in the following lemma.
\begin{lemma}\label{lem:direct_product}
  The following equation holds:
 $$i^*(\delta_{I_1,I_2})=p^*_1(\beta_1)+p^*_2(\beta_2).$$
\end{lemma}

\textbf{\em In order to prove the above lemma, we need to introduce some basic properties of the Chow group of 
a direct product of two varieties in~$\mathcal{M}_N$.} 
Let~$X$ and~$Y$ be two smooth projective subvarieties of~$\mathcal{M}_N$, then we have 
$$A^1_{num}(X\times Y)\cong A^1_{num}(X)\bigoplus A^1_{num}(Y).$$
In the sequel, we briefly give the references for the above formula. 
By~\cite{algebraic_cobordism}, the Chow ring modulo numerical equivalence is isomorphic to 
the ring of algebraic cobordism. By the same paper, algebraic cobordism is isomorphic to the cohomology ring.
The cohomology ring is equivalent to the homology, because of the Poincare duality for compact varieties (see~\cite{poincare_duality}).
Then, by the K$\ddot{u}$nneth formula (\cite{kuenneth_1},~\cite{kuenneth_2}) for the homology of a product of varieties, we see that the above formula 
holds.

Let~$\pi_l$,~$\pi_r$ be the projection from 
$X\times Y$ to~$X$ and~$Y$, respectively. We know that for any $y_0\in Y$, there exists a right 
inverse~$\sigma_l$ of~$\pi_l$ such that $\sigma_l(x):=(x,y_0)$. Let~$\sigma_l$ be any such inverse; the choice 
of the element in~$Y$ does not matter; we define~$\sigma_r$ analogously, as a right inverse for~$\pi_r$.
By specializing the K\"{u}nneth theorem to degree one, we have 
$$t=\pi^*_l\circ \sigma_l^*(t)+\pi^*_r\circ \sigma_r^*(t)$$
for any $t\in A^1_{num}(X\times Y)$.
Observe that we have the isomorphism 
$$D_{I_1,I_2}\cong \mathcal{M}_{I_1\cup\{x\}}\times \mathcal{M}_{I_2\cup\{x\}}$$
in~$\mathcal{M}_N$. Let~$q_{I_i,x,N}$ be any right inverse (as described above) of~$p_{I_i, x,N}$ for $i=1,2$. 
Denote by $q_i:=q_{I_i,x,N}$, for $i=1,2$.
Then, from the above analysis, we know that for any $a\in D_{I_1,I_2}$, we have 
$$a=p_1^*\circ q_1^*(a)+p_2^*\circ q_2^*(a).$$ 
Let $\alpha\in A^{\bullet}_{num}(X)$, $\beta\in A^{\bullet}_{num}(Y)$.
By the general properties of the degree (general fact of a proper map),
we have 
$$\int_{X\times Y}{\pi_l^*(\alpha)\cdot \pi_r^*(\beta)} = \int_X{\pi_{1*}(\pi_1^*(\alpha)\cdot \pi_2^*(\beta))}.$$
Then, because of the projection formula (\cite[Chapter 8]{Fulton}), we obtain that
$$\int_X{\pi_{1*}(\pi_1^*(\alpha)\cdot \pi_2^*(\beta))}=\int_X{\alpha\cdot \pi_{1*}(\pi_2^*(\beta))}.$$
The right hand side of the above equation then equals the following items:
$$\int_X{\alpha\cdot (\int_Y{\beta\cdot [X]})}=\int_Y{\beta\cdot \int_X{(\alpha\cdot [X])}}=\int_Y{\beta}\cdot \int_X{\alpha}.$$
Hence we have 
$$\int_{X\times Y}{\pi_l^*(\alpha)\cdot \pi_r^*(\beta)} = \int_X{\alpha}\cdot \int_Y{\beta}.$$

We need some more preparation before proving Lemma~\ref{lem:direct_product}.
Define $s_{k,l,N}:\mathcal{M}_{N\setminus \{k\}}\to \mathcal{M}_N$
as $$s_{k,l,N}:=i_{\{k,l\},N\setminus \{k,l\}}\circ p^{-1}_{N\setminus\{k,l\},l,N},$$ 
where~$k,l$ are two distinct labels of~$N$. 
Note that~$p_{N\setminus\{k,l\},l,N}$ is an isomorphism, hence it has an inverse. 
There is a surjective forgetful map $c_{a,N}:\mathcal{M}_N\to \mathcal{M}_{N\setminus\{a\}}$
for any $a\in N$. The above defined map~$s_{k,l,N}$ is a right inverse of~$c_{k,N}$. The image 
of~$s_{k,l,N}$ is the hypersurface~$D_{\{k,l\},N\setminus \{k,l\}}$ in~$\mathcal{M}_N$.

\begin{proof}[Proof of Lemma~\ref{lem:direct_product}]
Recall that~$i^*$ is the pullback map from~$A^{\bullet}_{num}(\mathcal{M}_N)$
 to~$A^{\bullet}_{num}(D_{I_1,I_2})$ and that $\delta_{I_1,I_2}\in A^1_{num}(\mathcal{M}_N)$.
 Since pullback is a degree-preserving ring homomorphism, we know that
 $i^*(\delta_{I_1,I_2})\in A^1_{num}(D_{I_1,I_2})$. Using the result from earlier
 analysis, we have Equation (a):
 $$i^*(\delta_{I_1,I_2})=p^*_1\circ q_1^*(i^*(\delta_{I_1,I_2}))+ p^*_2\circ q_2^*(i^*(\delta_{I_1,I_2})).$$
 We claim that \textbf{\em it suffices to prove 
 Equation (b): 
 $$q^*_1\circ i^*(\delta_{I_1,I_2})=\beta_1$$ 
 and Equation (c): 
 $$q^*_2\circ i^*(\delta_{I_1,I_2})=\beta_2.$$} 
 Suppose they hold,
 then from (b) we have: 
 $$p_1^*\circ q_1^*\circ i^*(\delta_{I_1,I_2})=p_1^*(\beta_1).$$ 
 Analogously, we obtain 
 $$(q_2\circ p_2)^*(i^*(\delta_{I_1,I_2}))=p^*_2(\beta_2)$$
 from (c). Substituting the equalities back to Equation (a), we obtain the wanted 
 equation. Since (b) and (c) are symmetric, it suffices to prove (b).
 
 \textbf{\em We prove Equation (b) by induction on~$|I_2|$.} Recall the definition of~$\beta_1$, in the base case, we have:
 $\beta_1=(p_1^{-1})^*\circ i^*(\delta_{I_1,I_2})$. It suffices to show that $(p_1^{-1})^*=q_1^*$. 
 Since~$q_1$ is a right inverse of~$p_1$, we have $p_1\circ q_1=id$.
 But in this case~$p_1$ is an isomorphism, so does~$q_1$; since $|I_2|=2$. Therefore, $q_1=p_1^{-1}$. Hence $(p_1^{-1})^*=q_1^*$.
 Suppose Equation (b) holds when $|I_2|=z-1$, now $|I_2|=z\geq 3$, $z\in \mathbb{N}$.
 Let $k,l\in I_2$ be two distinct labels and let $I'_2:=I_2\setminus \{k\}$, $N':=N\setminus \{k\}=I_1\cup I'_2$.
 Then we have the following equality: 
 $$i\circ q_1=s_{k,l,N}\circ i_{I_1,I'_2}\circ q_{I_1,x,N'}$$ of maps 
 from~$\mathcal{M}_{I_1\cup \{x\}}$ to~$\mathcal{M}_N$. Since pullback is a contravariant functor and
 $\delta_{I_1,I_2}\in A^{\bullet}_{num}(\mathcal{M}_N)$, we obtain:
    \begin{flalign*}
q^*_1\circ i^*(\delta_{I_1,I_2}) &=  
q^*_{I_1,x,N'}\circ i^*_{I_1,I'_2}\circ s^*_{k,l,N}(\delta_{I_1,I_2})&&\\
                      &=q^*_{I_1,x,N'}\circ i^*_{I_1,I'_2}(\delta_{I_1,I'_2}).&&
\end{flalign*}
 Now we can use the induction hypothesis,
 since $|I'_2|=|I_2|-1=z-1$. Hence we have $q^*_1\circ i^*(\delta_{I_1,I_2})=\beta_1$.
\end{proof}

Now we can come to the proof of Theorem~\ref{thm:reveal}.

\begin{proof}[Proof of Theorem~\ref{thm:reveal}]
 From earlier analysis, we have $$\mu_1\cdot \mu_2\cdot\delta_{I_1,I_2}=i_*(p_1^*(\gamma_1)\cdot p_2^*(\gamma_2)).$$
 Use the adjoint formula between pullback and pushforward, the result in Lemma~\ref{lem:direct_product}, 
 and the property that the pullback is a ring homomorphism. Then consider the isomorphism 
 $D_{I_1,I_2}\cong \mathcal{M}_{I_1\cup\{x\}}\times \mathcal{M}_{I_2\cup\{x\}}$, and use the 
fact on direct product of varieties mentioned earlier, we further get
   \begin{flalign*}
\int_{\mathcal{M}_N}{\mu_1\cdot \mu_2\cdot (\delta_{I_1,I_2})^r} &=  
\int_{\mathcal{M}_N}{(\delta_{I_1,I_2})^{r-1}\cdot i_*(p_1^*(\gamma_1)\cdot p_2^*(\gamma_2)) }&&\\
                      &=\int_{D_{I_1,I_2}}{i^*((\delta_{I_1,I_2})^{r-1})\cdot p_1^*(\gamma_1)\cdot p_2^*(\gamma_2)}&&\\
                       &= \int_{D_{I_1,I_2}}{(p^*_1(\beta_1)+p^*_2(\beta_2))^{r-1}\cdot p_1^*(\gamma_1)\cdot p_2^*(\gamma_2)}&&\\
     = \sum_{k=0}^{r-1}&{r-1\choose k}\cdot {\int_{D_{I_1,I_2}}{p_1^*(\beta_1^k\cdot \gamma_1)\cdot p_2^*(\beta_2^{r-1-k}\cdot \gamma_2)}}&&\\
                      = \sum_{k=0}^{r-1}&{r-1\choose k}\cdot \int_{\mathcal{M}_{I_1\cup\{x\}}}{\beta_1^k\cdot \gamma_1} \cdot 
                     \int_{\mathcal{M}_{I_2\cup\{x\}}}{\beta_2^{r-1-k}\cdot \gamma_2}&&
\end{flalign*}
\textbf{\em Recall that the integral value is defined to be zero if the monomial is not in the Chow group of codimension~$3$ of the ambient Chow ring.}
Therefore, we can already omit those summands that are zero in the above sum. Since $\beta_1\in A^1_{num}(\mathcal{M}_{I_1\cup\{x\}})$
and $\beta_2\in A^1_{num}(\mathcal{M}_{I_2\cup\{x\}})$, we see that we only need to consider the summands such that
$k+s_1=|I_1|+1-3$ and $r-1-k+s_2=|I_2|+1-3$ hold, that is, $k=|I_1|-s_1-2=r+1+s_2-|I_2|$. 
Hence there is only one summand left, we obtain the following
formula:
   \begin{flalign*}
&\int_{\mathcal{M}_N}{\mu_1\cdot \mu_2\cdot (\delta_{I_1,I_2})^r}&&\\
= {r-1\choose |I_1|-s_1-2}&\cdot
\int_{\mathcal{M}_{I_1\cup \{x\}}}{\beta_1^{|I_1|-s_1-2}\cdot \gamma_1}
\cdot \int_{\mathcal{M}_{I_2\cup \{x\}}}{\beta_2^{|I_2|-s_2-2}\cdot \gamma_2}.&&
\end{flalign*}
As a special case, let $I_2=\{a,b\}$, then we have $\mu_1=\nu_1$, $\mu_2=1$, $s_2=0$.
In this case, we see that $|I_1|+2-3=s_1+r$, that is, $r=|I_1|-s_1-1$. The formula becomes
$$\int_{\mathcal{M}_{I_1\cup\{a,b\}}}{\nu_1\cdot (\delta_{I_1,\{a,b\}})^{|I_1|-s_1-1}}
=\int_{\mathcal{M}_{I_1\cup\{x\}}}{\beta_1^{|I_1|-s_1-2}}\cdot \gamma_1.$$
Analogously, when $I_1=\{a,b\}$, we get the following equation:
$$\int_{\mathcal{M}_{I_2\cup\{a,b\}}}{\nu_2\cdot (\delta_{I_2,\{a,b\}})^{|I_2|-s_2-1}}
=\int_{\mathcal{M}_{I_2\cup\{x\}}}{\beta_2^{|I_2|-s_2-2}}\cdot \gamma_2.$$
The statement follows.
\end{proof}

\section{Pseudo code}
 
 In this section, we provide the missing pseudo codes for the forest algorithm. 
 Our algorithm is also implemented in Python \footnote{\url{https://github.com/muronghezi/integral-chow-ring-monomial}}, where the input is the monomial 
 in the Chow ring~$A^{\bullet}(\overline{\mathcal{M}}_{0,n})$, the output is its integral value. 
 Example~\ref{eg:forest} is also illustrated there.

\begin{algorithm}

\caption{sign}
\thispagestyle{empty}
\SetKwInOut{Input}{input}
\SetKwInOut{Output}{output}

\Input{a loaded tree~$LT$ with~$n$ labels and~$(n-3)$ fringes}
\Output{$1$ or~$-1$, which is the sign of tree value of~$LT$}
\texttt{\\}
$V \gets $ list of vertices  of~$LT$\;

$h\gets$ the labeling function of~$LT$\;
$w\gets$ empty function for the weight function of vertices of~$LT$\;
$a\gets 0$\;

\For{each element $v\in V$}
    {$w(v):=\deg(v)+|h(v)|-3$\;
      $a:=a+w(v)$}
\Return{$(-1)^a$}      

\end{algorithm}

\begin{algorithm}

\caption{redundancy$\_$forest}
\thispagestyle{empty}
\SetKwInOut{Input}{input}
\SetKwInOut{Output}{output}

\Input{a loaded tree~$LT$ with~$n$ labels and~$(n-3)$ fringes}
\Output{the corresponding redundancy forest of~$LT$}
\texttt{\\}
$V \gets $ list of vertices  of~$LT$\;
$E \gets $ list of edges of~$LT$\;
$h\gets$ the labeling function of~$LT$\;
forest$\_$node$_1 \gets$ empty list for the vertices of the redundancy forest which
come from the vertices of~$LT$\;
forest$\_$node$_2 \gets$ empty list for the vertices of the redundancy forest which
come from the edges of~$LT$\;
forest$\_$node$\gets$ empty list for the vertices of the redundancy forest\;
forest$\_$edge$\gets$ empty list for the edges of the redundancy forest\;
$w\gets$ empty function for the weight of vertices and edges of~$LT$\;
$m\gets$ multiplicity function for edges of~$LT$\;
\For{each element $v\in V$}
     {Assign weight to~$v$ by $w(v):=\deg(v)+|h(v)|-3$\;
     \If{$w(v)\neq 0$}
         {Append~$v$ to forest$\_$node$_1$} }
\For{each element $e\in E$}
     {Assign weight to~$e$ by $w(e):= m(e)-1$\;
     \If{$w(e)\neq 0$}
        {Append~$e$ to forest$\_$node$_2$} }
\For{each element $x\in$ forest$\_$node$_1$}
    {\For{each element $y \in$ forest$\_$node$_2$}
    {\If{$x$ is an element of~$y$}
    {Append $\{x,y\}$ to forest$\_$edge\;
    Append~$x$ to forest$\_$node\;
    Append~$y$ to forest$\_$node}}}
\Return{(forest$\_$node, forest$\_$edge,~$w$)}
    
\end{algorithm}

\begin{algorithm}\label{absolute}

\caption{absolute$\_$value}
\thispagestyle{empty}
\SetKwInOut{Input}{input}
\SetKwInOut{Output}{output}

\Input{a redundancy forest~$RF$}
\Output{ a natural number, which is the absolute value of~$RF$}
\texttt{\\}

$RF_0\gets$ a list of degree-zero vertices of~$RF$\;
$RF_1\gets$ a list of degree-one vertices of~$RF$\;
$w \gets $ the weight function of vertices of~$RF$\;

\If{$RF_0$ is not an empty list }
   {\Return $0$}
 \ElseIf{$RF_1$ is an empty list}
         {\Return $1$}
 \Else 
   {\For{any element~$x\in RF_1$}
       {$x_1$ is the unique parent of~$x$\;
         $w_0:= w(x)$\;
         $w_1:= w(x_1)$ \;
         \If{$w_0=w_1$}
             {Delete vertex~$x$ from~$RF$\;
              Delete vertex~$x_1$ from~$RF$\;
              \Return absolute$\_$value($RF$)}
          \ElseIf{$w_0 > w_1$}
                  {\Return 0}
           \Else
                {$w(x_1):=w_1 - w_0$\;
                Delete vertex~$x$ from~$RF$\;
                \Return binomial$\_$coefficient($w_1,w_0$)$\cdot$ absolute$\_$value($RF$)}}}
\end{algorithm}

\begin{algorithm}
\thispagestyle{empty}
\caption{tree$\_$value}
\SetKwInOut{Input}{input}
\SetKwInOut{Output}{output}

\Input{a loaded tree~$LT$ with~$n$ labels and~$(n-3)$ fringes}
\Output{an integer, which is the tree value of~$LT$}
\texttt{\\}
$RF \gets$ the redundancy forest of~$LT$\;
$s \gets$ the sign of tree value of~$LT$\;
$ab\gets$ the absolute value of~$RF$\;
\Return{$s\cdot ab$}

\end{algorithm}

\newpage

\section{Proof of the identity}\label{sec:proof_identity}

Continuing with the notations before, we give a combinatorics proof for Theorem~\ref{thm:identity} 
in the sequel
--- which then also leads to the correctness of both Theorem~\ref{thm:identity} and Theorem~\ref{thm:sun_like}.

We slightly modify the notations, so that they serve well for our proof ---
namely we add an index~$r$ for many of them, indicating that we are considering~$r$ many sums for the multinomial 
coefficient. We will see later on that this index is helpful. We denote:
\begin{itemize}
\item $m_1,m_2,\ldots, m_r$: $r$-many positive-integer parameters.
\item $s_r:=\sum_{i=1}^r{m_i}$.
\item $X_r:=\{x_1,x_2,\ldots,x_r\}$: a set of~$r$-many indeterminates. This set is introduced so that we can 
consider all the bipartitions of the values $\{m_1-1,m_2,\ldots,m_r\}$. 
In this way, we are able to formally go through
all combinations.
\item $T_r:=\{B\mid B\subset X_r, x_1\in B\}$. The elements in~$T_r$ indicates one part of the bipartition and we always
put~$x_1$ in it, so as to avoid repetition.
\item $\mathcal{B}_r:=\{(B_1,B_2)\mid B_1\in T_r, B_2 = X_r\setminus B_1\}$. This set is exactly the collection of all 
the bipartition of~$X_r$. 
\item $g_r: X_r\to \{m_1-1,m_2,\ldots, m_r\}$, $x_1\mapsto m_1-1$, $x_i\mapsto m_i$ for $i\neq 1$. 
\item $S(B):=\sum_{x\in B}{g_r(x)}$, for $B\subset X_r$. 
\item $${S(B) \choose B}:= \frac{S(B) !}{\prod_{x\in B}{(g_r(x)!)}},$$ for $B\subset X_r$. 
%Note that
%this notation is just a generalized definition of multinomial coefficient.
\end{itemize}
Besides these notations from earlier, we still need several more.
\begin{itemize}
 \item Define $$S_r:=\{(P_1, P_2, \ldots, P_r)\mid  \cup_{i=1}^r{P_i} = \{1,2,\ldots,s_r\},\;  |P_i|=m_i \}.$$
 With this set, we collect all partitions of the set $\{1,2,\ldots,s_r\}$ into~$r$ parts such that the~$i$-th 
 part has cardinality~$m_i$.
 \item Let $L_r:=\{2,3,\ldots,r\}$. These elements are {\bf special} elements in $\{1,2,\ldots,s_r\}$. 
 Later on we will see why or how they are special, in the definition of the function~$\varphi_r$.
 The next two notations are also there to serve the definition of the function~$\varphi_r$.
 \item For $A\subset \{1,2,\ldots,r\}$ and a partition $(P_1,\ldots,P_r)$, define $P_A:=\cup_{i\in A}{P_i}$. The set~$P_A$ is the union of the parts which have
 index in~$A$.
 \item For $A\subset \{1,2,\ldots,r\}$, define $X_A:=\{x_i \mid i\in A\}$. The set~$X_A$ collect the indeterminates that have
 index in~$A$.
\end{itemize}

Let us see an example, so that we do not get lost among the ocean of notations.
\begin{example}
 Given $r=3$, the following facts are already clear:
\begin{itemize} 
 \item $X_3=\{x_1,x_2,x_3\}.$
 \item $T_3=\{\{x_1\},\{x_1,x_2\},\{x_1,x_3\}, \{x_1,x_2,x_3\}\}.$ This is the collection of one part of the bipartition
 of~$X_3$ that contains~$x_1$.
 \item $\mathcal{B}_3=\{(\{x_1\},\{x_2,x_3\}), 
 (\{x_1,x_2\},\{x_3\}), (\{x_1,x_3\},\{x_2\}),
 (\{x_1,x_2,$
 
 $x_3\},\emptyset)\}.$ This is the collection of all bipartitions
 of~$X_3$.
 \item $L_3=\{2,3\}$. The elements~$2$ and~$3$ are special.
 \item Take $A=\{1,2\}\subset \{1,2,3\}$ for instance, then $X_A=\{x_1,x_2\}$ --- the collection of indeterminate
 with index in~$A$.
\end{itemize} 
However, in order to figure out those remaining notations, we should know the values of~$m_i$ for $1\leq i\leq r$.
Let $m_1=2$, $m_2=2$ and $m_3=1$ for instance, then we also obtain the following facts:
\begin{itemize}
 \item $s_3=\sum_{i=1}^3{m_i}=2+2+1=5.$ 
 \item $g_3:X_3\to \{1,2\}$ is defined as $g_3(x_1)=m_1-1=1$, $g_3(x_2)=m_2= 2$ and $g_3(x_3)=m_3= 1$. 
 \item Take $B=\{x_2,x_3\}\subset X_3$ for instance, then $$S(B)=g_3(x_2)+g_3(x_3)=m_2+m_3=3.$$ 
 Note that in this case~$S(B)$ is just the sum of~$m_2$ and~$m_3$, since $x_1\notin B$.
 \item Take $B=\{x_2,x_3\}\subset X_3$ for instance, then 
 $${S(B)\choose B}=\frac{S(B) !}{\prod_{x\in B}{(g_3(x)!)}}=\frac{3!}{g_3(x_2)\cdot g_3(x_3)}=\frac{6}{m_2\cdot m_3}=\frac{6}{2\cdot 1}
 =3.$$ 
 %Note that in this case, ${S(B)\choose B}$ is just the normal multinomial coefficient ${m_2+m_3 \choose m_2,m_3}$, since
 %$x_1\notin B$.
 \item $S_3$ is the set of all partitions $(P_1,P_2,P_3)$ of the set $\{1,2,3,4,5\}$ into three parts $P_1,P_2,P_3$
such that $|P_1|=m_1=2$, $|P_2|=m_2=2$ and $|P_3|=m_3=~1$.
\item Take $A=\{1,2\}\subset \{1,2,3\}$ for instance, then $P_A=P_1\cup P_2$ for some $(P_1,P_2,P_3)\in S_3$.
 \end{itemize}
 We leave it to the readers to check that the identity in Theorem~\ref{thm:identity} holds in this example.
 \end{example}

Now we define a function $\varphi_r: S_r\to T_r$, $(P_1, \ldots, P_r)\mapsto B$ by
Algorithm~\ref{alg:propoganda}. We will prove that it is indeed an algorithm later on.

\begin{algorithm}
\caption{the function~$\varphi_r$}\label{alg:propoganda}
\SetKwInOut{Input}{input}
\SetKwInOut{Output}{output}

\Input{$(P_1,\ldots, P_r)\in S_r$.}
\Output{$B\in T_r$.}
\texttt{\\}

$B \gets$ $\{x_1\}$\;
$A \gets$ $L_r\cap P_1$\;

\While{$A\neq \emptyset$}
        {$B = B\cup X_A$\;
         $A:= L_r\cap P_A$.}
\Return $B$
\end{algorithm}

For a better understanding, let us see how is this function defined in our running example.
\begin{example}
 $\varphi_3: S_3\to T_3$, $(P_1, P_2, P_3)\mapsto B\in T_3$. 
 Let us go through Algorithm~\ref{alg:propoganda} with the input $(P_1,P_2,P_3)=(\{1,3\},\{4,5\},\{2\})$. 
 First, we see that $L_3=\{2,3\}$.
 \begin{enumerate}
  \item Input: $(P_1,P_2,P_3)=(\{1,3\},\{4,5\},\{2\})$.
  \item Initial values: $B=\{x_1\}$, $A=\{2,3\}\cap\{1,3\}=\{3\}$.
  \item First loop: 
  Since $A=\{3\}\neq \emptyset$, we have $B=\{x_1\}\cup X_{\{3\}}=\{x_1\}\cup \{x_3\}=\{x_1,x_3\}$, 
  and then $A=\{2,3\}\cap \{2\}=\{2\}$.
  \item Second loop: 
  Since $A=\{2\}\neq \emptyset$, we have $B=\{x_1,x_3\}\cup X_{\{2\}}= \{x_1, x_3\} \cup\{x_2\}=\{x_1,x_2,x_3\}$,
  and then $A=\{2,3\}\cap \{4,5\}=\emptyset$.
  \item Since $A=\emptyset$,
  return $B=\{x_1,x_2,x_3\}$.
  \item Output: $B=\{x_1,x_2,x_3\}$.
 \end{enumerate}

\end{example}

For the convenience of later analysis, we introduce an extra subscript for the variable~$A$ in Algorithm~\ref{alg:propoganda}. 
Let $A_0:=\{1\}$ and define~$A_i$ to be~$A$ in the~$i$-th loop.

\begin{lemma}
 In the above defined process (Algorithm~\ref{alg:propoganda}), $A_i\cap A_j = \emptyset$ for all $i\neq j$. 
\end{lemma}
\begin{proof}
When $i=0,\; j\neq 0$, we have $A_0\cap A_j=\emptyset$ since $A_j\subset L_r$,
$A_0=\{1\}$ and $1\notin L_r$.
Suppose $A_i\cap A_j\neq \emptyset$ when $i,j > 0$ and assume
w.l.o.g. that $i<j$. Since $A_i:=L_r\cap P_{A_{i-1}}$, we obtain that $L_r\cap P_{A_{i-1}}\cap P_{A_{j-1}}\neq \emptyset$,
hence $A_{i-1}\cap A_{j-1}\neq \emptyset$. Repeating the 
similar process, after finite steps we reach a situation where 
$A_0\cap A_{j-i}\neq \emptyset$. This is a contradiction.
\end{proof}

\begin{proposition}\label{prop:empty}
 Algorithm~\ref{alg:propoganda} is an algorithm.
\end{proposition}
\begin{proof}
 Since $|L_r| < \infty$, $A_i\subset L_r$ for all $1\leq i\leq r$, and $A_i\cap A_j=\emptyset$ for all $i\neq j$, 
 there must exist $i\in\mathbb{N}^+$ such that $A_i=\emptyset$. Therefore,
 this process terminates. And it is well-defined --- once the input is given, the output is uniquely determined 
 via the process
 and clearly $B\in T_r$ ---
 thence it is indeed an algorithm. 
\end{proof}

\begin{proposition}
 The function $\varphi_r:S_r\to T_r$ is a surjection.
\end{proposition}
\begin{proof}
For any $B\in T_r$, define $(P_1,\ldots,P_r)$ as follows:
\begin{itemize}
\item Define $Q:=\{i\mid i\in L_r,\; x_i\in B\}$. Since $|Q|<\infty$, we can list its elements as: $q_1,\ldots,q_t$, assuming that 
$|Q|=t$. 
\item Let $P_1:=\{q_1\}$, $P_{q_j}=\{q_{j+1}\}$ for $1\leq j\leq t-1$. 
\item Now, we already defined all~$P_i$ when $i\in L_r$ and $x_i\in B$ except for~$P_{q_t}$, 
and we have defined~$P_1$ as well.
\item Let $P_i=\{i\}$ if $i\in L_r$ and $x_i\notin B$. Let 
$$P_{q_t}:=\{1,\ldots,s_r\}\setminus (\cup_{i\neq q_t}{P_i}).$$  
\end{itemize}
It is clear that $1\in P_{q_t}$, hence all parts defined above are non-empty. One can check that they indeed form
a partition of $\{1,\ldots,s_r\}$ into~$r$ parts. Hence the input is well-defined. Furthermore, one can check that 
$\varphi_r(P_1,\ldots,P_r)=B$, since going through the process in Algorithm~\ref{alg:propoganda} outputs
$\{x_{q_1},\ldots, x_{q_t}, x_1\}=B$.
\end{proof}
We see that~$\varphi_r$ is a well-defined surjective function. 
 Then, 
 $$\bigcup_{B\in T_r}{\varphi_r^{-1}(B)}=S_r$$ 
 and 
 $$|S_r|=\sum_{B\in T_r}{|\varphi_r^{-1}(B)|}=\sum_{(B,X\setminus B)\in \mathcal{B}_r}{|\varphi_r^{-1}(B)|}.$$ 
 In order to prove the identity, we only need to show one thing and it is formulated as the following lemma.
\begin{lemma}\label{lem:no_sum}
 For any $B_1\in T_r$, define $B_2:= X_r\setminus B_1$, then 
 $$|\{x\in S_r \mid \varphi_r(x)=B_1\}|={{s_r-r+1} \choose {S(B_2)-|B_2|}}{S(B_1) \choose B_1}{S(B_2)
 \choose B_2}.$$
\end{lemma}
 In order to prove the above lemma, we need to introduce the following proposition.
 \begin{proposition}\label{prop:observation}
  If $\varphi_r(P_1,\ldots, P_r)= B_1$ for some $(P_1,\ldots,P_r)\in S_r$ and
  $B_1\in T_r$; denote $B_2:=X_r\setminus B_1$.
  Then $$P_{F_{B_1}} \cap L_r = F_{B_1}\setminus \{1\},$$ where $F_B:=\{i\mid x_i\in B\}$. 
  Consequently, we have $P_{F_{B_2}} \cap L_r = F_{B_2}$ and $|P_{F_{B_2}}\cap L_r|=|B_2|$.
 \end{proposition}
\begin{proof}
 From Proposition~\ref{prop:empty} we know that there exists $k\in \mathbb{N}^+$ 
 such that $A_k=~\emptyset$. Let $t:=k-1$, we claim that $\cup_{i=0}^t{A_i}=F_{B_1}$.
 It is clear from Algorithm~\ref{alg:propoganda} that $B_1=\cup_{i=0}^t{X_{A_i}}$,
 which is equivalent to $F_{B_1}=\cup_{i=0}^t{A_i}$. Hence we know that 
 $F_{B_1}\setminus \{1\}=\cup_{i=1}^t{A_i}$. We only need to show 
 $P_{F_{B_1}} \cap L_r=\cup_{i=1}^t{A_i}$.
 
 For any $m\in \cup_{i=1}^t{A_i}$, there exists $1\leq j\leq t$ such that $m\in A_j$.
 From Algorithm~\ref{alg:propoganda} we know that $A_j:=L_r\cap P_{A_{j-1}}$. Hence
 $m\in L_r$. And $X_{A_{j-1}}\subset B_1$ is equivalent to $A_{j-1}\subset F_{B_1}$, which implies
  $P_{A_{j-1}}\subset P_{F_{B_1}}$. Hence $m\in P_{A_{j-1}}\subset P_{F_{B_1}}$. We obtain that 
  $m\in P_{F_{B_1}} \cap L_r$. For any $m\in P_{F_{B_1}} \cap L_r$, equivalently we have 
  $m\in \cup_{i=0}^t{A_i}\cap L_r=\cup_{i=1}^t{A_i}\cap L_r$. Since $A_i\subset L_r$ for 
  any $i\neq 0$, we obtain that $m\in \cup_{i=1}^t{A_i}$.
  
  Since $P_{F_{B_1}}\cup P_{F_{B_2}}=\{1,\ldots,s_r\}$ and $L_r\subset \{1,\ldots,s_r\}$, 
  we know that $(P_{F_{B_1}}\cap L_r) \cup (P_{F_{B_1}}\cap L_r)=L_r$. Therefore
  $$P_{F_{B_2}}\cap L_r=L_r\setminus (P_{F_{B_1}}\cap L_r)=L_r\setminus (F_{B_1}\setminus \{1\})=F_{B_2}.$$
 Then, $|P_{F_{B_2}}\cap L_r|=|B_2|$ follows naturally, since $|F_{B_2}|=|B_2|$.
\end{proof}

 Let $B_1\in T_r$ and $P\in\varphi_r^{-1}(B_1)$, by Proposition~\ref{prop:observation}, we know that $P_{F_{B_1}}\cap L_r=F_{B_1}\setminus~\{1\}$,
 hence the special elements in~$P_{F_{B_1}}$ are fixed (they are just~$F_{B_1}\setminus~\{1\}$, otherwise the image of~$P$ would not be~$B_1$) --- so do
 the special elements in~$P_{F_{B_2}}$ 
 since $P_{F_{B_2}}\cap L_r=F_{B_2}$, where $B_2:=X_r\setminus B_1$. Let $K_r:=\{1,\ldots,s_r\}$. Also, it is evident that 
 $P_{F_{B_1}}\cup P_{F_{B_2}}=K_r$.  Inspired by Proposition~\ref{prop:observation}, in order to figure out what configurations in~$S_r$ 
 are mapped to~$B_1$ by the function~$\varphi_r$, we view the problem in the following way.
 
 Given $B_1\in T_r$, we know that the special elements in~$P_{F_{B_1}}$ are fixed, where $P\in\varphi_r^{-1}(B_1)$. 
 We only need to choose a proper amount of elements 
 in~$K_r\setminus L_r$ and put them into~$P_{F_{B_1}}$. We call the elements in~$K_r\setminus L_r$ {\bf non-special}.
  We need to choose $$|P_{F_{B_1}}|-(|B_1|-1)= S(B_1)-|B_1|+1$$ many elements from 
 $$|K_r\setminus L_r|=s_r-|L_r|=s_r-r+1$$ many 
 elements, and put them in the group of~$P_{F_{B_1}}$. There are $${s_r-r+1 \choose S(B_1)-|B_1|+1}$$ many ways to do so. Since 
 $(S(B_1)-|B_1|+1)+(S(B_2)-|B_2|) =(S(B_1)+S(B_2))-(|B_1|+|B_2|)+1 =s_r-r+1$,
 we can also say that there are 
 $${s_r-r+1 \choose S(B_2)-|B_2|}$$ many ways to arrange the non-special elements.
 Considering the definition of~$\varphi_r$,
 we see that no matter how we arrange the elements in~$P_{F_{B_2}}$, the image of~$P$ under the function~$\varphi_r$ is not influenced. 
  Therefore, there are~${S(B_2) \choose B_2}$ many configurations for the elements in~$P_{F_{B_2}}$.
  As for the arrangements in~$P_{F_{B_1}}$, they need to obey certain rules in order to guarantee that the value of~$\varphi_r$ is~$B_1$. 
  
  From the analysis above, the first and third factors on the right hand side of the equation in Lemma~\ref{lem:no_sum} are 
  both explained well in a combinatorics way. In order to prove this lemma, we only need to prove that 
  given~$B_1\in T_r$, the number of configurations for the elements in~$P_{F_{B_1}}$ is 
  exactly~${S(B_1)\choose B_1}$ --- which would then conclude the proof of Lemma~\ref{lem:no_sum}. 
  Recalling the definition of~${S(B_1)\choose B_1}$, one can see that the remaining work for the 
  proof of Lemma~\ref{lem:no_sum} is equivalent to
  proving the following proposition.

 \begin{proposition}\label{prop:induction}
 Recall that $s_k:=\sum_{i=1}^k{m_i}$ and that $X_k:=\{x_1,\ldots,x_k\}$. Then we have
   $$f_k(m_1,m_2,\ldots, m_k)={s_k-1 \choose {m_1-1,m_2,\ldots,m_k}},\, k\in \mathbb{N^+},
   m_i\in\mathbb{N^+},$$ where  $f_k:(\mathbb{N}^+)^k\to \mathbb{N}$, 
 $(m_1,m_2,\ldots,m_k)\mapsto 
 |\{(P_1,P_2,\ldots,P_k)\in S_k\mid |P_i|=m_i, \varphi_k(P_1,P_2,\ldots,P_k)=X_k\}|.$
 \end{proposition}
 To explain the above defined function~$f_k$ in another way, we have 
 $$f_k(m_1,\ldots,m_k)=\#\{ \varphi^{-1}_k(X_k)\mid |P_i|=m_i\};$$ it is the cardinality of 
 the fiber of~$\varphi_k$ given that $|P_i|=m_i$. 
 
 In order to prove Theorem~\ref{thm:identity}, we only need to 
 prove Proposition~\ref{prop:induction}.
 Before we approach the proof, we need to introduce a known identity on
 multinomial coefficients. We also provide a proof for it since we did not 
 find a good reference.
 \begin{lemma}\label{lem:multinomial}
  For all $s,k,m_1,\ldots,m_k\in \mathbb{N}$ with $m_1+\cdots +m_k=s$, $s\geq 1$
  and $k\geq 2$, we have
  $${s \choose {m_1,\ldots,m_k}}=\sum_{i=1}^k{s-1 \choose {m_1,\ldots,m_i-1,\ldots,m_k}}.$$
 \end{lemma}
 \begin{proof}
  We want to partition~$s$ many apples into~$k$ piles, each has cardinality~$m_i$. We 
  can pick one apple and trace the position of it: this apple can be in the~$i$-th pile, then 
  we need to partition the remaining~$s-1$ many apples into~$k$ piles such that
  all other piles have cardinality~$m_j$ and the~$i$-th pile has cardinality~$m_i-1$. This concludes
  the above identity.
 \end{proof}

\begin{proof}[Proof of Proposition~\ref{prop:induction}]
 Prove by induction. When $k=1$, $L_1=\emptyset$, for any $m_1\in \mathbb{N^+}$, we have
 $$|\{(P_1)\in S_1\mid \varphi_1(P_1)=\{x_1\}\}|=1={s_1-1 \choose m_1-1}$$
 since $s_1=m_1$ in this case.
 Assume that the proposition holds whenever the number of parameters --- here parameters refer to~$m_i$ ($1\leq i\leq k$) --- 
 is less or equal to~$k-1$, where $k\geq 2$.
 
 When the number of parameters is~$k$, we start the inner induction on~$s_k$. 
 Obviously 
 $s_k\geq k$. When $s_k=k$, we know that $$m_1=m_2=\cdots=m_k=1.$$
 In the configurations that
 is mapped to~$X_k$ under~$\varphi_k$, we can choose any element in~$L_k$ for~$P_1$, 
 say~$i_1$; there are $|L_k|=k-1$ many possibilities. Then we can choose an element in~$L_k\setminus \{i_1\}$ 
 for~$P_{i_1}$, and so on. Until we choose the element $i_{k-1}\in L_k$ for~$P_{i_{k-2}}$. Then we already 
 arranged~$k-2$ many parts, then the only remaining part~$P_{i_{k-1}}$ can only be~$\{1\}$. 
 In total there are~$(k-1)!$ many configurations. Hence we have 
 $$f_k(m_1,m_2,\ldots,m_k)=(k-1)!={k-1 \choose {1,\ldots,1}}
 ={k-1 \choose {0,1,\ldots,1}},$$
 which equals to $${s_k-1 \choose {m_1-1,m_2,\ldots,m_k}}.$$

 Assume that the proposition holds whenever the sum of these parameters is less or equal to~$s_k-1$, 
 where we can assume $s_k-1\geq k$, i.e., $s_k\geq k+1$. When the sum of these parameters equals~$s_k$,
 recall that~$f_k(m_1,m_2,\ldots,m_k)$ is the number of configurations~$P$ that
 are mapped to~$X_k$ under~$\varphi_k$, given that $|P_i|=m_i$. We focus on the position of the element~$1$ among the piles~$P_i$, for $1\leq i\leq k$.
 Since~$1\notin L_k$, 
 it does not influence the value of~$\varphi_k$ on any configuration. So in the case when 
 $m_i\geq 2$ for all $1\leq i\leq k$, there are~$k$
 cases for the distribution of~$1$, i.e., it can belong to any part~$P_i$ for $1\leq i\leq k$.
 Hence in this case we obtain the following identity: 
  \begin{flalign*}
 f_k(m_1,m_2,\ldots,m_k) 
= \; &f_k(m_1-1,m_2,\ldots,m_k)+f_k(m_1,m_2-1,\ldots,m_k)\\\nonumber
 &+\cdots+f_k(m_1,m_2,\ldots,m_k-1).&&
 \end{flalign*}
 Now we can apply the induction hypothesis on the sum of the parameters. Then we obtain the 
 following equation: 
 \begin{flalign*}
 &f_k(m_1,m_2,\ldots,m_k) = \; {s_k-2 \choose {m_1-2,m_2,\ldots,m_k}}\\\nonumber
 &+ {s_k-2 \choose m_1-1,m_2-1,\ldots,m_k}+\cdots+{s_k-2 \choose m_1-1,m_2,\ldots,m_k-1}.&&
 \end{flalign*}
 Then by Lemma~\ref{lem:multinomial}, we obtain 
 $$f_k(m_1,m_2,\ldots,m_k)={s_k-1 \choose m_1-1,m_2,\ldots,m_k}.$$
 If $m_i=1$ for some $i\neq 1$. When we put~$1$ into~$P_i$, the problem can be reduced to counting the
 number of corresponding configurations of~$P_j$ for $j\neq i$, since $1\notin L_k$. 
 Therefore, in this case, considering the distribution of the element~$1$ gives us the following identity:
 \begin{flalign*}
 f_k(m_1,\ldots,m_k)=\; &f_k(m_1-1,\ldots,m_k)+\\\nonumber
 &\cdots+f_k(m_1,\ldots, m_{i-1}-1,m_i,\ldots,m_k)\\\nonumber
 &+f_{k-1}(m_1,\ldots,m_{i-1},m_{i+1},\ldots,m_k)\\\nonumber
 &+f_k(m_1,\ldots,m_i,m_{i+1}-1,\ldots,m_k)+\\\nonumber
 & \cdots+f_k(m_1,\ldots,m_{i-1},m_i,m_{i+1},\ldots, m_k-1).&&
 \end{flalign*}
 By induction hypothesis on~$k$ we obtain 
 \begin{flalign*}
 f_{k-1}&(m_1,\ldots,m_{i-1},m_{i+1},\ldots,m_k) \\\nonumber
 &=\; {(s_k-m_i)-1 \choose {m_1-1,\ldots,m_{i-1},m_{i+1},\ldots,m_k}}\\\nonumber
 &=\; {(s_k-1)-1 \choose {m_1-1,\ldots,m_{i-1},m_{i+1},\ldots,m_k}}\\\nonumber
 &=\; {s_k-2 \choose {m_1-1,\ldots,m_{i-1},0,m_{i+1},\ldots,m_k}}\\\nonumber
 &= {s_k-2 \choose {m_1-1,\ldots,m_{i-1},m_i-1,m_{i+1},\ldots,m_k}}.&&
 \end{flalign*}
 Substituting back this term, we get the same recurrence for~$f_k(m_1,\ldots,m_k)$ as in the case where $m_i\geq 2$ for $1\leq i\leq k$. 
 By Lemma~\ref{lem:multinomial}
 we as well obtain $$f_k(m_1,m_2,\ldots,m_k)={s_k-1 \choose m_1-1,m_2,\ldots,m_k}.$$
 With the same idea, it is not 
 hard to prove that the statement holds however many parameters except for~$m_1$ equals one.

 If $m_1=1$, from the definition of the function~$f_k$ and~$\varphi_k$, we know that $1\notin P_1$. 
 Hence considering the distribution of the element~$1$, the recurrence formula becomes 
 \begin{flalign*}
 f_r(m_1,m_2,\ldots,m_r) =&f_r(m_1,m_2-1,\ldots,m_r)+\\\nonumber
&\cdots+f_r(m_1,m_2,\ldots,m_r-1).&&
\end{flalign*}

 Then by induction hypothesis on the sum of the parameters, we obtain
 \begin{flalign*}
 &f_k(m_1,m_2,\ldots,m_k)\\\nonumber
 &={s_k-2 \choose {m_1-1,m_2-1,\ldots,m_k}}+\cdots+{s_k-2 \choose {m_1-1,m_2,\ldots,m_k-1}}\\\nonumber
 &={s_k-2 \choose {0,m_2-1,\ldots,m_k}}+\cdots +{s_k-2 \choose {0,m_2,\ldots,m_k-1}}\\\nonumber
 &={s_k-2 \choose {m_2-1,\ldots,m_k}}+\cdots +{s_k-2 \choose {m_2,\ldots,m_k-1}}.&&
 \end{flalign*}
 Now we can apply Lemma~\ref{lem:multinomial} and then obtain
 \begin{flalign*}
 f_k(m_1,m_2,\ldots,m_k)&={s_k-1 \choose {m_2,\ldots,m_k}}\\\nonumber
 &={s_k-1 \choose {0,m_2,\ldots,m_k}}\\\nonumber
 &={s_k-1 \choose {m_1-1,m_2,\ldots,m_k}}.&&
 \end{flalign*}
 By induction, the proposition holds.
 \end{proof}
 With this, we finished proving Theorem~\ref{thm:identity}.

 \newpage

\section*{Acknowledgements}
The research was funded by the Austrian Science Fund (FWF): W1214-N15, project DK9.

I am truly grateful to Josef Schicho for the references in Section~\ref{subsec:correctness_forest_algorithm},
providing me with Theorem~\ref{thm:reveal}
and the proof of it, the proof of Theorem~\ref{thm:edge-cutting} in algebraic geometric way, 
  which to a great extend realized the correctness proof, hence further 
drew a perfect full stop on the forest algorithm. 
I thank Nicolas Allen Smoot for 
helping me formulate in a good way the existence of star-cuts.

I genuinely thank
Cristian-Silviu Radu for helping me formulate the identity and as well its proof in a mathematically proper way.
 Also, I am sincerely thankful to Dongsheng Wu for the instructive discussions on the proof of the identity,
specifically for providing me with the main idea of the proof of Proposition~\ref{prop:induction}.

I authentically thank Matteo Gallet for the valuable suggestions on how to formulate a good introduction, 
for providing me with the references 
in Section~\ref{subsec:correctness_forest_algorithm}, and for answering 
to me some general questions in intersection theory, which contributes a lot to the background 
part of the introduction.

\end{document}